



\input amssym.def


\font \titel=cmr10 scaled\magstep2



\font \kapitel=cmbx10 scaled\magstephalf

\font \kapitel=cmbx10 scaled\magstephalf

\font \tenrm=cmr10
\font \mathbf=cmmib10
\font \kpt=cmcsc10
\font \mathit=cmmi10

\font \teneusm=eusm10
\font \blackboard=dsrom10
\font \boldsym=cmbsy10

\font \ninemathbf=cmmib9
\font \ninemsa=msam9
\font \ninemsb=msbm9
\font \ninecyr=wncyr9
\font \ninerm=cmr9
\font \ninei=cmmi9
\font \ninesy=cmsy9
\font \nineit=cmti9
\font \ninebf=cmbx9
\font \ninekpt=cmcsc9
\font \nineeusm=eusm9
\font \nineeufm=eufm9
\font \ninecyrit=wncyi9


\font \sevenrm=cmr7
\font \sevenit=cmti7
\font \seveneusm=eusm7

\font \sixrm=cmr6
\font \sixi=cmmi6
\font \sixsy=cmsy6

\font \fiveeusm=eusm5

\scriptfont4=\sevenit

\newfam\eusmfam
\textfont\eusmfam=\teneusm
\scriptfont\eusmfam=\seveneusm
\scriptscriptfont\eusmfam=\fiveeusm
\def\eusm#1{{\fam\eusmfam\relax#1}}

\def\eufm#1{{\fam\eufmfam\relax#1}}

\newfam \varmathbffam
\textfont \varmathbffam=\mathbf
\scriptfont \varmathbffam=\mathbf
\scriptscriptfont \varmathbffam=\mathbf

\def\div{\mathord{\rm div}\,}

\def\e{\varepsilon}

\def\la{\lambda}


\def \vvdots{\vbox{\baselineskip 3pt\lineskiplimit 0pt\kern3pt
         \hbox{.}\hbox{.}\hbox{.}}}
\def\ymk{\hbox{\boldsym Y\kern1pt}(\K)}
\def\gmk{\hbox{\boldsym J\kern1pt}(\K)}

\def \Gf #1 #2{\mathop{\mathstrut{\mathit G}}
        \nolimits^{#1}_{#2}}
\def \Rf #1 #2{\mathop{\mathstrut{\mathit R}}
        \nolimits^{#1}_{#2}}
\def \Rfc #1 #2{\mathstrut{\mathop{\mathit H}\limits^c}
\!\mathop{\vphantom{\mathit H}}\nolimits^{#1}_{#2}}
\def \Gfhat #1 #2{\hat{\mathit G}\kern-1.5pt\mathop{\vphantom{\mathit G}}\nolimits^{#1}_{#2}}
\def \Rfhat #1 #2{\hat{\mathit H}\kern-1.5pt\mathop{\vphantom{\mathit H}}\nolimits^{#1}_{#2}}


\catcode`\@=11 \newsymbol
\backepsilon 237F \newsymbol
\complement 107B \newsymbol
\curvearrowright 2379 \newsymbol
\diagdown 231F \newsymbol
\diagup 231E \newsymbol
\dotplus 1275 \newsymbol
\geqslant 133E \newsymbol
\leftarrowtail 131B \newsymbol
\leqslant 1336 \newsymbol
\nexists 2040 \newsymbol
\ngtr 2305 \newsymbol
\nless 2304 \newsymbol
\nparallel 232C \newsymbol
\rightarrowtail 131A \newsymbol
\rightrightarrows 1313 \newsymbol
\sphericalangle 105E \newsymbol
\square 1003 \newsymbol
\subsetneqq 2324 \newsymbol
\supsetneqq 2325 \newsymbol
\upuparrows 1314 \newsymbol
\vartriangleleft 1343 \newsymbol
\vartriangleright 1342


\def \blackbox{\mathord{\hbox{\vrule height 4pt width 4pt depth
        0pt}}}

\def \ccup{\mathop{\lower 2pt\hbox{$\textstyle\cup$}}\nolimits}

\def \cl{\mathord{ {\rm cl}} \,}

\def \cont #1 #2{\mathop{\mathstrut{\mathit C}}\nolimits^{#1}_{#2}}



\def \fforall{(\forall)}
\def \folge #1{\{\,#1\,\}\,}
\def \follows{\Longrightarrow \,}
\def \ge{\geqslant}

\def \home{(\Omega)}
\def \homer{(\Omega, \R {})}

\def \int{\intop \nolimits}
\def \L #1 #2{\mathop{\mathstrut{\mathit L}}
        \nolimits^{#1}_{#2}}
\def \le{\leqslant}
\def \lim{\mathop{\rm lim\vphantom{p}} \limits}

\def \lldots{\,{...}\,}

\def \lq{\char'134}

\def \Max{\mathop{\rm Max\vphantom{p}} \limits}

\def \norm #1{\mathop{\Vert\,#1\,\Vert} \nolimits}
\def \Norm #1{\mathop{\bigl\Vert\,#1\,\bigr\Vert} \nolimits}

\def \p{{\prime}}
\def \pp{{\prime\prime}}

\def \paragraph{\mathhexbox278}

\def \R #1{\mathop{\hbox{\blackboard R}}\nolimits^{#1}}

\def \section{\mathhexbox 278}
\def \sob #1 #2{\mathop{\mathstrut{\mathit W}}
         \nolimits^{#1}_{#2}}
\def \sobnull #1 #2{\mathop{\mathstrut{\mathit W}} \limits^\circ
         \hbox{\kern-2pt}^{#1}_{#2}}

\def \sum{\mathop{{\textstyle\mathchar"1350}}\limits}
\def \sup{\mathop{{\rm sup}} \limits}

\def \zero{\mathord{\eufm o}}
\def \A{\mathord{ {\rm A}} }

\def \C{\mathord{ {\rm C}} }

\def \K{\mathord{ {\rm K}} }

\def \P{\mathord{ {\rm P}} }


\def \V{\mathord{ {\rm V}} }

\def \X{\mathord{ {\rm X}} }
\def \Y{\mathord{ {\rm Y}} }
\def \T{\mathord{ \hbox{\sevenrm T}} }


\newdimen\mathindent \mathindent=20pt
\def\eqno{$\hfill$}

\long\def\leftdisplay#1$${\line{\hskip\mathindent
                                 $\displaystyle#1$\hfil}$$}
\everydisplay{\leftdisplay}
\catcode`\@=11
\def\eqalignno#1{%
    \displ@y \tabskip=0pt
    \advance\displaywidth by -\mathindent
    \vbox{%
        \halign to \displaywidth{%
        \hfil$\displaystyle{##}$\tabskip=0pt
        &$\displaystyle{{}##}$\tabskip=\centering
        &\llap{$##$}\tabskip=0pt\crcr#1\crcr}}}

\catcode`\@=12


\catcode`\@=11
\def\ninepoint{\def\rm{\fam0\ninerm}\def\kpt{\ninekpt}
\textfont0=\ninerm    \scriptfont0=\sixrm
                                      \scriptscriptfont0=\fiverm
\textfont1=\ninei     \scriptfont1=\sixi
                                      \scriptscriptfont1=\fivei
\textfont2=\ninesy    \scriptfont2=\sixsy
                                      \scriptscriptfont2=\fivesy
\textfont6=\ninebf
\textfont \varmathbffam=\ninemathbf
\font\blackboard=dsrom10 scaled900
\def\mathbf{\fam\varmathbffam\ninemathbf}
\textfont\itfam=\nineit \def\it{\fam\itfam\nineit}
\textfont\eufmfam=\nineeufm \def\eufm{\fam\eufmfam\nineeufm}
\textfont\eusmfam=\nineeusm \def\eusm{\fam\eusmfam\nineeusm}
\textfont\msafam=\ninemsa \textfont\msbfam=\ninemsb
\setbox\strutbox=\hbox{\vrule height8pt depth3pt width0pt}%
\normalbaselines\rm}
\def \footnote#1{\edef\@sf{\spacefactor\the\spacefactor}{\sixi
#1}\@sf \insert\footins\bgroup\ninepoint\interlinepenalty100
\let\par=\endgraf \leftskip=0pt \rightskip=0pt
\splittopskip=10pt plus 1pt minus 1pt \floatingpenalty=20000
\smallskip \item{#1} \bgroup \strut \aftergroup \@foot \let
\next}
\skip\footins=12pt plus 2pt minus 4pt \dimen\footins=30pc
\catcode`\@=12


\hyphenation{Rashevsky}

\def\makeheadline{\vbox to 0pt{\vskip-28.5pt%
\baselineskip=7pt\line{\vbox to 8.5pt{}\the\headline}%
\line{\hrulefill}\vss}\nointerlineskip}%
\headline={\hss {\tenrm\folio}}%
\newtoks \rightheadline\rightheadline={\hfil}
\newtoks \leftheadline\leftheadline={\hfil}

\input cyracc.def \def\cyr{\tencyr\cyracc}

\nopagenumbers \parindent=0pt \baselineskip=14.4pt

\overfullrule=4pt
\font\ninett=cmtt9

\def\potin{\mathord{\Phi_i}}
\def\potex{\mathord{\Phi_e}}
\def\pottr{\mathord{\Phi_{\hbox{\sevenit tr}}}}

\def\pottrhat{\mathord{{\hat\Phi}_{\hbox{\sevenit tr}}}}
\def\potexhat{\mathord{{\hat\Phi}_e}}
\def\gat{\mathord{W}}

\def\gathat{\mathord{{\hat W}}}
\def\condin{\mathord{M_i}}
\def\condex{\mathord{M_e}}
\def\excin{\mathord{I_i}}
\def\excinex{\mathord{I_{ie}}}
\def\excex{\mathord{I_e}}
\def\excexhat{\mathord{{\hat I}_e}}
\def\span{\mathord{\rm span}\,}
\def\normal{\mathord{\eufm n}}

\def\ioncurr{\mathord{I_{\hbox{\sevenit ion}}}}
\def\bochnerl #1 #2 #3 {\mathop{\mathstrut{\mathit L}}
        \nolimits^{#1} \bigl[\,{#2}\,,\,{#3}\,\bigr]}
\def\bochnerc #1 #2 #3 {\mathop{\mathstrut{\mathit C}}
        \nolimits^{#1} \bigl[\,{#2}\,,\,{#3}\,\bigr]}
\def\bochnerb #1 #2 #3 {\mathop{\mathstrut{\mathit B}}
        \nolimits^{#1} \bigl[\,{#2}\,,\,{#3}\,\bigr]}
\def\bochnerf #1 #2 #3 #4 {\mathop{\mathstrut{\mathit F}}
        \nolimits^{#1}_{#2} \bigl[\,{#3}\,,\,{#4}\,\bigr]}
\def\bochnersob #1 #2 #3 {\mathop{\mathstrut{\mathit W}}
        \nolimits^{#1} \bigl[\,{#2}\,,\,{#3}\,\bigr]}
\baselineskip=14.4pt\parindent=0pt\frenchspacing
\hyphenation{va-lues}
\hyphenation{models}


{\def\makeheadline{\relax}

\centerline{\titel Optimal control of the bidomain system (IV):}
\medskip
\centerline{\titel corrected proofs of the stability and regularity theorems}
\bigskip
\centerline{\it Karl Kunisch and Marcus Wagner}
\bigskip\bigskip
{\kapitel 1.~Introduction.}
\bigskip
{\bf a) Aim of the paper and main results.}
\medskip
We consider the full bidomain system, which represents a well-accepted description of the electrical activity of the heart, as given through$\,$\footnote{$^{01)}$}{The bidomain model has been considered first in $[\,${\ninekpt Tung 78}$\,]\,$. A detailed introduction may be found e.~g.~in $[\,${\ninekpt Sundnes/ Lines/Cai/Nielsen/Mardal/Tveito 06}$\,]\,$, pp.~21$\,-\,$56.}
$$\eqalignno{
{} & {\partial \pottr\over \partial t} + \ioncurr (\pottr, \gat) - \div \bigl(\, \condin\,\nabla \potin\,\bigr) \,=\, \phantom{-}\excin \quad \hbox{for a.~a.~} \, (x,t) \in \Omega \times [\, 0\,,\,T\,]\,; & (1.1)
\cr
{} & {\partial \pottr\over \partial t} + \ioncurr (\pottr, \gat) + \div \bigl(\,\condex\,\nabla \potex\,\bigr)\,=\, -  \excex \quad \hbox{for a.~a.~} \,(x,t) \in \Omega \times [\, 0\,,\,T\,]\,; & (1.2)
\cr
{} & {\partial \gat \over \partial t} + G(\pottr, \gat) \qquad\qquad\qquad\qquad\quad =\, \phantom{-} 0\quad \hbox{for a.~a.~} \,(x,t) \in \Omega \times [\, 0\,,\,T\,]\,; & (1.3)
\cr
{} & \normal^{\T} \condin\,\nabla \potin \,=\, 0 \quad \hbox{for all }\,(x,t) \in \partial \Omega \times [\,0\,,\, T\,]\,; \vphantom{\int} & (1.4)
\cr
{} & \normal^{\T} \condex\,\nabla \potex \,=\, 0 \quad \hbox{for all } \,(x,t) \in \partial \Omega \times [\,0\,,\, T\,]\,; \vphantom{\int} & (1.5)
\cr
{} & \pottr(x,0) \,= \,\potin(x,0) - \potex(x,0) \,=\, \Phi_0(x) \quad \hbox{and}  \quad \gat (x,0) \,=\, W_0 (x) \quad \hbox{for a.~a.~}\, x \in \Omega\, \vphantom{\int} & (1.6)
\cr}$$
on a bounded domain $\Omega \subset \R 3$ with the fixed time horizon $T > 0$. Further, we consider the monodomain system
$$\eqalignno{
{} & {\partial \pottr \over \partial t} + \ioncurr (\pottr, \gat) - {\la \over 1 + \la}\, \div \bigl(\,\condin\,\nabla \pottr\,\bigr) \,=\, {1 \over 1 + \la }\, \bigl(\, \la\,\excin - \excex\,\bigr) \quad \hbox{for a.~a.~}(x,t) \in \Omega \times [\, 0\,,\,T\,]\,; & (1.7)
\cr
{} & {\partial \gat \over \partial t} + G(\pottr, \gat) \,=\,  0\quad \hbox{for a.~a.~}(x,t) \in \Omega \times [\, 0\,,\,T\,]\,; & (1.8)
\cr
{} & \normal^{\T} \condin\,\nabla \pottr \,=\, 0 \quad \hbox{for all }(x,t) \in \partial \Omega \times [\,0\,,\, T\,]\,; \vphantom{\int} & (1.9)
\cr
{} & \pottr(x,0) \,= \,\Phi_0(x)\,\,\, \hbox{and} \,\,\, \gat (x,0) \,=\, W_0 (x) \quad \hbox{for a.~a.~} x \in \Omega  & (1.10)
\cr}$$
arising as a special case of (1.1)$\,-\,$(1.6) if the conductivity tensors satisfy $\condex = \lambda \,\condin$ with a constant parameter $\lambda > 0$, thus allowing to eliminate $\potex$ as an independent variable. In a series of papers,$\,$\footnote{$^{02)}$}{$[\,${\ninekpt Kunisch/Wagner 12}$\,]\,$, $[\,${\ninekpt Kunisch/Wagner 13a}$\,]$ and $[\,${\ninekpt Kunisch/Wagner 13b}$\,]\,$.} the authors investigated optimal control problems related to the dynamics $(1.1)\,-\,(1.6)$ and $(1.7)\,-\,(1.10)$ together with standard two-variable ionic models, namely the Rogers-McCulloch, FitzHugh-Nagumo and the linearized Aliev-Panfilov model (see Subsection 2.a) below). Using $\excex$ as control variable while $\excin = \zero$,$\,$\footnote{$^{03)}$}{This setting is due to physiological reasons.} and relying on a weak solution concept for the monodomain as well as for the bidomain system, the authors studied existence of minimizers and derived first-order necessary optimality conditions.
\par
The analysis of the control problems is based on a regularity discussion for the weak solutions, which leads to a stability estimate and a uniqueness theorem for the monodomain and bidomain system, respectively. The existence proof for the adjoint system is influenced by the regularity of the primal solutions as well. Unfortunately, the authors recognized a serious error within the proofs of these theorems. However, the present investigation shows that the assertions from $[\,${\kpt Kunisch/Wagner 12}$\,]$ and $[\,${\kpt Kunisch/Wagner 13a}$\,]$ can be maintained (with minor changes only) while the proofs must be subjected to considerable alterations. As a consequence, the optimization theorems from $[\,${\kpt Kunisch/Wagner 13b}$\,]$ allow for substantial improvements (see Theorems 4.2., 4.3.~and 4.5. below).
\par
In the present paper, we provide a refined regularity discussion and corrected proofs. For sake of completeness, all theorems from $[\,${\kpt Kunisch/Wagner 12}$\,]\,$, $[\,${\kpt Kunisch/Wagner 13a}$\,]$ and $[\,${\kpt Kunisch/Wagner 13b}$\,]$ will be repeated together with their possible corrections. In this updated version,$\,$\footnote{$^{04)}$}{A former version of this paper appeared at Sept 24, 2014.} further errors have been fixed, the regularity discussion is completed by Theorem 1.5.~below, and the consequences for the analysis of the optimal control problem in the bidomain case have been stated and proved.
\medskip
Our main results read as follows:
\medskip
{\bf Theorem 1.1.~(Stability estimate for weak solutions of the monodomain system)}$\,$\footnote{$^{05)}$}{Correction of $[\,${\ninekpt Kunisch/Wagner 12}$\,]\,$, p.~1533, Theorem 3.8.} {\it Consider the monodomain system in its weak formulation $(2.9)\,-\,(2.11)$, assuming that $\Omega \subset \R 3$ is a bounded Lipschitz domain, and $\condin \,\colon\,\,\cl(\Omega) \to \R {3 \times 3}$ is a symmetric, positive definite matrix function with $\L {\infty} {} \home$-coefficients, which obeys a uniform ellipticity condition with $\mu_1$, $\mu_2 > 0$:
$$  0 \le \mu_1 \,{\norm{\xi}}^2 \le \xi^{\T} \condin (x) \,\xi \le \mu_2 \,{\norm{\xi}}^2 \quad \forall\,\xi \in \R 3\,\,\, \forall\, x \in \Omega\,. \eqno (1.11)
$$
Let us specify either the Rogers-McCulloch, the FitzHugh-Nagumo or the linearized Aliev-Panfilov model.
\smallskip
If two weak solutions $(\pottr^\p, \gat^\p)$, $(\pottr^\pp, \gat^\pp) \in \bigl(\, \cont {0} {} \bigl[\, [\,0\,,\, T\,]\,,$ $ {\L {2} {} \home}\,\bigr] \,\cap\, \bochnerl {2} {(\,0\,,\,T\,)} {\sob {1,2} {} \home} $ $\cap\, \L {4} {} (\Omega_T) \,\bigr) \,\times \, \bochnerc {0} {[\,0\,,\, T\,]} {\L {2} {} \home} $ of the system correspond to initial values $\Phi^\p_0 = \Phi^\pp_0 = \Phi_0 \in \L {2} {} \home$, $W^\p_0 = W^\pp_0 = W_0 \in \L {4} {} \home$ and inhomogeneities $\excin^\p$, $\excex^\p$, $\excin^\pp $ and $\excex^\pp \in \L {2} {} \bigl[\, {(\,0\,,\,T\,)}\,,$ $ {\bigl(\,\sob {1,2} {} \home\,\bigr)^*}\,\bigr]\, $, whose norms are bounded by $R > 0$, then the following estimates hold:
$$\eqalignno{
{} & \hbox{\kern-20pt} {\norm{\pottr^\p - \pottr^\pp}}^2_{\bochnerc {0}
{[\,0\,,\,T\,]} {\L {2} {} \home} } + {\norm{\pottr^\p - \pottr^\pp }}^2_{\bochnerl {2} {(\,0\,,\,T\,)} {\sob {1,2} {} \home} } & (1.12)
\cr
{} & + \,{\norm{\gat^\p - \gat^\pp}}^2_{\bochnerc {0} {[\,0\,,\,T\,]} {\L {2} {} \home} }
+ {\norm{\gat^\p - \gat^\pp}}^2_{\bochnersob {1,2} {(\,0\,,\,T\,)} {\L {2} {} \home} } \phantom{\int}\cr
{} & & \displaystyle \le\,C \,\Bigl(\, {\norm{\excin^\p - \excin^\pp }}^2_{\bochnerl
{2} {(\,0\,,\,T\,)}  {\bigl(\,\sob {1,2} {} \home\,\bigr)^*} } + {\norm{\excex^\p - \excex^\pp}}^2_{\bochnerl {2} {(\,0\,,\,T\,)} {\bigl(\,\sob {1,2} {} \home\,\bigr)^*} } \,\Bigr)\,.
\cr}$$
The constant $C > 0$ does not depend on $\excin^\p$, $\excex^\p$, $\excin^\pp$ and $\excex^\pp$ but possibly on $\Omega$, $R$, $\Phi_0$ and $W_0$.}%
\medskip
}%
\pageno=2
{\bf Theorem 1.2.~(Stability estimate for weak solutions of the bidomain system)}$\,$\footnote{$^{06)}$}{Correction of $[\,${\ninekpt Kunisch/Wagner 13a}$\,]\,$, p.~959, Theorem 2.7., and $[\,${\ninekpt Kunisch/Wagner 13b}$\,]\,$, p.~1082, Theorem 2.4.} {\it Consider the bidomain system in its weak formulation $(2.75)\,-\,(2.78)$, assuming that $\Omega \subset \R 3$ is a bounded Lipschitz domain, $\condin$, $\condex\,\colon\,\,\cl(\Omega) \to \R {3 \times 3}$ are symmetric, positive definite matrix functions with $\L {\infty} {} \home$-coefficients, obeying uniform ellipticity conditions with $\mu_1$, $\mu_2 > 0$:
$$\eqalignno{
{} & \hbox{\kern-20pt} 0 \le \mu_1 \,{\norm{\xi}}^2 \le \xi^{\T} \condin (x) \,\xi \le \mu_2 \,{\norm{\xi}}^2 \,\,\hbox{ and } \,\, 0 \le \mu_1 \,{\norm{\xi}}^2 \le \xi^{\T} \condex (x) \,\xi \le \mu_2 \,{\norm{\xi}}^2 \quad \forall\,\xi \in \R 3\,\,\, \forall\, x \in \Omega\,. & (1.13)
\cr}$$
Let us specify either the Rogers-McCulloch, the FitzHugh-Nagumo or the linearized Aliev-Panfilov model.
\smallskip
If two weak solutions $(\pottr^\p, \potex^\p, \gat^\p)$, $(\pottr^\pp, \potex^\pp, \gat^\pp)$ correspond to initial values $\Phi^\p_0 = \Phi^\pp_0 = \Phi_0 \in \L {2} {} \home$, $W^\p_0 = W^\pp_0 = W_0 \in \L {4} {} \home$ and inhomogeneities $\excin^\p$, $\excex^\p$, $\excin^\pp $ and $\excex^\pp \in \bochnerl {2} {(\,0\,,\,T\,)} {\bigl(\,\sob {1,2} {} \home\,\bigr)^*} $, which satisfy the compatibility conditions
$$ \int_\Omega \Bigl(\,\excin^\p (x,t) + \excex^\p(x,t)\,\Bigr)\, dx \,=\, \int_\Omega \Bigl(\,\excin^\pp (x,t) + \excex^\pp(x,t)\,\Bigr)\, dx\,=\, 0 \quad \hbox{for a.~a.~} \, t \in (\,0\,,\, T\,)\,, \eqno (1.14)
$$
and whose norms are bounded by $R > 0$, then the following estimates hold:
$$\eqalignno{
{} & \hbox{\kern-20pt} {\norm{\pottr^\p - \pottr^\pp}}^2_{\bochnerc {0} {[\,0\,,\,T\,]} {\L {2} {} \home} } + {\norm{\pottr^\p - \pottr^\pp }}^2_{\bochnerl {2} {(\,0\,,\,T\,)} {\sob {1,2} {} \home} } + {\norm{\potex^\p - \potex^\pp }}^2_{\bochnerl {2} {(\,0\,,\,T\,)} {\sob {1,2} {} \home} } \phantom{\int} & (1.15)
\cr
{} & + \,{\norm{\gat^\p - \gat^\pp}}^2_{\bochnerc {0} {[\,0\,,\,T\,]} {\L {2} {} \home} } + {\norm{\gat^\p - \gat^\pp }}^2_{\bochnersob {1,2} {(\,0\,,\,T\,)} {\L {2} {} \home} } \phantom{\int}
\cr
{} & & \displaystyle \le\,C \,\Bigl(\, {\norm{\excin^\p - \excin^\pp }}^2_{\bochnerl {2} {(\,0\,,\,T\,)} {\bigl(\,\sob {1,2} {} \home\,\bigr)^*} } + {\norm{\excex^\p - \excex^\pp }}^2_{\bochnerl {2} {(\,0\,,\,T\,)} {\bigl(\,\sob {1,2} {} \home\,\bigr)^*} } \,\Bigr)\,. \phantom{\int}
\cr}$$
The constant $C > 0$ does not depend on $\excin^\p$, $\excex^\p$, $\excin^\pp$ and $\excex^\pp$ but possibly on $\Omega$, $R$, $\Phi_0$ and $W_0$.}%
\medskip
\medskip
{\bf Theorem 1.3.~(Uniqueness of weak solutions of the monodomain system)}$\,$\footnote{$^{07)}$}{Correction of $[\,${\ninekpt Kunisch/Wagner 12}$\,]\,$, p.~1529, Theorem 2.5.}  {\it Consider the mono\-domain system in its weak formulation $(2.9)\,-\,(2.11)$ under the assumptions of Theorem 1.1. Specifying either the Rogers-McCulloch, the FitzHugh-Nagumo or the linearized Aliev-Panfilov model, the system admits a unique weak solution
$$\eqalignno{
{} & (\pottr, \gat) \,\in \, \Bigl(\,\bochnerc {0} {[\,0\,,\, T\,]} {\L {2} {} \home} \,\cap\,\bochnerl
{2} {(\,0\,,\,T\,)} {\sob {1,2} {} \home} \,\cap\, \L {4} {} (\Omega_T) \,\Bigr)\,\times\, \bochnerc {0} {[\,0\,,\, T\,]} {\L {2} {} \home} & (1.16)
\cr}$$
in correspondence to initial values $\Phi_0 \in \L {2} {} \home$, $W_0 \in \L {4} {} \home$ and inhomogeneities $\excin$, $\excex \in \L {2} {} \bigl[\,{(\,0\,,\,T\,)}\,,$ $ {\bigl(\,\sob {1,2} {} \home\,\bigr)^*} \,\bigr]\,$.}
\medskip
{\bf Theorem 1.4.~(Uniqueness of weak solutions of the bidomain system)}$\,$\footnote{$^{08)}$}{Correction of $[\,${\ninekpt Kunisch/Wagner 13a}$\,]\,$, p.~959, Theorem 2.8., and $[\,${\ninekpt Kunisch/Wagner 13b}$\,]\,$, p.~1082, Theorem 2.3.} {\it Consider the bidomain system in its weak formulation $(2.75)\,-\,(2.78)$ under the assumptions of Theorem 1.2. Specifying either the Rogers-McCulloch, the FitzHugh-Nagumo or the linearized Aliev-Panfilov model, the system admits a unique weak solution
$$\eqalignno{
{} & (\pottr, \potex, \gat) \in \Bigl(\,\bochnerc {0} {[\,0\,,\, T\,]} {\L {2} {} \home} \,\cap\,\bochnerl {2} {(\,0\,,\,T\,)} {\sob {1,2} {} \home} \,\cap\, \L {4} {} (\Omega_T)\,\Bigr) & (1.17)
\cr
{} & & \displaystyle \times \, \bochnerl {2} {(\,0\,,\,T\,)} {\sob {1,2} {} \home} \,\times \, \bochnerc {0} {[\,0\,,\, T\,]} {\L {2} {} \home} \phantom{\int}
\cr}$$
in correspondence to initial values $\Phi_0 \in \L {2} {} \home$, $W_0 \in \L {4} {} \home$ and inhomogeneities $\excin$, $\excex \in \L {2} {} \bigl[\, {(\,0\,,\,T\,)}\,,$ $ {\bigl(\,\sob {1,2} {} \home\,\bigr)^*} \,\bigr]\,$.}%
\medskip
The main error to be corrected was the claim that, already under the assumptions of Theorems 1.1.$\,-\,$1.4., the transmembrane potential $\pottr$ within a weak solution of $(1.1)\,-\,(1.6)$ or $(1.7)\,-\,(1.10)$ can admit $\L {4} {} \bigl[\,{(\,0\,,\,T\,)}\,,$ $ {\L {6} {} \home} \,\bigr] $- or even $\bochnerl {\infty} {(\,0\,,\,T\,)} {\L {4} {} \home} $-regularity.$\,$\footnote{$^{09)}$}{This error can be traced back to $[\,${\ninekpt Kunisch/Wagner 12}$\,]\,$, p.~1534, (3.39), and p.~1544, (B.14), as well as to $[\,${\ninekpt Ku\-nisch/Wagner 13a}$\,]\,$, p.~964, (2.69) and (2.70), respectively.} In Section 3 below, we will see that the theorems can be proven without relying on this claim. However, as the following theorem states, in order to ensure the claimed $\L {4} {} \bigl[\,{(\,0\,,\,T\,)}\,,$ $ {\L {6} {} \home} \,\bigr] $-regularity of $\pottr$, it suffices to consider right-hand sides $\excin$ and $\excex$ belonging to the space $\bochnerl {4} {(\,0\,,\,T\,)} {\L {2} {} \home} $ instead of $\bochnerl {2} {(\,0\,,\,T\,)} {\bigl(\,\sob {1,2} {} \home\,\bigr)^*} $,
\medskip
{\bf Theorem 1.5.~(Higher regularity of the transmembrane potential in weak solutions of the monodomain and the bidomain system)} {\it 1) Consider the mono\-domain system in its weak formulation $(2.9)\,-\,(2.11)$ under the assumptions of Theorem 1.1.~with either the Rogers-McCulloch, FitzHugh-Nagumo or the linearized Aliev-Panfilov model and initial values $\Phi_0 \in \L {2} {} \home$, $W_0 \in \L {2} {} \home$. If a weak solution $(\pottr, \gat)$ corresponds to inhomogeneities $\excin$, $\excex \in \bochnerl {4} {(\,0\,,\,T\,)} {\L {2} {} \home} $ then the transmembrane potential admits the higher regularity $\pottr \in \bochnerl {4} {(\,0\,,\,T\,)} {\sob {1,2} {} \home} $.
\smallskip
2) Consider the bidomain system in its weak formulation $(2.75)\,-\,(2.78)$ under the assumptions of Theo\-rem 1.2.~with the either the Rogers-McCulloch, FitzHugh-Nagumo or the linearized Aliev-Panfilov model and initial values $\Phi_0 \in \L {2} {} \home$, $W_0 \in \L {2} {} \home$. If a weak solution $(\pottr, \potex, \gat)$ corresponds to inhomogeneities $\excin$, $\excex \in \bochnerl {4} {(\,0\,,\,T\,)} {\L {2} {} \home} $ then the transmembrane potential admits the higher regularity $\pottr \in \bochnerl {4} {(\,0\,,\,T\,)} {\sob {1,2} {} \home} $.}%
\medskip
The paper is structured as follows. We continue with a short collection of notations (Subsect.~1.b)$\,$) and repeat, for the reader's sake, the imbedding theorems for Bochner spaces used below (Subsect.~1.c)$\,$). In Section 2, we start with the desciption of the ionic models, which will be subsequently used (Subsect.~2.a)$\,$). Then we restate the monodomain system in its weak formulation and study the existence and regularity of the weak solutions for the different models (Subsect.~2.b)$\,-\,$d)$\,$. Subsequently, the weak formulation of the full bidomain system together with existence and regularity results for its weak solutions is provided (Subsect.~2.e)$\,-\,$g)$\,$). Even here, the proof of Theorem 1.5.~is included (Subsect.~2.h)$\,$). In Section 3, we deliver the corrected proofs of the stability estimates for the monodomain and bidomain system, respectively (Theorems 1.1.~and 1.2.), which imply the uniqueness of the weak solutions (Theorems 1.3.~and 1.4.). Some remarks and corollaries are added. Finally, in Section 4, we list in full detail the corrections to be made in the authors' previous papers.
\bigskip
{\bf b) Notations.}
\medskip
We abbreviate $\Omega \times [\,0\,,\,T\,]$ by $\Omega_T$. By $\L p {} \home$, we denote the space of functions, which are in the $p$th power integrable ($\,1 \le p < \infty$), or are measurable and essentially bounded $(p=\infty)$, and by $\sob {1,p} {} \home$ the Sobolev space of functions $\psi\,\colon \,\,\Omega \to \R {}$ which, together with their first-order weak partial derivatives, belong to the space $\L p {} \homer$ ($\,1 \le p < \infty$). Concerning spaces of Bochner integrable mappings, e.~g.~$\bochnerl {2} {(\,0\,,\,T\,)} {\sob {1,2} {} \home} $, we refer to $[\,${\kpt Kunisch/Wagner 12}$\,]\,$, p.~1542. The gradient $\nabla$ is always taken only with respect to the spatial variables $x$. Finally, we use the nonstandard abbreviation \lq$\fforall\, t \in \A$", which has to be read as \lq for almost all $t \in \A$" or \lq for all $t \in \A$ except for a Lebesgue null set". The symbol $\zero$ denotes, depending on the context, the zero element or the zero function of the underlying space.
\bigskip
{\bf c) Compact imbeddings of Bochner spaces.}
\medskip
{\bf Theorem 1.6.~(Aubin-Dubinskij lemma)}$\,$\footnote{$^{10)}$}{$[\,${\ninekpt Dubinskij 65}$\,]\,$, p.~612, Teorema 1, and p.~615, Teorema 2. Its formulation is not affected by the corrections, which have been presented recently in $[\,${\ninekpt Barrett/S\"uli 12}$\,]\,$. Note that this theorem has been cited incorrectly in $[\,${\ninekpt Fursikov}$\,]\,$, p.~8, Lemma 1.2., $[\,${\ninekpt Kunisch/Wagner 12}$\,]\,$, p.~1542, Theorem A.6, and the former version of this paper, p.~4, Theorem 1.5.} {\it Consider three normed spaces $\X_0 \subseteq \X \subseteq \X_1$ where the imbedding $\X_0 \hookrightarrow \X$ is compact and the imbedding $\X \hookrightarrow \X_1$ is continuous. If $p$, $p^\p \in (\,1 \,,\,\infty\,)$ then the space
$$ \Y \,=\, \bigl\{\,f \in \bochnerl {p} {(\,0\,,\, T\,)} {\X_0} \,\,\big\vert \,\, {df \over dt} \in \bochnerl {p^\p} {(\,0\,,\, T\,)} {\X_1} \,\bigr\} \eqno (1.18)
$$
is compactly imbedded into $\bochnerl {p} {(\,0\,,\, T\,)} {\X} $. $\blackbox$}
\medskip
{\bf Theorem 1.7.~(Generalization of the Aubin-Dubinskij lemma)}$\,$\footnote{$^{11)}$}{$[\,${\ninekpt Simon 87}$\,]\,$, p.~90, Corollary 8.} {\it Consider three Banach spaces $\X_0 \subseteq \X \subseteq \X_1$ where the imbeddings $\X_0 \hookrightarrow \X$ and $\X \hookrightarrow \X_1$ are continuous while $\X_0 \hookrightarrow \X_1$ is compact. Assume further that there exists a number $0 < \vartheta < 1$ such that
$$\eqalignno{
{} & {\norm{\psi}}_{\X} \,\le \, C \, {\norm{\psi}}^{1 - \vartheta}_{\X_0} \cdot {\norm{\psi}}^\vartheta_{\X_1} \quad \forall \,\psi \in \X_0 \,\cap\,\X_1\,. & (1.19)
\cr}$$
Consider for exponents $1 \le p < p^\p \le \infty$ the space
$$ \Y \,=\, \bigl\{\,f \in \bochnerl {p} {(\,0\,,\, T\,)} {\X_0} \,\,\big\vert \,\,
{df \over dt} \in \bochnerl {p^\p} {(\,0\,,\, T\,)} {\X_1} \,\bigr\}\,. \eqno (1.20)
$$
1) If $\Delta(p, p^\p, \vartheta) = (1 - \vartheta)/p - \vartheta\,(1 - 1/p^\p) \ge 0$ then the space $\Y$ is compactly imbedded into $\bochnerl {q} {(\,0\,,\, T\,)} {\X} $ for all $1 \le q < 1/\Delta$.
\smallskip
2) If $\Delta (p, p^\p, \vartheta)= (1 - \vartheta)/p - \vartheta\,(1 - 1/p^\p) < 0$ then the space $\Y$ is compactly imbedded into $\bochnerc {0} {[\,0\,,\, T\,]} {\X} $.  $\blackbox$}%
\bigskip\medskip
{\kapitel 2.~Regularity of weak solutions for the monodomain and bidomain system.}
\medskip
{\bf a) The ionic models.}
\medskip
The following models for the ionic current $\ioncurr$ and the function $G$ within the gating equation will be considered:
\medskip
{\it 1) The FitzHugh-Nagumo model.}$\,$\footnote{$^{12)}$}{$[\,${\ninekpt FitzHugh 61}$\,]\,$, together with $[\,${\ninekpt Nagumo/Arimoto/Yoshizawa 62}$\,]\,$.}
$$\eqalignno{
{} & \ioncurr(\varphi, w ) \,=\, \varphi\,(\varphi - a)\,(\varphi - 1) + w \,=\, \varphi^3 - (a+1)\,\varphi^2 + a\, \varphi + w\,; & (2.1)
\cr
{} & G(\varphi,w) \,=\, \e\,w - \e\,\kappa\,\varphi & (2.2)
\cr}$$
with $0 < a < 1$, $\kappa > 0$ and $\e > 0$. Thus the gating variable obeys the linear ODE
$$ \partial \gat / \partial t + \e\,\gat \,=\, \e\,\kappa\,\pottr\,. \eqno (2.3)
$$
\par
{\it 2) The Rogers-McCulloch model.}$\,$\footnote{$^{13)}$}{$[\,${\ninekpt Rogers/McCulloch 94}$\,]\,$.}
$$\eqalignno{
{} & \ioncurr(\varphi, w) \,=\, b \cdot \varphi\,(\varphi - a)\,(\varphi - 1) + \varphi \cdot w \,=\, b\,\varphi^3 - (a+1)\, b\,\varphi^2 + a\,b\, \varphi + \varphi\, w\,; & (2.4)
\cr
{} & G(\varphi, w) \,=\, \e\,w - \e\,\kappa\, \varphi & (2.5)
\cr}$$
with $0 < a < 1$, $b > 0$, $\kappa > 0$ and $\e > 0$. Consequently, the ODE for the gating variable is the same as before.
\medskip
{\it 3) The linearized Aliev-Panfilov model.}$\,$\footnote{$^{14)}$}{This model is taken from $[\,${\ninekpt Bourgault/Coudi\`ere/Pierre 09}$\,]\,$, p.~480. Instead, the original model from $[\,${\ninekpt Aliev/ Pan\-filov 96}$\,]$ contains a Riccati equation for the gating variable.}
$$\eqalignno{
{} & \ioncurr(\varphi, w) \,=\, b \cdot \varphi\,(\varphi - a)\,(\varphi - 1) + \varphi \cdot w \,=\, b\,\varphi^3 - (a+1)\, b\,\varphi^2 + a\,b\, \varphi + \varphi\, w\,; & (2.6)
\cr
{} & G(\varphi,w) \,=\, \e\,w - \e\,\kappa\,\bigl(\,(a+1)\,\varphi - \varphi^2\,\bigr) & (2.7)
\cr}$$
with $0 < a < 1$, $b > 0$, $\kappa > 0$ and $\e > 0$. The linear ODE for the gating variable is
$$ \partial \gat / \partial t + \e\,\gat \,=\, \e\,\kappa\,\bigl( \, (a+1)\,\pottr - \pottr^2\,\bigr)\,. \eqno (2.8)
$$
\medskip
{\bf b) Weak formulation of the monodomain system and known regularity of weak solutions.}
\medskip
The weak formulation of the monodomain system $(1.7)\,-\,(1.10)$ reads as follows:
$$\eqalignno{
{} & \int_\Omega \Bigl(\,{\partial \pottr \over \partial t }+ \ioncurr(\pottr,
\gat)\,\Bigr)\, \psi \, dx + \int_\Omega {\la \over 1 + \la}\, \nabla \psi^{\T}
\condin \,\nabla \pottr \, dx \,=\,\int_\Omega {1 \over 1 + \la }\, \bigl(\, \la\,\excin - \excex \,\bigr) \,\psi \,dx  & (2.9)
\cr
{} & & \displaystyle \quad \forall\, \psi \in \sob {1,2} {} \home \,\,\, \fforall \,
t \in [\, 0\,,\,T\,]\,;
\cr
{} & \int_\Omega \Bigl(\,{\partial \gat \over \partial t} + G(\pottr,
\gat)\,\Bigr)\, \psi \, dx\,=\, 0 \quad \forall\, \psi \in \L {2} {} \home\,\,\, \fforall \, t \in [\,0\,,\,T\,] & (2.10)
\cr
{} & \pottr(x,0) \,= \,\Phi_0(x)\quad \fforall\, x \in \Omega\,; \quad  \gat (x,0)
\,=\, W_0 (x) \quad \fforall\, x \in \Omega\, \vphantom{\int} & (2.11)
\cr}$$
\par
where $\la > 0$. Under the assumptions of Theorem 1.1., the system $(2.9)\,-\,(2.11)$ with either the FitzHugh-Nagumo, the Rogers-McCulloch or the linearized Aliev-Panfilov model admits for arbitrary initial values $\Phi_0$, $W_0 \in \L {2} {} \home$ and inhomogeneities $\excin$, $\excex \in \bochnerl {2} {(\,0\,,\, T\,)} {\bigl(\, \sob {1,2} {} \home\,\bigr)^*} $ at least one weak solution$\,$\footnote{$^{15)}$}{$[\,${\ninekpt Kunisch/Wagner 12}$\,]\,$, p.~1528 f., Theorem 2.2.}
$$\eqalignno{
{} & (\pottr, \gat ) \in \Bigl(\,\bochnerc {0} {[\,0\,,\,T\,]} {\L {2} {} \home } \,\cap \, \bochnerl {2} {(\,0\,,\,T\,)} {\sob {1,2} {} \home} \,\cap\, \L {4} {} (\Omega_T)\,\Bigr) \,\times \, \bochnerc {0} {[\,0\,,\,T\,]} {\L {2} {} \home} \,. & (2.12)
\cr}$$
Any weak solution satisfies the a-priori estimate$\,$\footnote{$^{16)}$}{$[\,${\ninekpt Kunisch/Wagner 12}$\,]\,$, p.~1529, Theorem 2.4.}
$$\eqalignno{
{} & \hbox{\kern-20pt} {\norm{\pottr}}^2_{ \bochnerc {0} {[\,0\,,\, T\,]} {\L {2} {} \home} } + {\norm{\pottr}}^2_{ \bochnerl {2} {(\,0\,,\,T\,)} {\sob {1,2} {} \home} } + {\norm{\pottr}}^4_{\L {4} {} (\Omega_T) \vphantom{\bigl[}} +  {\norm{\partial
\pottr/\partial t}}^{4/3}_{\bochnerl {4/3} {(\,0\,,\,T\,)} {\bigl(\,\sob {1,2} {} \home\,\bigr)^*} }
\cr
{} & & \displaystyle + \,{\norm{\gat}}^2_{\bochnerc {0} {[\,0\,,\, T\,]} {\L {2} {} \home} }  + {\norm{\partial \gat /\partial t}}^2_{\bochnerl {2} {(\,0\,,\,T\,)} {\bigl(\,\sob {1,2} {} \home\,\bigr)^*} }
\cr
{} & \hbox{\kern-20pt}  \le \, C\cdot \Bigl(\,1 + {\norm{\Phi_0}}^2_{\L {2} {} \home\vphantom{\bigl[}} + {\norm{W_0}}^2_{\L {2} {} \home \vphantom{\bigl[}} + {\norm{\excin}}^2_{\bochnerl {2} {(\,0\,,\,T\,)} {\bigl(\,\sob {1,2} {} \home\,\bigr)^*} } + \, {\norm{\excex}}^2_{\bochnerl {2} {(\,0\,,\,T\,)} {\bigl(\,\sob {1,2} {} \home\,\bigr)^*} }  \,\Bigr) & (2.13)
\cr}$$
with a constant $C >0$, which does not depend on $\Phi_0$, $W_0$, $\excin$ and $\excex$. We will investigate now how to improve the regularity of a given weak solution, depending on the model and the regularity of $W_0$. In the following, we will refer to
$$\eqalignno{
{} & {\cal M} (\psi^\p, \psi^\pp) \,=\, {\la \over 1 + \la}\, \int_\Omega \nabla (\psi^\p)^{\T} \condin \,\nabla \psi^\pp \, dx\, & (2.14)
\cr}$$
as the {\it monodomain form} ${\cal M} \,\colon \,\, \sob {1,2} {} \home \times \sob {1,2} {} \home \to \R {}$.
\bigskip
{\bf c) Monodomain system with Rogers-McCulloch model: improvement of regularity for the weak solutions.}
\medskip
{\bf Proposition 2.1.~(Gain of regularity for the gating variable)} {\it Consider the monodomain system with the Rogers-McCulloch model under the assumptions of Theorem 1.1. Let a weak solution $(\pottr, \gat)$ of it correspond to inital values $\Phi_0\in \L {2} {} \home$, $W_0 \in \L {4} {} \home$ and right-hand sides $\excin$, $\excex \in \L {2} {} \bigl[\, {(\,0\,,\,T\,)}\,,$ $ {\bigl(\,\sob {1,2} {} \home\,\bigr)^*}\,\bigr]\,$, thus belonging to the spaces in $(2.12)$.
\smallskip
1) Then $\gat$ belongs to $\bochnerc {1} {(\,0\,,\,T\,)} {\L {2} {} \home} $; consequently, $\partial \gat / \partial t$ belongs to $\L {2} {} (\Omega_T)$.
\smallskip
2) Moreover, $\,W_0 \in \L {4} {} \home$ implies that $\gat$ belongs even to $\bochnerc {0} {[\,0\,,\,T\,]} {\L {4} {} \home} $, and it holds that
$$\eqalignno{
{} & {\norm{\gat}}^4_{\bochnerc {0} {[\,0\,,\,T\,]} {\L {4} {} \home} }\,  \le\, C\,\Bigl(\,1 + {\norm{\Phi_0}}^2_{\L {2} {} \home} + {\norm{W_0}}^4_{\L {4} {} \home} & (2.15)
\cr
{} & & \displaystyle  + \, {\norm{\excin}}^2_{\bochnerl {2} {(\,0\,,\,T\,)} {\bigl(\,\sob {1,2} {} \home\,\bigr)^*} }  +  {\norm{\excex}}^2 _{\bochnerl {2} {(\,0\,,\,T\,)} {\bigl(\,\sob {1,2} {} \home\,\bigr)^*} }  \,\Bigr) \,.
\cr}$$}%
\par
Note that Part 1) of Proposition 2.1.~is still true for $ W_0 \in \L {2} {} \home$.
\medskip
{\bf Proposition 2.2.~(Gain of regularity for the transmembrane potential, strong solution of the monodomain system)} {\it Consider the monodomain system with the Rogers-McCulloch model under the assumptions of Theorem 1.1. If, moreover, $\partial\Omega$ is of $\cont {1,1} {}$-regularity, the coefficients of $\condin$ belong to $\sob {1, \infty} {} \home$, $\Phi_0 \in \sob {1,2} {} \home$, $W_0 \in \L {4} {}\home$ and $\excin$, $\excex \in \L {2} {} (\Omega_T)$ then the system admits even a strong solution $(\pottr, \gat)$  with
$$\eqalignno{
{} & \pottr \,\in \, \bochnerc {0} {[\,0\,,\,T\,]} {\sob {1,2} {} \home } \,\cap \, \bochnerl {2} {(\,0\,,\,T\,)} {\sob {2,2} {} \home} \,\cap\,\bochnersob {1, 4/3} {(\,0\,,\,T\,)} {\L {4/3} { } \home} \,. & (2.16)
\cr}$$}%
\par
If the Rogers-McCulloch model is replaced by the FitzHugh-Nagumo model then Propositions 2.1.~and 2.2.~hold accordingly. The most important regularity property of the solution, however, is stated in Theorem 1.5., 1). We continue with the proofs of Propositions 2.1.~and 2.2.~and postpone the proof of Theorem 1.5.~to Subsect.~2.h).
\medskip
{\bf Proof of Proposition 2.1.} {\it Part 1)\/} Observe first that $\gat$ admits the representation
$$\eqalignno{
{} & \gat (x,t) \,=\, W_0 (x)\, e^{- \e t} + \e\,\kappa \, e^{- \e t} \int^t_0 \pottr(x, \tau )\, e^{\e \tau} \, d\tau \,, & (2.17)
\cr}$$
from which the claimed $\bochnerc {1} {(\,0\,,\,T\,)} {\L {2} {} \home} $-regularity follows by differentiation.
\medskip
{\it Part 2)\/} Recall now that $W_0$ belongs to $\L {4} {} \home$. Then from $(2.13)$ and $(2.17)$, we derive the further estimate
$$\eqalignno{
{} & \hbox{\kern-20pt} \int_\Omega \gat (t)^4 \, dx \,\le\, C\,{\norm{W_0}}^4_{\L {4} {} \home} + C\,\e\,\kappa \,{\norm{\pottr }}^4_{\L {4} {} (\Omega_T)} & (2.18)
\cr
{} & \le\, C\,\Bigl(\,1 + {\norm{\Phi_0}}^2_{\L {2} {} \home} + {\norm{W_0}}^4_{\L {4} {} \home} +  {\norm{\excin}}^2_{\bochnerl {2} {(\,0\,,\,T\,)} {\bigl(\,\sob {1,2} {} \home\,\bigr)^*} }  +  {\norm{\excex}}^2 _{\bochnerl {2} {(\,0\,,\,T\,)} {\bigl(\,\sob {1,2} {} \home\,\bigr)^*} }  \,\Bigr) \,,
\cr}$$
which confirms that $\gat$ belongs to $\bochnerl {\infty} {(\,0\,,\,T\,)} {\L {4} {} \home} $. In order to prove that ${\norm{\gat (\,\cdot\,, t)}}^4_{\L {4} {} \home}$ depends continuously on $t$, consider for $s$, $t \in [\,0\,,\,T\,]\,$ the difference
$$\eqalignno{
{} & \Bigl\vert\,{\norm{\gat(s)}}^4_{\L {4} {} \home} - {\norm{\gat (t)}}^4_{\L {4} {} \home} \,\Bigr\vert\,=\,\Bigl\vert \, \int_\Omega \Bigl(\,\gat(s)^4 - \gat (t)^4\,\Bigr)\, dx \,\Bigr\vert & (2.19)
\cr
{} & =\, \Bigl\vert \,\int_\Omega \bigl(\, \gat(s) - \gat (t)\,\bigr)\,\bigl(\,\gat (s)^3 + \gat(s)^2 \,\gat (t) + \gat (s) \,\gat(t)^2 + \gat(t)^3\,\bigr)\, dx\,\Bigr\vert & (2.20)
\cr
{} & \le\, \Bigl(\, \int_\Omega \bigl\vert\, \gat(s) - \gat (t)\,\bigr\vert^4\, dx\,\Bigr)^{1/4} \,\Bigl(\, \int_\Omega \bigl\vert\, \gat (s)^3 + \gat(s)^2 \,\gat (t) + \gat (s) \,\gat(t)^2 + \gat(t)^3\,\bigr\vert^{4/3} \, dx\,\Bigr)^{3/4}\,. & (2.21)
\cr}$$
We rely on the following lemma:
\medskip
{\bf Lemma 2.3.~(Generalized Minkowski inequality)}$\,$\footnote{$^{17)}$}{$[\,${\ninekpt Stein 70}$\,]\,$, p.~271, A.1.} {\it Assume that $1 \le p < \infty$. Then for every measurable function $\varphi \,\colon \,\, (\,0\,,\,T\,) \times \Omega \to \R {}$, it holds that
$$\eqalignno{
{} & \Bigl(\, \int_\Omega \,\Bigl(\, \int^T_0 \bigl\vert \,\varphi(x,t)\,\bigr\vert \, dt\,\Bigr)^p \, dx\,\Bigr)^{1/p} \,\le\, \int^T_0 \,\Bigl(\, \int_\Omega\,\bigl\vert \,\varphi(x,t)\,\bigr\vert^p \, dx\,\Bigr)^{1/p}\, dt\,. \quad \blackbox & (2.22)
\cr}$$}%
\medskip
Applying Jensen's inequality and Lemma 2.3.~to the estimation of the first term from (2.21), we obtain
$$\eqalignno{
{} & \hbox{\kern-20pt}\Bigl(\,\int_\Omega \bigl\vert\, \gat(s) - \gat (t)\,\bigr\vert^4\, dx\,\Bigr)^{1/4} =\, \Bigl(\,\int_\Omega \,\Bigl\vert\, W_0\,\bigl(\, e^{-\e s} - e^{-\e t}\,\bigr) + \e\,\kappa\,e^{-\e s}\,\bigl(\, \int^s_0 \pottr \,e^{\e\tau}\, d\tau - \int^t_0 \pottr \, e^{\e\tau}\,d\tau\,\bigr) & (2.23)
\cr
{} & & \displaystyle + \,\e\,\kappa\,\bigl(\, e^{-\e s} - e^{-\e t} \,\bigr)\,\int^t_0 \pottr \,e^{\e\tau}\, d\tau\,\Bigr\vert^4 \, dx \,\Bigr)^{1/4}
\cr}$$
$$\eqalignno{
{} & \le\, C \,\bigl\vert\, e^{-\e s} - e^{-\e t}\,\bigr\vert\cdot {\norm{W_0}}_{\L {4} {} \home} + C\,\Bigl(\,\int_\Omega \Bigl\vert\,\int^s_t \pottr\,e^{\e\tau}\, d\tau\,\Bigr\vert^4\, dx\,\Bigr)^{1/4} & (2.24)
\cr
{} & & \displaystyle +\, C \,\bigl\vert \, e^{-\e s} - e^{-\e t}\,\bigr\vert \cdot \Bigl(\,\int_\Omega \Bigl\vert\,\int^t_0 \pottr\,e^{\e\tau}\, d\tau\,\Bigr\vert^4\, dx \,\Bigr)^{1/4}
\cr
{} & \le\, C \,\bigl\vert\, e^{-\e s} - e^{-\e t}\,\bigr\vert\cdot {\norm{W_0}}_{\L {4} {} \home} + C\,\vert\, s - t\,\vert ^{3/4}\,\Bigl(\,\int_\Omega \, \int^s_t \vert\, \pottr\,\vert^4 \,e^{4 \,\e\,\tau}\, d\tau\, dx\,\Bigr)^{1/4} & (2.25)
\cr
{} & & \displaystyle +\, C \,\bigl\vert \, e^{-\e s} - e^{-\e t}\,\bigr\vert \cdot \int^t_0 \Bigl(\, \int_\Omega \bigl\vert \, \pottr\,\bigr\vert^4 \,e^{4 \,\e \,\tau}\, dx\, \Bigr)^{1/4} \, d\tau
\cr
{} & \le\, C \,\bigl\vert\, e^{-\e s} - e^{-\e t}\,\bigr\vert\cdot {\norm{W_0}}_{\L {4} {} \home} + C\,\vert\, s - t\,\vert ^{3/4}\,\Bigl(\,\int_\Omega \, \int^T_0 \vert\, \pottr\,\vert^4 \,e^{4 \,\e\,T}\, d\tau\, dx\,\Bigr)^{1/4} & (2.26)
\cr
{} & & \displaystyle +\, C \,\bigl\vert \, e^{-\e s} - e^{-\e t}\,\bigr\vert \cdot \int^T_0 \Bigl(\, \int_\Omega \bigl\vert \, \pottr\,\bigr\vert^4 \,e^{4 \,\e \,T}\, dx\, \Bigr)^{1/4} \, d\tau
\cr
{} & \le \, C \,\bigl\vert\, e^{-\e s} - e^{-\e t}\,\bigr\vert\cdot \Bigl(\, {\norm{W_0}}_{\L {4} {} \home} + {\norm{\pottr}}_{\bochnerl {1} {(\,0\,,\,T\,)} {\L {4} {} (\Omega)} } \,\Bigr) + C\, \bigl\vert \,s - t\,\bigr\vert^{3/4} \cdot {\norm{\pottr}}_{\L {4} {} (\Omega_T)} & (2.27)
\cr
{} & \le \, C\, \Bigl(\,\bigl\vert \, s - t\,\bigr\vert^{3/4} + \bigl\vert\, e^{-\e s} - e^{-\e t}\,\bigr\vert\,\Bigr) \cdot \Bigl(\, {\norm{W_0}}_{\L {4} {} \home} + {\norm{\pottr}}_{\L {4} {} (\Omega_T)} \,\Bigr) \,. & (2.28)
\cr
}$$
Estimation of the second term yields:
$$\eqalignno{
{} & \hbox{\kern-20pt} \int_\Omega \bigl\vert\, \gat (s)^3 + \gat(s)^2 \,\gat (t) + \gat (s) \,\gat(t)^2 + \gat(t)^3\,\bigr\vert^{4/3} \, dx \, & (2.29)
\cr
{} & \le\, C\,\int_\Omega \Bigl(\, \vert \,\gat(s)\,\vert^4 + \vert \,\gat(s)\,\vert^{8/3}\, \vert \, \gat(t)\,\vert^{4/3} + \vert\, \gat(s)\,\vert^{4/3} \,\vert\, \gat(t)\,\vert^{8/3} + \vert\, \gat(t)\,\vert^4\,\Bigr)\, dx\,.
\cr}$$
From $(2.18)$ we already know that $\int_\Omega \vert \,\gat(s)\,\vert^4 \, dx$ and
$$\eqalignno{
{} & \int_\Omega \vert \,\gat(s)\,\vert^{8/3}\, \vert \, \gat(t)\,\vert^{4/3}\, dx \,\le \,\int_\Omega \Bigl(\, \Max\,\bigl(\,\vert \,\gat(s)\,\vert\,,\, \vert \, \gat(t)\,\vert\,\bigr)\,\Bigr)^4\, dx & (2.30)
\cr
{} & \le \, \int_\Omega \Max \,\Bigl(\,\vert \,\gat(s)\,\vert^4\,,\, \vert \,\gat(t)\,\vert^4\,\Bigr) \, dx
\, \le \,\int_\Omega \Bigl(\, \vert \,\gat(s)\,\vert^4 + \vert \,\gat(t)\,\vert^4\,\Bigr) \, dx & (2.31)
\cr}$$
remain uniformly bounded. Consequently, for every $\e > 0$ we may determine $\delta(\e) > 0$ such that $\vert \, s - t \,\vert \le \delta(\e)$ implies $\bigl\vert\,{\norm{\gat(s)}}^4_{\L {4} {} \home} - {\norm{\gat (t)}}^4_{\L {4} {} \home} \,\bigr\vert \le \e$, and $\gat$ belongs to $\bochnerc {0} {[\,0\,,\,T\,]} {\L {4} {} \home} $. $\blackbox$
\medskip
{\bf Proof of Proposition 2.2.} We abbreviate: $\excinex = \bigl(\, \la\,\excin - \excex\,\bigr) / (1 + \la)$.
\par
$\bullet$ {\bf Step 1.} {\it A-priori estimates for a strong solution.} Assume that $(\pottr, \gat) \in \Bigl(\,\bochnerc {0} {[\,0\,,\,T\,]} {\sob {1,2} {} \home } \,\cap \, \bochnerl {2} {(\,0\,,\,T\,)} {\sob {2,2} {} \home} \,\cap\, \bochnersob {1, 4/3} {(\,0\,,\,T\,)} {\L {4/3} {} \home} \,\Bigr) \times {\bochnerc {0} {[\,0\,,\, T\,]} {\L {4} {} \home} }$ solves the monodomain system in strong sense. Then, forming the inner product of (1.7) with the function $\psi = - \div \bigl(\,\condin\,\nabla \pottr\,\bigr) \in \L {2} {} \home$, we obtain
$$\eqalignno{
{} & \hbox{\kern-20pt} - \int_\Omega {\partial \pottr \over \partial t} \,\div \bigl(\,\condin\,\nabla \pottr\,\bigr)\, dx  - \int_\Omega \Bigl(\, b\,\pottr^3 - (a+1)\, b\,\pottr^2 + a\,b\, \pottr + \pottr\,\gat\,\Bigr)  \,\div \bigl(\,\condin\,\nabla \pottr\,\bigr)\, dx  & (2.32)
\cr
{} & + {\la \over 1 + \la}\, \int_\Omega \Bigl(\, \div \bigl(\,\condin\,\nabla \pottr\,\bigr) \,\Bigr)^2\, dx \,=\, \int_\Omega  \excinex \,\div \bigl(\,\condin\,\nabla \pottr\,\bigr)\, dx\quad \hbox{for a.~a.~} t \in (\, 0\,,\,T\,) \quad \iff
\cr
{} & \hbox{\kern-20pt} {1 \over 2}\,{\partial \over \partial t} \,\langle \,\nabla \pottr\,,\, \condin\,\nabla \pottr\,\rangle^2 + {\la \over 1 + \la}\,\langle \, \div \bigl(\,\condin\,\nabla \pottr\,\bigr) \,,\, \div \bigl(\,\condin\,\nabla \pottr\,\bigr) \,\rangle^2 + \langle \, 3\, b\,\pottr^2 & (2.33)
\cr
{} & - \,2\,(a+1)\, b\,\pottr + a\,b\, ,\, \nabla \pottr\, \condin\,\nabla \pottr\,\rangle \,=\, - \langle\,\excinex\,,\, \div \bigl(\,\condin\,\nabla \pottr\,\bigr)\,\rangle - \langle \, \pottr\,\gat \,,\, \div \bigl(\,\condin\,\nabla \pottr\,\bigr)\,\rangle \,.\phantom{\int}
\cr}$$
Completing the square on the left-hand side according to
$$\eqalignno{
{} & 3\, b\,  \varphi^2 - 2\,(a+1)\, b\, \varphi + a\,b \,=\, \bigl(\, \sqrt{\,3\,b\,}\,\varphi - \sqrt{\, b/3\,} \,(a+1)\,\bigr)^2  - b\,(a+1)^2 /3 + a\,b\,, & (2.34)
\cr}$$
thus yielding
$$\eqalignno{
{} & \hbox{\kern-20pt} {1 \over 2}\,{\partial \over \partial t} \,\bigl\vert\, {\cal M} ( \pottr, \pottr)\,\bigr\vert^2 + {\la \over 1 + \la}\,{\Norm{\div \bigl(\,\condin\,\nabla \pottr\,\bigr)}}^2_{\L {2} {} \home} + \langle \, \bigl(\, \sqrt{\,3\,b\,}\,\pottr - \sqrt{\, b/3\,} \,(a+1)\,\bigr)^2 \,,\, \nabla \pottr\, \condin\,\nabla \pottr\,\rangle
\cr
{} &  =\, - \langle\,\excinex\,,\, \div \bigl(\,\condin\,\nabla \pottr\,\bigr)\,\rangle - \langle \, \pottr\,\gat \,,\, \div \bigl(\,\condin\,\nabla \pottr\,\bigr)\,\rangle + \langle \, {b\,(a+1)^2\over 3} - a\,b\,,\, \nabla \pottr\, \condin\,\nabla \pottr\,\rangle\,, & (2.35)
\cr}$$
we arrive at the estimate
$$\eqalignno{
{} & \hbox{\kern-20pt} {1\over 2}\,{\partial \over \partial t} \,\bigl\vert\, {\cal M} ( \pottr, \pottr)\,\bigr\vert^2 + {\la \over 1 + \la}\,{\Norm{\div \bigl(\,\condin\,\nabla \pottr\,\bigr)}}^2_{\L {2} {} \home} \, & (2.36)
\cr
{} & & \displaystyle \le \,C\,\Bigl(\,1 + \bigl\vert\, {\cal M} (\pottr, \pottr)\,\bigr\vert^2\, +\, (\e_1 + \e_2)\,{\Norm{\div \bigl(\,\condin\,\nabla \pottr\,\bigr)}}^2_{\L {2} {} \home} + {1 \over 4\, \e_1} \,{\norm{\excinex}}^2_{\L {2} {} \home}  + {1 \over 4\,\e_2} \,{\Norm { \pottr\,\gat}}^2_{\L {2} {} \home}  \,\Bigr)
\cr}$$
which holds true for arbitrary $\e_1$, $\e_2 > 0$. Choosing for these parameters sufficiently small values, we find
$$\eqalignno{
{} & \hbox{\kern-20pt} {1\over 2}\,{\partial \over \partial t} \,\bigl\vert\, {\cal M} ( \pottr, \pottr)\,\bigr\vert^2 + {\la \over 2 + 2\,\la}\,{\Norm{\div \bigl(\,\condin\,\nabla \pottr\,\bigr)}}^2_{\L {2} {} \home} \, & (2.37)
\cr
{} & & \displaystyle \le \,C\,\Bigl(\,1 + \bigl\vert\, M(\pottr, \pottr)\,\bigr\vert^2  + {\norm{\excinex}}^2_{\L {2} {} \home}  +  {\Norm { \pottr\,\gat}}^2_{\L {2} {} \home}  \,\Bigr)\quad \follows
\cr
{} & \hbox{\kern-20pt} {\partial \over \partial t} \,\bigl\vert\, {\cal M} ( \pottr, \pottr)\,\bigr\vert^2 + {\Norm{\div \bigl(\,\condin\,\nabla \pottr\,\bigr)}}^2_{\L {2} {} \home} \, \cr
{} &  \le \,C\,\Bigl(\,1 + \bigl\vert\, {\cal M} ( \pottr, \pottr)\,\bigr\vert^2  + {\norm{\excinex}}^2_{\L {2} {} \home}  +  {\Norm { \pottr}}^4_{\L {4} {} \home}  +  {\Norm {\gat}}^4_{\L {4} {} \home}\,\Bigr) \phantom{\int} & (2.38)
\cr
{} & \le \,C\,\Bigl(\,1 + {\Norm { \pottr}}^2_{\sob {1,2} {} \home}   + {\Norm { \pottr}}^4_{\L {4} {} \home}  + {\norm{\excinex}}^2_{\L {2} {} \home}+ {\Norm {\gat}}^4_{\L {4} {} \home}\,\Bigr)\,, & (2.39)
\cr}$$
using the estimate for $\vert \, {\cal M}(\pottr, \pottr)\,\vert$ from $[\,${\kpt Kunisch/Wagner 12}$\,]\,$, p.~1529, Lemma 2.3. In view of (2.13) and Proposition 2.1., 2), this implies for any $t \in (\,0\,,\,T\,]$
$$\eqalignno{
{} & \hbox{\kern-20pt} \bigl\vert\, {\cal M} ( \,\pottr (t), \pottr(t)\,)\,\bigr\vert^2 \,\le\, C\,\int^t_0 \,\Bigl(\, 1 + {\Norm { \pottr}}^2_{\sob {1,2} {} \home}   + {\Norm { \pottr}}^4_{\L {4} {} \home}  + {\Norm{\excinex}}^2_{\L {2} {} \home} + {\Norm {\gat}}^4_{\L {4} {} \home}\,\Bigr)\, d\tau & (2.40)
\cr
{} & \le\, C\,\Bigl( \, 1 + {\Norm {\pottr}}^2_{ \bochnerl {2} {(\,0\,,\,T\,)} {\sob {1,2} {} \home} } + {\Norm { \pottr}}^4_{\L {4} {} (\Omega_T) } +  {\Norm {\gat}}^4_{\bochnerc {0} {[\,0\,,\,T\,]} {\L {4} {} \home} } & (2.41)
\cr
{} & & \displaystyle + \,{\Norm{\excin}}^2_{\L {2} {} (\Omega_T)} + {\Norm{\excex}}^2_{\L {2} {} (\Omega_T)} \,\Bigr)
\cr
{} & \le\, C\,\Bigl( \, 1 + {\Norm {\Phi_0}}^2_{\sob {1,2} {} \home } +  {\Norm {W_0}}^4_{ \L {4} {} \home } + {\Norm{\excin}}^2_{\L {2} {} (\Omega_T)} + {\Norm{\excex}}^2_{\L {2} {} (\Omega_T)} \,\Bigr)\,. & (2.42)
\cr}$$
Analogously, we get
$$\eqalignno{
{} & {\Norm{\div \bigl(\,\condin\,\nabla \pottr\,\bigr)}}^2_{\L {2} {} (\Omega_T)} \,\le\, C\,\int^T_0 \,\Bigl(\, 1 + {\Norm { \pottr}}^2_{\sob {1,2} {} \home}   + {\Norm { \pottr}}^4_{\L {4} {} \home}  + {\norm{\excinex}}^2_{\L {2} {} \home} + {\Norm {\gat}}^4_{\L {4} {} \home}\,\Bigr)\, d\tau
\cr
{} & \le \, C\,\Bigl( \, 1 + {\Norm {\Phi_0}}^2_{\sob {1,2} {} \home } +  {\Norm {W_0}}^4_{ \L {4} {} \home } + {\Norm{\excin}}^2_{\L {2} {} (\Omega_T)} + {\Norm{\excex}}^2_{\L {2} {} (\Omega_T)} \,\Bigr)\,. & (2.43)
\cr}$$
Since ${\cal M}$ is coercive, both together give
$$\eqalignno{
{} & {\Norm{\nabla \pottr}}^2_{\bochnerl {\infty} {(\,0\,,\,T\,)} {\L {2} {} \home} }+ {\Norm{\div \bigl(\,\condin\,\nabla \pottr\,\bigr)}}^2_{\L {2} {} (\Omega_T)} & (2.44)
\cr
{} & & \displaystyle \le \, C\,\Bigl( \, 1 + {\Norm {\Phi_0}}^2_{\sob {1,2} {} \home } +  {\Norm {W_0}}^4_{ \L {4} {} \home } + {\Norm{\excin}}^2_{\L {2} {} (\Omega_T)} + {\Norm{\excex}}^2_{\L {2} {} (\Omega_T)} \,\Bigr)\,.
\cr}$$
Observe now that $\bigl(\, {\norm{\phi}}^2_{\L {2} {}} + {\norm { \div \bigl(\,\condin\,\nabla \pottr\,\bigr)}}^2_{\L {2} {}} \,\bigr)^{1/2}$ is equivalent to the canonical norm within $\sob {2,2} {} \home$. Consequently, the last estimate implies
$$\eqalignno{
{} & {\Norm { \pottr}}^2_{\bochnerl {\infty} {(\,0\,,\,T\,)} {\sob {1,2} {} \home} } + {\Norm { \pottr}}^2_{\bochnerl {2} {(\,0\,,\,T\,)} {\sob {2,2} {} \home} } \, & (2.45)
\cr
{} & & \displaystyle \le\, C\,\Bigl( \, 1 + {\Norm {\Phi_0}}^2_{\sob {1,2} {} \home } +  {\Norm {W_0}}^4_{ \L {4} {} \home } + {\Norm{\excin}}^2_{\L {2} {} (\Omega_T)} + {\Norm{\excex}}^2_{\L {2} {} (\Omega_T)} \,\Bigr)\,.
\cr}$$
Combining this with (1.7) again, we may comfirm ourselves that
$$\eqalignno{
{} & {\Norm{{\partial \pottr \over \partial t} }}^{4/3}_{\L {4/3} {} (\Omega_T)} \,\le\, C\,\Bigl( \, 1 + {\Norm {\Phi_0}}^2_{\sob {1,2} {} \home } +  {\Norm {W_0}}^4_{ \L {4} {} \home } + {\Norm{\excin}}^2_{\L {2} {} (\Omega_T)} + {\Norm{\excex}}^2_{\L {2} {} (\Omega_T)} \,\Bigr) & (2.46)
\cr}$$
holds as well, and the a-priori estimate for the solution is established.
\smallskip
$\bullet$ {\bf Step 2.} {\it Existence of a strong solution.} Consider an orthonormal sequence $\folge{\psi_i}^\infty_{i=1}$ of eigenfunctions of the operator $-\div \bigl(\,\condin  \,\nabla\,(\,\cdot\,)\,\bigr)$ with Neumann boundary conditions. The aditional regularity assumptions about $\partial \Omega$ and $\condin$ imply that eigenfunctions $\psi \in \sob {2,2} {} \home$ may be chosen. Now, applying the standard Galerkin technique, we approximate $\pottr$ by
$$\eqalignno{
{} & \pottr^N (t) = \sum^N_{i=1}\, c^N_i (t)\,\psi_i (x) \,, & (2.47)
\cr}$$
obtaining $c^N_i (t)$ from the Galerkin approximation for (1.7) within the space $\V^N = \span \folge{\psi_i}^N_{i=1}$. The functions $\pottr^N $ still satisfy the uniform estimates from Step 1. Since the norms are related to reflexive Banach spaces, we can pass to the limit with standard subsequential arguments. $\blackbox$
\bigskip
{\bf d) Monodomain system with linearized Aliev-Panfilov model: improvement of regularity for the weak solutions.}
\medskip
{\bf Proposition 2.4.~(Gain of regularity for the gating variable)} {\it Consider the monodomain system with the linearized Aliev-Panfilov model under the assumptions of Theorem 1.1., and let a weak solution $(\pottr, \gat)$ of it correspond to inital values $\Phi_0\in \L {2} {} \home$, $W_0 \in \L {4} {} \home$ and right-hand sides $\excin$, $\excex \in \L {2} {} \bigl[\, {(\,0\,,\,T\,)}\,,$ $ {\bigl(\,\sob {1,2} {} \home\,\bigr)^*}\,\bigr]\,$, thus belonging to the spaces in $(2.12)$.
\smallskip
1) Then $\gat$ belongs to $\bochnerc {1} {(\,0\,,\,1\,)} {\L {1} {} \home} $, and $\partial \gat / \partial t$ belongs to $\L {2} {} (\Omega_T)$.
\smallskip
2) Moreover, $\,W_0 \in \L {4} {} \home$ implies that $\gat$ belongs even to $\bochnerc {0} {[\,0\,,\,T\,]} {\L {3} {} \home} $, and it holds that
$$\eqalignno{
{} & {\norm{\gat}}^3_{\bochnerc {0} {[\,0\,,\,T\,]} {\L {3} {} \home} }\,  \le\, C\,\Bigl(\,1 + {\norm{\Phi_0}}^2_{\L {2} {} \home} + {\norm{W_0}}^3_{\L {3} {} \home} & (2.48)
\cr
{} & & \displaystyle  + \, {\norm{\excin}}^2_{\bochnerl {2} {(\,0\,,\,T\,)} {\bigl(\,\sob {1,2} {} \home\,\bigr)^*} }  +  {\norm{\excex}}^2 _{\bochnerl {2} {(\,0\,,\,T\,)} {\bigl(\,\sob {1,2} {} \home\,\bigr)^*} }  \,\Bigr) \,.
\cr}$$}%
\par
Note that Part 1) of Proposition 2.4.~still holds true for $W_0 \in \L {2} {} \home$.
\medskip
{\bf Proposition 2.5.~(Gain of regularity for the transmembrane potential, strong solution of the monodomain system)} {\it Consider the monodomain system with the linearized Aliev-Panfilov model under the assumptions of Theorem 1.1. If, moreover, $\partial\Omega$ is of $\cont {1,1} {}$-regularity, the coefficients of $\condin$ belong to $\sob {1, \infty} {} \home$, $\Phi_0 \in \sob {1,2} {} \home$, $W_0 \in \L {4} {} \home$ and $\excin$, $\excex \in \L {2} {} (\Omega_T)$ then the system admits even a strong solution $(\pottr, \gat)$  with
$$\eqalignno{
{} & \pottr \,\in \, \bochnerc {0} {[\,0\,,\,T\,]} {\sob {1,2} {} \home } \,\cap \, \bochnerl {2} {(\,0\,,\,T\,)} {\sob {2,2} {} \home} \,\cap\,\bochnersob {1, 4/3} {(\,0\,,\,T\,)} {\L {4/3} { } \home} \,. & (2.49)
\cr}$$}%
\par
{\bf Proof of Proposition 2.4.} {\it Part 1)\/} Note first that, under the assumptions of the proposition, $\gat$ is represented as
$$\eqalignno{
{} & \gat (x,t) \,=\, W_0(x) \, e^{-\e t} + \e\,\kappa\, e ^{-\e t} \int^t_0 \Bigl(\,(a+1)\,\pottr(x,\tau) - \pottr^2 (x, \tau)\,\Bigr) \,e^{\e \tau}\, d\tau\,, & (2.50)
\cr}$$
and $(2.13)$ still holds true. By differentiation of $(2.50)$, we get
$$\eqalignno{
{} & {\partial \gat \over \partial t} (x,t) \,=\, - W_0(x) \, \e\,e^{-\e t} - \e\,\kappa\, e ^{-\e t} \int^t_0 \Bigl(\,(a+1)\,\pottr(x,\tau) - \pottr^2 (x, \tau)\,\Bigr) \,e^{\e \tau}\, d\tau\, & (2.51)
\cr
{} & & \displaystyle + \,\e\,\kappa\, \Bigl(\,(a+1)\,\pottr(x,t) - \pottr^2 (x, t)\,\Bigr) \quad \follows
\cr
{} & \hbox{\kern-20pt} \int_\Omega \bigl\vert \,{\partial \gat \over \partial t} (x,t)\,\bigr\vert \, dx\,\le \, C \, e^{-\e t}\,{\norm{W_0}}_{\L {1} {} \home} + C\, e ^{-\e t} \int^t_0 \int_\Omega \Bigl(\,\bigl\vert\, \pottr(x,\tau)\,\bigr\vert + \bigl\vert \,\pottr (x, \tau)\,\bigr\vert^2 \,\Bigr) \,e^{\e \tau}\, dx\, d\tau & (2.52)
\cr
{} & & \displaystyle  + \,C\, \Bigl(\,{\norm{\pottr(t)}}_{\L {1} {} \home} + {\norm{\pottr(t)}}^2 _{\L {2} {} \home}\,\Bigr)
\cr
{} & \hbox{\kern-20pt} \le \, C\, e^{-\e t}\,\Bigl(\, {\norm{W_0}}_{\L {1} {} \home} + \int^t_0  \Bigl(\,{\norm{\pottr(\tau)}}_{\L {1} {} \home}  + {\norm{\pottr (\tau)}}^2_{\L {2} {} \home} \,\Bigr) \,e^{\e \tau}\, d\tau \,\Bigr) & (2.53)
\cr
{} & & \displaystyle + \,C\, \Bigl(\,{\norm{\pottr(t)}}_{\L {2} {} \home} + {\norm{\pottr(t)}}^2 _{\L {2} {} \home}\,\Bigr)\,,
\cr}$$
and $\pottr \in \bochnerc {0} {[\,0\,,\,T\,]} {\L {2} {} \home} $ implies $\partial \gat / \partial t \in \bochnerc {0} {(\,0\,,\,T\,)} {\L {1} {} \home} $. Thus the claimed $\bochnerc {1} {(\,0\,,\,1\,)} {\L {1} {} \home} $-regularity of $\gat$ is proved.
\par
Using Lemma 2.3., let us estimate
$$\eqalignno{
{} & \hbox{\kern-20pt} \Bigl(\,\int_\Omega \bigl\vert \,{\partial \gat \over \partial t} (x,t)\,\bigr\vert^2 \, dx\,\Bigr)^{1/2} \,\le \, C\,{\norm{W_0}}_{\L {2} {} \home} +  C\, \Bigl(\, \int_\Omega \Bigl(\, \int^T_0 \bigl(\, \bigl\vert \pottr(x,\tau) \,\bigr\vert + \bigl\vert\, \pottr(x, \tau) \,\bigr\vert^2  \,\bigr) \, d\tau \,\Bigr)^2\, dx \,\Bigr)^{1/2} & (2.54)
\cr
{} & & \displaystyle + \,C\, \Bigl(\,\int_\Omega \Bigl(\,\bigl\vert \,\pottr(x,t)\,\bigr\vert^2 + \bigl\vert \, \pottr(x,t) \,\bigr\vert^{4} \,\Bigr)\, dx\,\Bigr)^{1/2}
\cr
{} & \hbox{\kern-20pt}\le \, C\, {\norm{W_0}}_{\L {2} {} \home} +  C\, \int^T_0 \Bigl(\, \int_\Omega \bigl(\, \bigl\vert \pottr(x,\tau) \,\bigr\vert^2 + \bigl\vert\, \pottr(x, \tau) \,\bigr\vert^4  \,\bigr) \, dx \,\Bigr)^{1/2} \, d\tau & (2.55)
\cr
{} & & \displaystyle + \,C\, \Bigl(\,\int_\Omega \Bigl(\,\bigl\vert \,\pottr(x,t)\,\bigr\vert^2 + \bigl\vert \, \pottr(x,t) \,\bigr\vert^{4} \,\Bigr)\, dx\,\Bigr)^{1/2}
\cr
{} & \hbox{\kern-20pt} \le \,C\, \Bigl(\, {\norm{W_0}}_{\L {2} {} \home} + {\norm{\pottr}}_{\bochnerl {1} {(\,0\,,\,T\,)} {\L {2} {} (\Omega)} } + {\norm{\pottr}}^2_{\bochnerl {2} {(\,0\,,\,T\,)} {\L {4} {} \home} } + {\norm{\pottr(t)}}_{\L {2} {} \home} + {\norm{\pottr(t)}}^{2}_{\L {4} {} \home}\,\Bigr)\,. &
\cr}$$
Consequently, we get\hfill(2.56)
$$\eqalignno{
{} & \bigl\Vert\,{\partial\gat \over \partial t}\,\bigr\Vert_{\L {2} {} (\Omega_T)} \,\le\, C\, \Bigl(\, {\norm{W_0}}_{\L {2} {} \home} + 2\,{\norm{\pottr}}_{\L {2} {} (\Omega_T)} + 2\,{\norm{\pottr}}^2_{\L {4} {} (\Omega_T) } \,\Bigr)\,, & (2.57)
\cr}$$
and by (2.13), the right-hand side is finite.
\medskip
{\it Part 2)\/} Recall now that $W_0 \in \L {4} {} \home$ (as we will see, it would suffice to assume that $W_0 \in \L{3} {} \home\,$). Then $(2.13)$ and Lemma 2.3.~imply
$$\eqalignno{
{} & \hbox{\kern-20pt} \Bigl(\,\int_\Omega \gat (t)^3 \, dx \,\Bigr)^{1/3} \,\le\, C\,{\norm{W_0}}_{\L {3} {} \home} + C\,\Bigl(\,\int_\Omega \,\Bigl(\, \int^T_0 \,\bigl(\, \bigl\vert \,\pottr (x, \tau)\,\bigr\vert + \bigl\vert\,\pottr (x, \tau)\,\bigr\vert^2\,\bigr)\, d\tau \,\Bigr)^3 \,\Bigr)^{1/3} & (2.58)
\cr
{} & \le\, C\, {\norm{W_0}}_{\L {3} {} \home} + C\,\int^T_0 \,\Bigl(\,\int_\Omega \,\bigl(\, \bigl\vert \pottr (x, \tau)\,\bigr\vert^3 + \bigl\vert\,\pottr (x, \tau)\,\bigr\vert^6\,\bigr)\, dx \,\Bigr)^{1/3} \,d\tau & (2.59)
\cr
{} & \le \,C \,\Bigl(\,  {\norm{W_0}}_{\L {3} {} \home} + {\norm{\pottr}}_{\bochnerl {1} {(\,0\,,\,T\,)} {\L {3} {} (\Omega)} } + {\norm{\pottr}}^2_{\bochnerl {2} {(\,0\,,\,T\,)} {\L {6} {} \home} } \,\Bigr) & (2.60)
\cr
{} & \le \,C \,\Bigl(\,  {\norm{W_0}}_{\L {3} {} \home} + {\norm{\pottr}}_{\L {3} {} (\Omega_T)} + {\norm{\pottr}}^2_{\bochnerl {2} {(\,0\,,\,T\,)} {\sob {1,2} {} \home} } \,\Bigr) \phantom{\int} & (2.61)
\cr
{} & \hbox{\kern-20pt} \le\, C\,\Bigl(\,1 + {\norm{\Phi_0}}^2_{\L {2} {} \home} + {\norm{W_0}}^3_{\L {3} {} \home} +  {\norm{\excin}}^2_{\bochnerl {2} {(\,0\,,\,T\,)} {\bigl(\,\sob {1,2} {} \home\,\bigr)^*} }  +  {\norm{\excex}}^2 _{\bochnerl {2} {(\,0\,,\,T\,)} {\bigl(\,\sob {1,2} {} \home\,\bigr)^*} }  \,\Bigr) \,, & (2.62)
\cr}$$
which confirms that $\gat$ belongs to $\bochnerl {\infty} {(\,0\,,\,T\,)} {\L {3} {} \home} $. In order to confirm that ${\norm{\gat (\,\cdot\,, t)}}^3_{\L {3} {} \home}$ depends even continuously on $t$, consider for $s$, $t \in [\,0\,,\,T\,]\,$
$$\eqalignno{
{} & \Bigl\vert\,{\norm{\gat(s)}}^3_{\L {3} {} \home} - {\norm{\gat (t)}}^3_{\L {3} {} \home} \,\Bigr\vert\,=\,\Bigl\vert \, \int_\Omega \Bigl(\,\gat(s)^3 - \gat (t)^3\,\Bigr)\, dx \,\Bigr\vert & (2.63)
\cr
{} & =\, \Bigl\vert \,\int_\Omega \bigl(\, \gat(s) - \gat (t)\,\bigr)\,\bigl(\,\gat (s)^2 + \gat(s)\,\gat (t) + \gat(t)^2\,\bigr)\, dx\,\Bigr\vert & (2.64)
\cr
{} & \le\, \Bigl(\, \int_\Omega \bigl\vert\, \gat(s) - \gat (t)\,\bigr\vert^3\, dx\,\Bigr)^{1/3} \,\Bigl(\, \int_\Omega \bigl\vert\, \gat (s)^2 + \gat(s)\,\gat (t) + \gat(t)^2\,\bigr\vert^{3/2} \, dx\,\Bigr)^{2/3}\,. & (2.65)
\cr}$$
Estimating the first term with Jensen's inequality and Lemma 2.3., we obtain
$$\eqalignno{
{} & \hbox{\kern-20pt} \Bigl(\,\int_\Omega \bigl\vert\, \gat(s) - \gat (t)\,\bigr\vert^3\, dx\,\Bigr)^{1/3} \,=\, \Bigl(\,\int_\Omega \,\Bigl\vert\, W_0\,\bigl(\, e^{-\e s} - e^{-\e t}\,\bigr) & (2.66)
\cr
{} & + \, \e\,\kappa\,(a+1)\,e^{-\e s}\,\bigl(\, \int^s_0 \pottr \,e^{\e\tau}\, d\tau - \int^t_0 \pottr \, e^{\e\tau}\,d\tau\,\bigr) \,+ \,\e\,\kappa\,(a+1) \,\bigl(\, e^{-\e s} - e^{-\e t} \,\bigr)\,\int^t_0 \pottr \,e^{\e\tau}\, d\tau\,
\cr
{} & - \, \e\,\kappa\,e^{-\e s}\,\bigl(\, \int^s_0 \pottr^2 \,e^{\e\tau}\, d\tau - \int^t_0 \pottr^2 \, e^{\e\tau}\,d\tau\,\bigr) \,- \,\e\,\kappa\,\,\bigl(\, e^{-\e s} - e^{-\e t} \,\bigr)\,\int^t_0 \pottr^2 \,e^{\e\tau}\, d\tau\,\Bigr\vert^3 \, dx \,\Bigr)^{1/3}
\cr
{} & \hbox{\kern-20pt}\le\, C\,\bigl\vert \, e^{-\e s} - e^{-\e t}\,\bigr\vert \cdot {\norm{W_0}}_{\L {3} {} \home} + C\,\Bigl(\,\int_\Omega \Bigl\vert\,\int^s_t \pottr\,e^{\e\tau}\, d\tau\,\Bigr\vert^3\, dx \,\Bigr)^{1/3}  + C\,\Bigl(\,\int_\Omega \Bigl\vert\,\int^s_t \pottr^2 \,e^{\e\tau}\, d\tau\,\Bigr\vert^3\, dx \,\Bigr)^{1/3} & (2.67)
\cr
{} & & \displaystyle +\, C\,\bigl\vert \, e^{-\e s} - e^{-\e t}\,\bigr\vert \cdot \Bigl(\, \int_\Omega \Bigl\vert\,\int^t_0 \pottr\,e^{\e\tau}\, d\tau\,\Bigr\vert^3\, dx \,\Bigr)^{1/3} + C\,\bigl\vert \, e^{-\e s} - e^{-\e t}\,\bigr\vert \cdot \Bigl(\, \int_\Omega \Bigl\vert\,\int^t_0 \pottr^2\,e^{\e\tau}\, d\tau\,\Bigr\vert^3\, dx \,\Bigr)^{1/3}
\cr
{} & \hbox{\kern-20pt} \le\, C\,\bigl\vert \, e^{-\e s} - e^{-\e t}\,\bigr\vert \cdot {\norm{W_0}}_{\L {3} {} \home} \phantom{\int} & (2.68)
\cr
{} &  + \,C\, \bigl\vert \,s - t\,\bigr\vert^{2/3}\,\Bigl(\,\Bigl(\, \int_\Omega \int^s_t \vert \,\pottr\,\vert^3 \,e^{ 3\, \e\, \tau}\, d\tau\, dx\,\Bigr)^{1/3} + \Bigl(\,\int_\Omega \int^s_t \vert \,\pottr\,\vert^6 \, e^{3\,\e\,\tau} \,d\tau \,dx \,\Bigr)^{1/3}\,\Bigr)
\cr
{} & & \displaystyle +  \,C\,\bigl\vert \, e^{-\e s} - e^{-\e t}\,\bigr\vert \cdot \Bigl(\,\int^t_0 \Bigl(\, \int_\Omega\,\vert \, \pottr \,\vert ^3\, e^{3\,\e \, \tau} dx\,\Bigr)^{1/3} \, d\tau + \int^t_0 \Bigl(\,\int_\Omega\,\vert \, \pottr \,\vert ^6\, e^{3\,\e \, \tau} dx\,\Bigr)^{1/3} \, d\tau\,\Bigr)
\cr
{} & \hbox{\kern-20pt}\le\, C\,\bigl\vert \, e^{-\e s} - e^{-\e t}\,\bigr\vert \cdot {\norm{W_0}}_{\L {3} {} \home} \phantom{\int} & (2.69)
\cr
{} &  + \,C\, \bigl\vert \,s - t\,\bigr\vert^{2/3}\,\Bigl(\,\Bigl(\, \int_\Omega \int^T_0 \vert \,\pottr\,\vert^3 \,e^{3\, \e\, T}\, d\tau\,dx\,\Bigr)^{1/3} + \Bigl(\,\int_\Omega \int^T_0 \vert \,\pottr\,\vert^6 \, e^{3\,\e\,T} \,d\tau \,dx \,\Bigr)^{1/3}\,\Bigr)
\cr
{} & & \displaystyle +  \,C\,\bigl\vert \, e^{-\e s} - e^{-\e t}\,\bigr\vert \cdot \Bigl(\,\int^T_0 \Bigl(\, \int_\Omega\,\vert \, \pottr \,\vert ^3\, e^{3\,\e \, T} \,dx\,\Bigr)^{1/3} \, d\tau + \int^T_0 \Bigl(\,\int_\Omega\,\vert \, \pottr \,\vert ^6\, e^{3\,\e \, T} \,dx\,\Bigr)^{1/3} \, d\tau\,\Bigr)\,.
\cr}$$
Applying Lemma 2.3.~again, we may continue
$$\eqalignno{
{} & \hbox{\kern-20pt} \lldots \le\, C\,\bigl\vert \, e^{-\e s} - e^{-\e t}\,\bigr\vert \cdot {\norm{W_0}}_{\L {3} {} \home} \phantom{\int} & (2.70)
\cr
{} &  + \,C\, \bigl\vert \,s - t\,\bigr\vert^{2/3}\,\Bigl(\,\int^T_0 \Bigl(\, \int_\Omega \vert \,\pottr\,\vert^3 \,dx\,\Bigr)^{1/3} \, d\tau + \int^T_0 \Bigl(\,\int_\Omega \vert \,\pottr\,\vert^6 \, dx \,\Bigr)^{1/3}\,d\tau\,\Bigr)
\cr
{} & & \displaystyle +  \,C\,\bigl\vert \, e^{-\e s} - e^{-\e t}\,\bigr\vert \cdot \Bigl(\,\int^T_0 \Bigl(\, \int_\Omega\,\vert \, \pottr \,\vert ^3\, dx\,\Bigr)^{1/3} \, d\tau + \int^T_0 \Bigl(\,\int_\Omega\,\vert \, \pottr \,\vert ^6\, dx\,\Bigr)^{1/3} \, d\tau\,\Bigr)
\cr
{} &  \le\,C\,\bigl\vert \, e^{-\e s} - e^{-\e t}\,\bigr\vert \cdot\Bigl(\, {\norm{W_0}}_{\L {3} {} \home} + {\norm{\pottr}}_{\bochnerl {1} {(\,0\,,\,T\,)} {\L {3} {} (\Omega_T)} } + {\norm{\pottr}}^2_{\bochnerl {2} {(\,0\,,\,T\,)} {\L {6} {} \home} } \,\Bigr) & (2.71)
\cr
{} & & \displaystyle + \,C\,\bigl\vert \, s - t\,\bigr\vert^{2/3} \cdot \Bigl(\,  {\norm{\pottr}}_{\L {3} {} (\Omega_T)} + {\norm{\pottr}}^2_{\bochnerl {2} {(\,0\,,\,T\,)} {\L {6} {} \home} }\,\Bigr) \,.
\cr}$$
Estimation of the second term yields:
$$\eqalignno{
{} &  \int_\Omega \bigl\vert\, \gat (s)^2 + \gat(s) \,\gat (t) + \gat(t)^2\,\bigr\vert^{3/2} \, dx \, \le \, C\,\int_\Omega \, \Bigl(\,\vert \,\gat(s)\,\vert^3 + \vert \,\gat(t)\,\vert^3 \,\Bigr)\, dx\,.& (2.72)
\cr}$$
Since $\gat \in \bochnerl {\infty} {(\,0\,,\,T\,)} {\L {3} {} \home} $, the second term remains uniformly bounded. Consequently, for every $\e > 0$ we find a $\delta(\e) > 0$ such that $\vert \, s - t \,\vert \le \delta(\e)$ implies $\bigl\vert\,{\norm{\gat(s)}}^3_{\L {3} {} \home} - {\norm{\gat (t)}}^3_{\L {3} {} \home} \,\bigr\vert \le \e$, and $\gat$ belongs to $\bochnerc {0} {[\,0\,,\,T\,]} {\L {3} {} \home} $. $\blackbox$
\medskip
{\bf Proof of Proposition 2.5.} With the following exception, the proof of Proposition 2.2.~may be repeated. In view of Proposition 2.4., 2), we estimate the last summand in (2.37) as follows:
$$\eqalignno{
{} & \hbox{\kern-20pt} {\Norm { \pottr\,\gat}}^2_{\L {2} {} \home} \,\le\, \Bigl(\,\int_\Omega \bigl\vert \,\pottr\,\bigr\vert^6 \, dx \,\Bigr)^{1/3} \,\Bigl(\,\int_\Omega \bigl\vert \,\gat \,\bigr\vert^3 \, dx\,\Bigr)^{2/3} \, =\, {\Norm { \pottr}}^2_{\L {6} {} \home} \cdot {\Norm { \gat}}^2_{\L {3} {} \home}  \, & (2.73)
\cr
{} & & \displaystyle \le \, {\Norm { \pottr}}^2_{\sob {1,2} {} \home} \cdot {\Norm { \gat}}^2_{\L {3} {} \home} \,\le \, {\Norm { \pottr}}^2_{\sob {1,2} {} \home} \cdot {\Norm { \gat}}^2_{\L {3} {} \home} \,.
\cr}$$
Consequently, we arrive at the estimate
$$\eqalignno{
{} & {\Norm{{\partial \pottr \over \partial t} }}^{4/3}_{\L {4/3} {} (\Omega_T)} \,\le\, C\,\Bigl( \, 1 + {\Norm {\Phi_0}}^2_{\sob {1,2} {} \home } +  {\Norm {W_0}}^2_{ \L {3} {} \home } + {\Norm{\excin}}^2_{\L {2} {} (\Omega_T)} + {\Norm{\excex}}^2_{\L {2} {} (\Omega_T)} \,\Bigr)\,, & (2.74)
\cr}$$
and the proof may be continued as above. $\blackbox$
\bigskip
{\bf e) Weak formulation of the bidomain system and known regularity of weak solutions.}
\medskip
The full bidomain system $(1.1)\,-\,(1.6)$ can be equivalently stated in a parabolic-elliptic form.$\,$\footnote{$^{18)}$}{$[\,${\ninekpt Kunisch/Wagner 13a}$\,]\,$, p.~954, (2.1)$\,-\,$(2.6).} To this, the following weak formulation corresponds.
$$\eqalignno{
{} &  \int_\Omega \Bigl(\,{\partial \pottr \over \partial t }\,\cdot \psi + \nabla \psi^{\T} \condin (\,\nabla \pottr + \nabla \potex\,) + \ioncurr(\pottr, \gat) \,\psi\,\Bigr)\, dx \,=\,\int_\Omega \excin \,\psi \,dx & (2.75)
\cr
{} & & \displaystyle \quad \forall\, \psi \in \sob {1,2} {} \home\,,\,\,\, \hbox{for a.~a.~$t \in (\,0\,,\,T\,)$} \,;
\cr
{} & \int_\Omega \Bigl(\,\nabla \psi^{\T} \condin \nabla \pottr + \nabla \psi^{\T} (\condin + \condex)\,\nabla \potex \,\Bigr)\, dx\,=\, \int_\Omega \Bigl(\,\excin + \excex \,\Bigr)\, \psi \, dx \quad & (2.76)
\cr
{} & & \displaystyle \forall \, \psi \in \sob {1,2} {} \home\,\,\,\hbox{with}\,\,\, \int_\Omega \psi(x)\, dx = 0\,,\,\,\, \hbox{for a.~a.~$t \in (\,0\,,\,T\,)$} \,;
\cr
{} & \int_\Omega \Bigl(\,{\partial \gat \over \partial t} + G(\pottr, \gat)\,\Bigr)\, \psi \, dx\,=\, 0 \quad \forall\, \psi \in \L {2} {} \home\,,\,\,\, \hbox{for a.~a.~$t \in (\,0\,,\,T\,)$} \,; & (2.77)
\cr
{} & \pottr(x,0) \,= \,\Phi_0(x)\quad \hbox{and}\quad  \gat (x,0) \,=\, W_0 (x) \quad \hbox{for almost all } x \in \Omega\,. \vphantom{\int} & (2.78)
\cr}$$
Under the assumptions of Theorem 1.2., the system $(2.75)\,-\,(2.78)$ with either the FitzHugh-Nagumo, the Rogers-McCulloch or the linearized Aliev-Panfilov model admits for arbitrary initial values $\Phi_0$, $W_0 \in \L {2} {} \home$ and inhomogeneities $\excin$, $\excex \in \bochnerl {2} {(\,0\,,\, T\,)} {\bigl(\, \sob {1,2} {} \home\,\bigr)^*} $, which satisfy the compatibility condition
$$\eqalignno{
{} & \int_\Omega \Bigl(\,\excin (x,t) + \excex (x,t)\,\Bigr)\, dx\,=\, 0 \quad \hbox{for a.~a.~} \, t \in (\,0\,,\, T\,)\,, & (2.79)
\cr}$$
at least one weak solution$\,$\footnote{$^{19)}$}{$[\,${\ninekpt Kunisch/Wagner 13a}$\,]\,$, p.~958, Theorem 2.5.}
$$\eqalignno{
{} & (\pottr, \potex, \gat)\,\in \,
\Bigl(\, \bochnerc {0} {[\,0\,,\, T\,]} {\L {2} {} \home} \,\cap\,\bochnerl {2} {(\,0\,,\,T\,)} {\sob {1,2} {} \home} \,\cap\, \L {4} {} (\Omega_T)\,\Bigr) & (2.80)
\cr
{} & & \displaystyle \times \, \bochnerl {2} {(\,0\,,\,T\,)} {\sob {2,2} {} \home} \,\times \, \bochnerc {0} {[\,0\,,\, T\,]} {\L {2} {} \home}
\cr}$$
with $\int_\Omega \potex(x,t)\, dx\,=\,0$ for almost all $t \in (\,0\,,\,T\,)$. Any weak solution obeys the a-priori estimate$\,$\footnote{$^{20)}$}{$[\,${\ninekpt Kunisch/Wagner 13a}$\,]\,$, p.~958, Theorem 2.6.}
$$\eqalignno{
{} & \hbox{\kern-20pt} {\norm{\pottr}}^2_{ \bochnerc {0} {[\,0\,,\, T\,)} {\L {2} {} \home} } + {\norm{\pottr}}^2_{ \bochnerl {2} {(\,0\,,\,T\,)} {\sob {1,2} {} \home} } + {\norm{\pottr}}^4_{\L {4} {} (\Omega_T) } +  {\norm{\partial \pottr/\partial t}}^{4/3}_{\bochnerl {4/3} {(\,0\,,\,T\,)} {\bigl(\,\sob {1,2} {} \home\,\bigr)^\ast} }
\cr
{} & + {\norm{\potex}}^2_{\bochnerl {2} {(\,0\,,\,T\,)} {\sob {1,2} {} \home} } + {\norm{\gat}}^2_{\bochnerc {0} {[\,0\,,\, T\,)} {\L {2} {} \home} }  + {\norm{\partial \gat /\partial t}}^2_{\bochnerl {2} {(\,0\,,\,T\,)} {\bigl(\,\sob {1,2} {} \home\,\bigr)^\ast} }
\cr}$$
$$\eqalignno{
{} & \hbox{\kern-20pt} \le \, C\cdot \Bigl(\,1 + {\norm{\Phi_0}}^2_{\L {2} {} \home} + {\norm{W_0}}^2_{\L {2} {} \home} + {\norm{\excin}}^2_{\bochnerl {2} {(\,0\,,\,T\,)} {\bigl(\,\sob {1,2} {} \home\,\bigr)^\ast} } + \, {\norm{\excex}}^2_{\bochnerl {2} {(\,0\,,\,T\,)} {\bigl(\,\sob {1,2} {} \home\,\bigr)^\ast} }  \,\Bigr) & (2.81)
\cr}$$
with a constant $C > 0$ not depending on $\Phi_0$, $W_0$, $\excin$ and $\excex$. It turns out that, under the assumptions of Theorem 1.2., a triple $(\pottr, \potex, \gat)$ forms a weak solution of the bidomain system $(2.75)\,-\,(2.78)$ iff the pair $(\pottr, \gat)$ solves the {\it reduced bidomain system\/}$\,$\footnote{$^{21)}$}{$[\,${\ninekpt Kunisch/Wagner 13a}$\,]\,$, p.~956 f., Theorem 2.4., 1).}
$$\eqalignno{
{} & \int_\Omega \Bigl(\,{\partial \pottr \over \partial t }+ \ioncurr (\pottr(t), \gat (t))\,\Bigr)\, \psi \, dx \,+ \, {\cal A} \bigl(\, \pottr(t)\,,\, \psi\,\bigr)  \,=\, \int_\Omega S(t)\, \psi\,dx \quad & (2.82)
\cr
{} & & \displaystyle \forall\, \psi \in \sob {1,2} {} \home\,\,\, \hbox{for a.~a.~$t \in (\,0\,,\,T\,)$} \,;
\cr
{} & \int_\Omega \Bigl(\,{\partial \gat \over \partial t} + G(\pottr, \gat)\,\Bigr)\, \psi \, dx\,=\, 0 \quad \forall\, \psi \in \L {2} {} \home\,\,\,\, \hbox{for a.~a.~$t \in (\,0\,,\,T\,)$} & (2.83)
\cr
{} & \pottr(x,0) \,= \,\Phi_0(x)\quad \hbox{and}\quad  \gat (x,0) \,=\, W_0 (x) \quad \hbox{for a.~a.~} x \in \Omega\,. \vphantom{\int} & (2.84)
\cr}$$
where ${\cal A}\,\colon \,\, \sob {1,2} {} \home \times \sob {1,2} {} \home \to \R {}$ is the {\it bidomain bilinear form}, and $S(t)$ is defined with the aid of $\excin$ and $\excex$, cf.~$[\,${\kpt Kunisch/Wagner 13a}$\,]\,$, p.~956 f., $(2.22)\,-\,(2.25)$. The structure of the reduced bidomain system is in complete analogy to those of the monodomain system, and the respective solutions obey the same type of a-priori estimates. Consequently, the regularity of $\pottr$ and $\gat$ within a weak solution $(\pottr, \potex, \gat)$ of $(2.75)\,-\,(2.78)$ can be improved in the same way as in Subsections 2.c) and d). The most important regularity property of $\pottr$, however, is stated again in Theorem 1.5., 2).
\bigskip
{\bf f) Bidomain system with Rogers-McCulloch model: improvement of regularity for the weak solutions.}
\medskip
{\bf Proposition 2.6.~(Gain of regularity for the gating variable)} {\it Consider the bidomain system $\,(2.75)\,-\,(2.79)$ with the Rogers-McCulloch model under the assumptions of Theorem 1.2., and let a weak solution $(\pottr, \potex, \gat)$ of it correspond to inital values $\Phi_0\in \L {2} {} \home$, $W_0 \in \L {4} {} \home$ and right-hand sides $\excin$, $\excex \in \L {2} {} \bigl[\, {(\,0\,,\,T\,)}\,,$ $ {\bigl(\,\sob {1,2} {} \home\,\bigr)^*}\,\bigr]\,$, thus belonging to the spaces in $(2.80)$.
\smallskip
1) Then $\gat$ belongs to $\bochnerc {1} {(\,0\,,\,T\,)} {\L {2} {} \home} $.
\smallskip
2) Moreover, $\,W_0 \in \L {4} {} \home$ implies that $\gat$ belongs even to $\bochnerc {0} {[\,0\,,\,T\,]} {\L {4} {} \home} $, and it holds that
$$\eqalignno{
{} & {\norm{\gat}}^4_{\bochnerc {0} {[\,0\,,\,T\,]} {\L {4} {} \home} }\,  \le\, C\,\Bigl(\,1 + {\norm{\Phi_0}}^2_{\L {2} {} \home} + {\norm{W_0}}^4_{\L {4} {} \home} & (2.85)
\cr
{} & & \displaystyle  + \, {\norm{\excin}}^2_{\bochnerl {2} {(\,0\,,\,T\,)} {\bigl(\,\sob {1,2} {} \home\,\bigr)^*} }  +  {\norm{\excex}}^2 _{\bochnerl {2} {(\,0\,,\,T\,)} {\bigl(\,\sob {1,2} {} \home\,\bigr)^*} }  \,\Bigr) \,.
\cr}$$}%
\par
{\bf Proof.} Since the proof of Proposition 2.1.~relies exclusively on the structure of the weak gating equation, which is the same in the monodomain and the reduced bidomain system, as well as on the a-priori estimate for $\gat$, we may carry over the argumentation without alterations. $\blackbox$
\medskip
Note that Part 1) of Proposition 2.6.~holds already true for $W_0 \in \L {2} {} \home$.
\bigskip
{\bf g) Bidomain system with linearized Aliev-Panfilov model: improvement of regularity for the weak solutions.}
\medskip
{\bf Proposition 2.7.~(Gain of regularity for the gating variable)} {\it Consider the bidomain system $\,(2.75)\,-\,(2.79)$ with the linearized Aliev-Panfilov model under the assumptions of Theorem 1.2., and let a weak solution $(\pottr, \potex, \gat)$ of it correspond to inital values $\Phi_0\in \L {2} {} \home$, $W_0 \in \L {4} {} \home$ and right-hand sides $\excin$, $\excex \in \L {2} {} \bigl[\, {(\,0\,,\,T\,)}\,,$ $ {\bigl(\,\sob {1,2} {} \home\,\bigr)^*}\,\bigr]\,$, thus belonging to the spaces in $(2.80)$.
\smallskip
1) Then $\gat$ belongs to $\bochnerc {1} {(\,0\,,\,1\,)} {\L {1} {} \home} $.
\smallskip
2) Moreover, $\,W_0 \in \L {4} {} \home$ implies that $\gat$ belongs even to $\bochnerc {0} {[\,0\,,\,T\,]} {\L {3} {} \home} $, and it holds that
$$\eqalignno{
{} & {\norm{\gat}}^3_{\bochnerc {0} {[\,0\,,\,T\,]} {\L {3} {} \home} }\,  \le\, C\,\Bigl(\,1 + {\norm{\Phi_0}}^2_{\L {2} {} \home} + {\norm{W_0}}^3_{\L {3} {} \home} & (2.86)
\cr
{} & & \displaystyle  + \, {\norm{\excin}}^2_{\bochnerl {2} {(\,0\,,\,T\,)} {\bigl(\,\sob {1,2} {} \home\,\bigr)^*} }  +  {\norm{\excex}}^2 _{\bochnerl {2} {(\,0\,,\,T\,)} {\bigl(\,\sob {1,2} {} \home\,\bigr)^*} }  \,\Bigr) \,.
\cr}$$}%
\medskip
{\bf Proof.} The proof of Proposition 2.4.~may be repeated without changes. $\blackbox$
\medskip
Even here, Part 1) of Proposition 2.8.~holds true for $W_0 \in \L {2} {} \home$.
\bigskip
{\bf h) Proof of Theorem 1.5.}
\medskip
{\it Part 1)\/} Within the monodomain system, we specify the Rogers-McCulloch model.
\smallskip
$\bullet$ {\bf Step 1.} {\it An estimate for $\bigl\vert \, {1 \over 2} \,{d \over dt}\, {\norm{\pottr (t)}}^2_{\L {2} {} \home} \,\bigr\vert$.} Inserting into equation (2.9) the feasible test function $\psi = \pottr (t) \in \sob {1,2} {} \home$ and applying the lower estimate for the monodomain form, we obtain
$$\eqalignno{
{} & \hbox{\kern-20pt} \langle\,{d \over dt} \,\pottr(t) \,,\, \pottr(t)\,\rangle + {\cal M} \bigl(\, \pottr(t) \,,\, \pottr(t) \,\bigr) \,=\, - \int_\Omega \ioncurr (\pottr(t), \gat (t)) \,\pottr(t) \, dx & (2.87)
\cr
{} & & \displaystyle +\, \langle \, {1 \over 1 + \la}\, \Bigl(\,\la\,\bigl(\, \excin(t) - \excex (t) \,\Bigr) \,,\,\pottr(t)\,\rangle \quad \follows
\cr
{} & \hbox{\kern-20pt}{1 \over 2} \,{d \over dt}\, {\norm{\pottr}}^2_{\L {2} {} \home} + \beta \,{\norm{\pottr}}^2_{\sob {1,2} {} \home} \,\le \, -  \int_\Omega \Bigl(\,b\,(\pottr)^4 - (a+1)\,b\, (\pottr)^3 + a\,b\, (\pottr)^2 + (\pottr)^2\, \gat\,\Bigr) \, dx & (2.88)
\cr
{} & & \displaystyle \,+  \,\bigl\vert \, \langle\, {1 \over 1 + \la}\, \Bigl(\,\la\, \excin - \excex \,\Bigr) \,,\,\pottr \,\rangle \,\bigr\vert  + \beta \,{\norm{\pottr}}^2_{\L {2} {} \home} \quad \follows
\cr
{} & \hbox{\kern-20pt}{1 \over 2} \,{d \over dt}\, {\norm{\pottr}}^2_{\L {2} {} \home} + \beta \,{\norm{\pottr}}^2_{\sob {1,2} {} \home} \, & (2.89)
\cr
{} & & \displaystyle \le \, C\,\int_\Omega \bigl\vert \,\pottr\,\bigr\vert^2\,\bigl(\,\bigl\vert \,\pottr\,\bigr\vert + \bigl\vert\,\gat\,\bigr\vert\,\bigr)\, dx
\,+  \,\bigl\vert \, \langle\, {1 \over 1 + \la}\, \Bigl(\,\la\, \excin - \excex \,\Bigr) \,,\,\pottr \,\rangle \,\bigr\vert + \beta \,{\norm{\pottr}}^2_{\L {2} {} \home} \,.
\cr}$$
In order to estimate the first term on the right-hand side, we apply the generalized Cauchy inequality, thus obtaining
$$\eqalignno{
{} &  C\, \int_\Omega \bigl\vert \,\pottr \,\bigr\vert^2\,\Bigl(\, \bigl\vert \,\pottr\,\bigr\vert + \bigl\vert \,\gat\,\bigr\vert \,\Bigr) \, dx \, \le\,C\,\e_1(t)\, {\norm{\pottr}}^4_{\L {4} {} \home} + {C \over \e_1(t)} \, \Bigl(\, {\norm {\pottr}}^2_{\L {2} {} \home} + {\norm{\gat}}^2_{\L {2} {} \home}\,\Bigr) & (2.90)
\cr}$$
for arbitrary $\e_1(t) > 0$. Specifying within (2.90) $\e_1(t) = \e^\p_1 / \bigl(\,1 + {\norm{\pottr (t)}}^2_{\L {4} {} \home}\,\bigr)$, we may continue
$$\eqalignno{
{} & \hbox{\kern-20pt} C\, \int_\Omega \bigl\vert \,\pottr \,\bigr\vert^2\,\Bigl(\, \bigl\vert \,\pottr\,\bigr\vert + \bigl\vert \,\gat\,\bigr\vert \,\Bigr) \, dx \, & (2.91)
\cr
{} & \le\,C\,\e^\p_1   \,{\norm{\pottr}}^2_{\L {4} {} \home} + {C \over \e^\p_1 } \,\bigl(\,1 + {\norm{\pottr}}^2_{\L {4} {} \home}\,\bigr) \, \Bigl(\, {\norm {\pottr}}^2_{\L {2} {} \home} + {\norm{\gat}}^2_{\L {2} {} \home} \,\Bigr) &
\cr
{} & \hbox{\kern-20pt}\le\,C\,\e^\p_1 \,{\norm{\pottr}}^2_{\sob {1,2} {} \home} + {C \over \e^\p_1 } \,\bigl(\,1 + {\norm{\pottr}}^2_{\L {4} {} \home} \,\bigr) \, \Bigl(\, {\norm {\pottr}}^2_{\bochnerc {0} {[\,0\,,\,T\,]} {\L {2} {} \home} }+ {\norm{\gat}}^2_{\bochnerc {0} {[\,0\,,\,T\,]} {\L {2} {} \home} }\,\Bigr) & (2.92)
\cr}$$
$$\eqalignno{
{} & \hbox{\kern-20pt}\le \,C\,\e^\p_1  \,{\norm{\pottr}}^2_{\sob {1,2} {} \home} & (2.93)
\cr
{} & & \displaystyle + \,{C \over \e^\p_1 } \,\bigl(\,1 + {\norm{\pottr}}^2_{\L {4} {} \home}\,\bigr) \cdot \Bigl(\, 1 + {\norm{\Phi_0}}^2_{\L {2} {} \home } + {\norm{W_0}}^2_{\L {2} {} \home } + {\norm{\excin}}^2_{ \L {2} {} (\Omega_T) } +  {\norm{\excex}}^2_{ \L {2} {} (\Omega_T)} \,\Bigr)
\cr
{} & \hbox{\kern-20pt}\le \,C\,\e^\p_1  \,{\norm{\pottr}}^2_{\sob {1,2} {} \home} + {C \over \e^\p_1 } \,\bigl(\,1 + {\norm{\pottr}}^2_{\L {4} {} \home}\,\bigr) \cdot \Bigl(\, 1 + {\norm{\Phi_0}}^2_{\L {2} {} \home } + {\norm{W_0}}^2_{\L {2} {} \home } + 2\,R^2\,\Bigr) & (2.94)
\cr}$$
where (2.13) and the additional regularity of $\excin$ and $\excex$ have been employed. The generalized Cauchy inequality will be applied to the second term as well, thus getting
$$\eqalignno{
{} &  \bigl\vert \, \langle\, {1 \over 1 + \la}\, \bigl(\,\la\, \excin
 - \excex \,\bigr) \,,\,\pottr \,\rangle \,\bigr\vert \,\le\, {\norm{\, {1 \over 1 + \la}\,\bigl(\,\la \,\excin  - \excex \,\bigr)}}{}_{\bigl(\,\sob {1,2} {} \home\,\bigr)^*} \, {\norm{\pottr }}_{\sob {1,2} {} \home} & (2.95)
\cr
{} & & \displaystyle \le\, C\,\e^\p_2 \,{\norm{\pottr}}^2_{\sob {1,2} {} \home} + {C \over \e^\p_2}\, \Bigl(\, {\norm{\excin }}^2_{\L {2} {} \home } + {\norm{\excex }}^2_{\L {2} {} \home } \,\Bigr)
\cr}$$
with arbitrary $\e^\p_2 > 0$. The third term is estimated again by
$$\eqalignno{
{} & \beta \,{\norm{\pottr}}^2_{\L {2} {} \home} \,\le\, C \,{\norm{\pottr}}^2_{\bochnerc {0} {[\,0\,,\,T\,]} {\L {2} {} \home} } \, & (2.96)
\cr
{} & \le\, C \,\Bigl(\, 1 + {\norm{\Phi_0}}^2_{\L {2} {} \home } + {\norm{W_0}}^2_{\L {2} {} \home } + {\norm{\excin}}^2_{ \L {2} {} (\Omega_T) } +  {\norm{\excex}}^2_{ \L {2} {} (\Omega_T)} \,\Bigr)\, & (2.97)
\cr
{} & \le\, C \,\Bigl(\, 1 + {\norm{\Phi_0}}^2_{\L {2} {} \home } + {\norm{W_0}}^2_{\L {2} {} \home } + 2\,R^2\,\Bigr)\,, & (2.98)
\cr}$$
using the bound $R> 0$ for the norms of $\excin$ and $\excex$. Combining (2.94), (2.95) and (2.98), we arrive at
$$\eqalignno{
{} & \hbox{\kern-20pt}{1 \over 2} \,{d \over dt}\, {\norm{\pottr}}^2_{\L {2} {} \home} +\beta\, {\norm{\pottr}}^2_{\sob {1,2} {} \home} \, & (2.99)
\cr
{} & \le \, C\,\e^\p_1 \,{\norm{\pottr}}^2_{\sob {1,2} {} \home} \,+ \,{C \over \e^\p_1 } \,\bigl(\,1 + {\norm{\pottr}}^2_{\L {4} {} \home}\,\bigr) \, \Bigl(\, 1 + {\norm{\Phi_0}}^2_{\L {2} {} \home }  +\, {\norm{W_0}}^2_{\L {2} {} \home } + 2\,R^2 \,\Bigr)
\cr
{} & + \,C\,\e^\p_2 \,{\norm{\pottr}}^2_{\sob {1,2} {} \home} + \,{C \over \e^\p_2}\, \Bigl(\, {\norm{\excin }}^2_{\L {2} {} \home } + {\norm{\excex }}^2_{\L {2} {} \home } \,\Bigr) +  C \,\Bigl(\, 1 + {\norm{\Phi_0}}^2_{\L {2} {} \home } + {\norm{W_0}}^2_{\L {2} {} \home } + 2\,R^2\,\Bigr)\,.
\cr}$$
Observing that $1 + {\norm{\Phi_0}}^2_{\L {2} {} \home } + {\norm{W_0}}^2_{\L {2} {} \home } + 2\,R^2 \le C$, we find
$$\eqalignno{
{} & {1 \over 2} \,{d \over dt}\, {\norm{\pottr}}^2_{\L {2} {} \home} + \beta\, {\norm{\pottr}}^2_{\sob {1,2} {} \home} \,\le \, C + C\,\bigl(\,\e^\p_1 + \e^\p_2\,\bigr) \,{\norm{\pottr}}^2_{\sob {1,2} {} \home} \, & (2.100)
\cr
{} & & \displaystyle + \,{C \over \e^\p_1 } \,\bigl(\,1 + {\norm{\pottr}}^2_{\L {4} {} \home}\,\bigr) \,+ \,{C \over \e^\p_2}\, \Bigl(\, {\norm{\excin }}^2_{\L {2} {} \home } + {\norm{\excex }}^2_{\L {2} {} \home } \,\Bigr) \,.
\cr}$$
Choosing now $\e^\p_1$, $\e^\p_2 > 0$ in such a way that the terms with ${\norm{\pottr}}^2_{\sob {1,2} {} \home}$ on both sides of (2.100) annihilate, we obtain the inequality
$$\eqalignno{
{} & {1 \over 2} \,{d \over dt}\, {\norm{\pottr (t( }}^2_{\L {2} {} \home} \,\le \, C \,\bigl(\,1 + {\norm{\pottr(t)}}^2_{\L {4} {} \home}\,\bigr) + C \, \Bigl(\, {\norm{\excin (t)}}^2_{\L {2} {} \home } + {\norm{\excex(t) }}^2_{\L {2} {} \home } \,\Bigr) \,. & (2.101)
\cr}$$
Inserting now the reverse {\it test function\/} $\psi = - \pottr(t)$ into (2.9), we get instead
$$\eqalignno{
{} & \hbox{\kern-20pt} - \langle\,{d \over dt} \,\pottr(t) \,,\, \pottr(t)\,\rangle - {\cal M} \bigl(\, \pottr(t) \,,\, \pottr(t) \,\bigr) \,=\, \int_\Omega \ioncurr (\pottr(t), \gat (t)) \,\pottr(t) \, dx & (2.102)
\cr
{} & & \displaystyle -\, \langle \, {1 \over 1 + \la}\, \Bigl(\,\la\,\bigl(\, \excin(t) - \excex (t) \,\Bigr) \,,\,\pottr(t)\,\rangle \quad \follows
\cr
{} & \hbox{\kern-20pt} - {1 \over 2} \,{d \over dt}\, {\norm{\pottr}}^2_{\L {2} {} \home} - \beta \,{\norm{\pottr}}^2_{\sob {1,2} {} \home} \,\ge \, \int_\Omega \Bigl(\,b\,(\pottr)^4 - (a+1)\,b\, (\pottr)^3 + a\,b\, (\pottr)^2 + (\pottr)^2\, \gat\,\Bigr) \, dx & (2.103)
\cr
{} & & \displaystyle \,-  \,\bigl\vert \, \langle\, {1 \over 1 + \la}\, \Bigl(\,\la\, \excin - \excex \,\Bigr) \,,\,\pottr \,\rangle \,\bigr\vert  - \beta \,{\norm{\pottr}}^2_{\L {2} {} \home} \quad \follows
\cr
{} & \hbox{\kern-20pt} - {1 \over 2} \,{d \over dt}\, {\norm{\pottr}}^2_{\L {2} {} \home} - \beta \,{\norm{\pottr}}^2_{\sob {1,2} {} \home} \, & (2.104)
\cr
{} & & \displaystyle \ge \, - C\,\int_\Omega \bigl\vert \,\pottr\,\bigr\vert^2\,\bigl(\,\bigl\vert \,\pottr\,\bigr\vert + \bigl\vert\,\gat\,\bigr\vert\,\bigr)\, dx
\,-  \,\bigl\vert \, \langle\, {1 \over 1 + \la}\, \Bigl(\,\la\, \excin - \excex \,\Bigr) \,,\,\pottr \,\rangle \,\bigr\vert - \beta \,{\norm{\pottr}}^2_{\L {2} {} \home} \,.
\cr}$$
Using again (2.94), (2.95) and (2.98), (2.104) implies the reverse inequality
$$\eqalignno{
{} & \hbox{\kern-20pt}- {1 \over 2} \,{d \over dt}\, {\norm{\pottr}}^2_{\L {2} {} \home} - \beta\, {\norm{\pottr}}^2_{\sob {1,2} {} \home} \, & (2.105)
\cr
{} & \ge \, - C\,\e^\p_1 \,{\norm{\pottr}}^2_{\sob {1,2} {} \home} \,- \,{C \over \e^\p_1 } \,\bigl(\,1 + {\norm{\pottr}}^2_{\L {4} {} \home}\,\bigr) \, \Bigl(\, 1 + {\norm{\Phi_0}}^2_{\L {2} {} \home } +\, {\norm{W_0}}^2_{\L {2} {} \home } + 2\,R^2\,\Bigr)
\cr
{} & - \,C\,\e^\p_2 \,{\norm{\pottr}}^2_{\sob {1,2} {} \home} - \,{C \over \e^\p_2}\, \Bigl(\, {\norm{\excin }}^2_{\L {2} {} \home } + {\norm{\excex }}^2_{\L {2} {} \home } \,\Bigr) -\,  C \,\Bigl(\, 1 + {\norm{\Phi_0}}^2_{\L {2} {} \home } + {\norm{W_0}}^2_{\L {2} {} \home } + 2\,R^2\,\Bigr)\,,
\cr}$$
and we may choose again $\e^\p_1$, $\e^\p_2 > 0$ in such a way that the summands with $- {\norm{\pottr}}^2_{\sob {1,2} {} \home}$ annihilate. Thus (2.101) is reversed as
$$\eqalignno{
{} &  -{1 \over 2} \,{d \over dt}\, {\norm{\pottr (t)}}^2_{\L {2} {} \home} \,\ge \,- C \,\bigl(\,1 + {\norm{\pottr (t)}}^2_{\L {4} {} \home}\,\bigr) \, - \,C \, \Bigl(\, {\norm{\excin(t) }}^2_{\L {2} {} \home } + {\norm{\excex (t)}}^2_{\L {2} {} \home } \,\Bigr) \,, & (2.106)
\cr}$$
and we arrive at the desired estimate
$$\eqalignno{
{} & \Bigl\vert\,{1 \over 2} \,{d \over dt}\, {\norm{\pottr(t)}}^2_{\L {2} {} \home} \,\Bigr\vert\,\le \, C \,\bigl(\,1 + {\norm{\pottr(t)}}^2_{\L {4} {} \home}\,\bigr) + C \, \Bigl(\, {\norm{\excin(t) }}^2_{\L {2} {} \home } + {\norm{\excex (t)}}^2_{\L {2} {} \home } \,\Bigr) \,. & (2.107)
\cr}$$
\par
$\bullet$ {\bf Step 2.} {\it An estimate for ${\norm{\pottr}}^4_{\bochnerl {4} {(\,0\,,\,T\,)} {\sob {1,2} {} \home} }$.} We return to (2.100) and choose now $\e^\p_1$, $\e^\p_2 > 0$ in such a way that $C\,\e^\p_1 + C\,\e^\p_2 = \beta/2$. Then (2.100) and (2.107) imply
$$\eqalignno{
{} & \hbox{\kern-20pt}{ \beta\over 2} \, {\norm{\pottr}}^2_{\sob {1,2} {} \home} \,\le \,- {1 \over 2} \,{d \over dt}\, {\norm{\pottr}}^2_{\L {2} {} \home} + C \,+ \,C \,\bigl(\,1 + {\norm{\pottr}}^2_{\L {4} {} \home}\,\bigr) \, + \,C \, \Bigl(\, {\norm{\excin }}^2_{\L {2} {} \home } + {\norm{\excex }}^2_{\L {2} {} \home } \,\Bigr) & (2.108)
\cr
{} & \le\,\Bigl\vert \,{1 \over 2} \,{d \over dt}\, {\norm{\pottr}}^2_{\L {2} {} \home} \,\Bigr\vert \,+ \,C \,\bigl(\,1 + {\norm{\pottr}}^2_{\L {4} {} \home}\,\bigr) \, + \,C\, \Bigl(\, {\norm{\excin }}^2_{\L {2} {} \home } + {\norm{\excex }}^2_{\L {2} {} \home } \,\Bigr) & (2.109)
\cr
{} & \le\, 2\,C \,\bigl(\,1 + {\norm{\pottr}}^2_{\L {4} {} \home}\,\bigr) \, + 2\,C\, \Bigl(\, {\norm{\excin }}^2_{\L {2} {} \home } + {\norm{\excex }}^2_{\L {2} {} \home } \,\Bigr)\,. & (2.110)
\cr}$$
Consequently, we find
$$\eqalignno{
{} & {\norm{\pottr}}^4_{\bochnerl {4} {(\,0\,,\,T\,)} {\sob {1,2} {} \home} } \,=\, \int^T_0 {\norm{\pottr (t)}}^4_{\sob {1,2} {} \home} \, dt & (2.111)
\cr
{} & \le \,C \,\int^T_0 \Bigl(\,1 + 2\, {\norm{\pottr}}^2_{\L {4} {} \home} + {\norm{\pottr}}^4_{\L {4} {} \home}\,\Bigr)\, dt \, + \,C\, \int^T_0 \Bigl(\, {\norm{\excin }}^4_{\L {2} {} \home } + {\norm{\excex }}^4_{\L {2} {} \home } \,\Bigr)\, dt\, & (2.112)
\cr
{} & \le\,C \,\Bigl(\,1 + {\norm{\pottr}}^4_{\L {4} {} (\Omega_T) } + {\norm{\excin}}^4_{\bochnerl {4} {(\,0\,,\,T\,)} {\L {2} {} \home} } + {\norm{\excex }}^4_{\bochnerl {4} {(\,0\,,\,T\,)} {\L {2} {} \home} } \,\Bigr)\, , & (2.113)
\cr}$$
and the right-hand side is bounded by assumption about $\excin$, $\excex$ and (2.13).
\medskip
$\bullet$ {\bf Step 3.} {\it The other ionic models.} If the Rogers-McCulloch model is replaced by the FitzHugh-Nagumo or the linearized Aliev-Panfilov model, all arguments may be repeated since (2.13) holds still true, $\ioncurr$ admits the same or an analogous structure, and the gating equation is not involved.
\medskip
{\it Part 2)\/} We consider the bidomain system with either the Rogers-McCulloch, the FitzHugh-Nagumo or the linearized Aliev-Panfilov model and rely on the a-priori estimate (2.81) and the reduced system $(2.82)\,-\,(2.84)$ instead of (2.13) and $(2.9)\,-\,(2.11)$. Since both systems and estimates admit an identical structure, we may repeat all arguments from Part 1 provided that an analogue of (2.95) can be proven. Indeed, let us recall that the linear functionals $S(t) \,\colon \,\, \sob {1,2} {} \home \to \R {}$ in (2.82) are defined through$\,$\footnote{$^{22)}$}{\rm Cf.~$[\,${\ninekpt Bourgault/Coudi\`ere/Pierre 09}$\,]\,$, p.~464, Definition 5.}
$$\eqalignno{
{} & \langle\, S(t)\,,\,\psi \,\rangle \,=\, \langle \, \excin(t)\,,\, \psi\,\rangle - \int_\Omega \nabla {\overline {\psi}}^{\T}_e \condin \,\nabla \psi\, dx \quad \forall\, \psi \in \sob {1,2} {} \home & (2.114)
\cr}$$
where ${\overline{\psi}}_e \in \sob {1,2} {} \home$ is the uniquely determined solution of the variational equation
$$\eqalignno{
{} & \int_\Omega \nabla {\overline{\psi}}^{\T}_e (\condin + \condex) \,\nabla \varphi \, dx \, =\, \langle \, \excin(t) + \excex(t)\,,\, \varphi\,\rangle \quad \forall \,\varphi \in \sob {1,2} {} \home \,,& (2.115)
\cr}$$
which satisfies $\int_\Omega  {\overline{\psi}}_e (x,t)\, dx = 0$ $\fforall\,t \in (\,0\,,\,T\,)$. Repeating now the arguments from $[\,${\kpt Kunisch/Wagner 13a}$\,]\,$, p.~960 f., Proof of Lemma 2.9., we insert the solution $ {\overline{\psi}}_e \in \sob {1,2} {} \home$ itself as a test function into (2.115) and find by application of the Poincar\'e inequality$\,$\footnote{$^{23)}$}{Together with the Rellich-Kondrachov theorem, $[\,${\ninekpt Evans 98}$\,]\,$, p.~275, Theorem 1, holds true even on a bounded strongly Lipschitz domain, cf.~$[\,${\ninekpt Adams/Fournier 07}$\,]\,$, p.~168, Theorem 6.3. Note that ${\overline{\psi}}_e$ admits a zero spatial mean.} and the generalized Cauchy inequality \hfill (2.116)
$$\eqalignno{
{} & \hbox{\kern-20pt} C\,{\norm{{\overline{\psi}}_e }}^2_{\sob {1,2} {} \home} \,\le\, \int_\Omega \nabla {\overline{\psi}}^{\T}_e (\condin + \condex) \,\nabla {\overline{\psi}}_e  \, dx \, = \, \langle \, \excin(t) + \excex(t)\,,\, {\overline{\psi}}_e \,\rangle \,\le\, \bigl\vert \,\langle \, \excin(t)\,,\, {\overline{\psi}}_e \,\rangle \,\bigr\vert \,+ \,\bigl\vert \,\langle \, \excex(t)\,,\, {\overline{\psi}}_e \,\rangle \,\bigr\vert
\cr
{} & \le \,{1 \over 2\,\tilde\e_1} \,\Bigl(\, {\norm{ \excin(t)}}^2_{\bigl(\, \sob {1,2} {} \home \,\bigr)^*}  + {\norm{ \excex(t)}}^2_{\bigl(\, \sob {1,2} {} \home\,\bigr)^*} \,\Bigr) + {\tilde \e_1 \over 2} \,{\norm{ {\overline{\psi}}_e }}^2_{\sob {1,2} {} \home}\, & (2.117)
\cr}$$
for arbitrary $\tilde \e_1 > 0$. Using the continuous imbedding $\L {2} {} \home \hookrightarrow \bigl(\, \sob {1,2} {} \home\,\bigr)^*$ and specifying $\tilde\e_1 = C/2$, we arrive after a normalization of the constants at
$$\eqalignno{
{} & {\norm{{\overline{\psi}}_e }}^2_{\sob {1,2} {} \home} \,\le\, C \,\Bigl(\, {\norm{ \excin(t)}}^2_{\L {2} {} \home }  + {\norm{ \excex(t)}}^2_{ \L {2} {} \home} \,\Bigr) \,. & (2.118)
\cr}$$
Going back to (2.114) and inserting the feasible test function $\psi = \pottr(t)$ into this equation, we estimate
$$\eqalignno{
{} & \bigl\vert \,\langle \, S(t)\,,\, \pottr\,\rangle\,\bigr\vert \,\le\, \bigl\vert \,\langle \, \excin(t)\,,\, \pottr\,\rangle\,\bigr\vert + \bigl\vert \,\langle \, \nabla {\overline {\psi}}^{\T}_e \,,\, \condin \,\nabla \pottr\,\rangle\,\bigr\vert \,& (2.119)
\cr
{} & \le \,{1 \over 2\,\tilde \e_2} \, {\norm{\excin(t)}}^2_{\L {2} {} \home} + {\tilde \e_2 \over 2}\, {\norm{\pottr}}^2_{\sob {1,2} {} \home} + {1 \over 2\,\tilde \e_3} \, {\norm{{\overline{\psi}}_e }}^2_{\sob {1,2} {} \home}  + {\tilde \e_3\over 2} \, {\norm{M_i}}^2_{\L {\infty} {} \home} \cdot {\norm{\pottr}}^2_{\sob {1,2} {} \home}\,, & (2.120)
\cr}$$
which holds true for arbitrary $\tilde \e_2$, $\tilde \e_3 > 0$. Specifying $\tilde \e_2 \tilde \e_3 = \e^\p_2$ with arbitrary $\e^\p_2 > 0$ and normalizing the constants in an appropriate way, (2.118) and (2.120) yield the desired estimate
$$\eqalignno{
{} & \bigl\vert \,\langle \, S(t)\,,\, \pottr\,\rangle\,\bigr\vert \,\le\, C\,\e^\p_2 \,{\norm{\pottr}}^2_{\sob {1,2} {} \home} + {C \over \e^\p_2} \,\Bigl(\, {\norm{ \excin(t)}}^2_{\L {2} {} \home }  + {\norm{ \excex(t)}}^2_{ \L {2} {} \home} \,\Bigr) \,, & (2.121)
\cr}$$
and the proof is complete. $\blackbox$
\bigskip\medskip
{\kapitel 3.~Correction of the proof of the stability estimate.}
\bigskip
{\bf a) Proof of Theorem 1.1.}
\medskip
{\it Part A. The Rogers-McCulloch model.} Let us specify within (2.9)$\,-\,$(2.11) the Rogers-McCulloch model. The proof, however, has been organized in such a way that the estimates work in the case of the linearized Aliev-Panfilov model as well. Throughout the following, $C$ denotes a generical positive constant, which may appropriately change from line to line. $C$ will never depend on the data $\Phi_0$, $W_0$, $\excin$ and $\excex$ but, possibly, on $\Omega$ and $p = 4$.
\medskip
$\bullet$ {\bf Step 1.} {\it The difference of the parabolic equations.} From the parabolic equations, satisfied by the pairs $(\pottr^\p, \gat^\p)$ and $(\pottr^\pp,\gat^\pp)$ for almost all $t \in [\,0\,,\,T\,]\,$, we obtain the difference
$$\eqalignno{
{} &  \langle\,{d \over dt} \,\bigl(\,\pottr^\p(t)- \pottr^\pp(t)\,\bigr)\,,\, \psi\,\rangle + {\cal M} \bigl(\, \pottr^\p(t)- \pottr^\pp(t)\,,\, \psi\,\bigr) & (3.1)
\cr
{} & \qquad + \,\int_\Omega \Bigl(\,\ioncurr (\pottr^\p(t), \gat^\p (t)) - \ioncurr (\pottr^\pp(t), \gat^\pp (t))\,\Bigr)\,\psi\, dx
\cr
{} & & \displaystyle \,=\, \langle \, {1 \over 1 + \la}\, \Bigl(\,\la\,\bigl(\, \excin^\p(t) - \excin^\pp(t)\,\bigr) - \bigl(\,\excex^\p (t)-
\excex^\pp(t)\,\bigr)\,\Bigr) \,,\,\psi\,\rangle \quad \forall\, \psi \in \sob {1,2}
{} \home\, .
\cr}$$
Inserting into $(3.1)$ the feasible test function $\psi = \pottr^\p(t)- \pottr^\pp(t)
\in \sob {1,2} {} \home$ and applying the lower estimate for the monodomain bilinear form, we arrive at
$$\eqalignno{
{} & \hbox{\kern-20pt} {1 \over 2} \,{d \over dt}\, {\norm{\pottr^\p - \pottr^\pp}}^2_{\L {2} {} \home} + \beta\, {\norm{\pottr^\p - \pottr^\pp }}^2_{\sob {1,2} {} \home}  +  \int_\Omega \Bigl(\,\ioncurr (\pottr^\p, \gat^\p )-  \ioncurr (\pottr^\pp, \gat^\pp )\,\Bigr)\,\bigl(\, \pottr^\p - \pottr^\pp\,\bigr)\, dx
\cr
{} & & \displaystyle \,\le\, \bigl\vert \, \langle\, {1 \over 1 + \la}\, \Bigl(\,\la\,\bigl(\, \excin^\p
- \excin^\pp \,\bigr) - \bigl(\,\excex^\p - \excex^\pp \,\bigr)\,\Bigr) \,,\,\pottr^\p - \pottr^\pp \,\rangle \,\bigr\vert + \beta\, {\norm{\pottr^\p - \pottr^\pp}}^2_{\L {2} {} \home} \,. \quad (3.2)
\cr}$$
The first term on the right-hand side will be estimated with the help of the generalized Cauchy inequality as follows:
$$\eqalignno{
{} & \hbox{\kern-20pt} \bigl\vert \, \langle\, {1 \over 1 + \la}\, \Bigl(\,\la\, \bigl(\, \excin^\p (t) - \excin^\pp (t) \,\bigr) - \bigl(\,\excex^\p (t)- \excex^\pp(t)\,\bigr)\,\Bigr)  \,,\,\pottr^\p - \pottr^\pp \,\rangle \,\bigr\vert \,
\cr
{} & \le\, {\norm{\, {1 \over 1 + \la}\,\Bigl(\,\la \,\bigl(\, \excin^\p (t) - \excin^\pp (t) \,\bigr) - \bigl(\,\excex^\p (t)- \excex^\pp(t)\,\bigr)\,\Bigr) \,}}{}_{\bigl(\,\sob {1,2} {} \home\,\bigr)^*} \, {\norm{\pottr^\p -\pottr^\pp}}_{\sob {1,2} {} \home} & (3.3)
\cr
{} & \le\, \e^\p_1 \,{\norm{\pottr^\p - \pottr^\pp}}^2_{\sob {1,2} {} \home} + {C \over \e^\p_1}\,
\Bigl(\, {\norm{\excin^\p- \excin^\pp}}^2_{\bigl(\,\sob {1,2} {} \home\,\bigr)^*} + {\norm{\excex^\p - \excex^\pp}}^2_{\bigl(\,\sob {1,2} {} \home\,\bigr)^*} \,\Bigr) & (3.4)
\cr}$$
with arbitrary $\e^\p_1 > 0$. For the term with the difference of the ionic currents on the left-hand side, we get a lower estimate with the help of the following lemma.
\medskip
{\bf Lemma 3.1.} {\it For all $\varphi_1$, $\varphi_2 \in \R {}$, the following identity holds:
$$\eqalignno{
{} &  \bigl(\,\varphi^3_1 - (a+1)\,\varphi^2_1 + a\, \varphi_1\,\bigr) -
\bigl(\,\varphi^3_2 - (a+1)\,\varphi^2_2 + a\, \varphi_2\,\bigr) \,
\cr
{} & & \displaystyle =\, (\varphi_1 - \varphi_2) \cdot\bigl(\, \varphi^2_1 +
\varphi_1\,\varphi_2 + \varphi^2_2 - (a+1)\,(\varphi_1 + \varphi_2) + a\,\bigr)\,.
\quad \blackbox \quad  (3.5)
\cr}$$}%
\par
Consequently, we find
$$\eqalignno{
{} & \hbox{\kern-20pt} \int_\Omega \Bigl(\,\ioncurr (\pottr^\p, \gat^\p )-  \ioncurr
(\pottr^\pp, \gat^\pp )\,\Bigr)\,\bigl(\, \pottr^\p - \pottr^\pp\,\bigr)\, dx  \, &
(3.6)
\cr
{} & =\, \int_\Omega (\pottr^\p - \pottr^\pp)\, b\,\bigl(\,(\pottr^\p)^2 + \pottr^\p
\,\pottr^\pp + (\pottr^\pp)^2  + a\,\bigr)\,(\pottr^\p - \pottr^\pp)\, dx
\cr
{} & & \displaystyle - \, (a+1)\,b\, \int_\Omega (\pottr^\p - \pottr^\pp)\,\bigl(\,
\pottr^\p + \pottr^\pp\,\bigr) \, (\pottr^\p - \pottr^\pp)\, dx + \int_\Omega
\bigl(\, \pottr^\p\,\gat ^\p - \pottr^\pp \, \gat^\pp\,\bigr)\,(\pottr^\p -
\pottr^\pp)\,dx\,.
\cr}$$
Since $\pottr^\p(x,t)^2 + \pottr^\p(x,t)\,\pottr^\pp(x,t) + \pottr^\pp(x,t)^2 \ge 0$
for almost all $(x,t) \in \Omega_T$ and $a$, $b > 0$, the inequalities (3.2), (3.4)
and (3.6) imply
$$\eqalignno{
{} &  \hbox{\kern-20pt} {d \over dt}\, {\norm{\pottr^\p - \pottr^\pp}}^2_{\L {2} {}
\home} + 2\, \beta\, {\norm{\pottr^\p - \pottr^\pp }}^2_{\sob {1,2} {} \home}
\, & (3.7)
\cr
{} & \le\,2\,C \int_\Omega \bigl\vert \,\pottr^\p - \pottr^\pp\,\bigr\vert\cdot
\bigl\vert\, \pottr^\p + \pottr^\pp\,\bigr\vert \cdot \bigl\vert\, \pottr^\p -
\pottr^\pp\,\bigr\vert\, dx
\cr
{} &  + 2\,\int_\Omega (\pottr^\p - \pottr^\pp)\,\gat^\p\,(\pottr^\p -
\pottr^\pp)\, dx \,+ \, 2\, \int_\Omega (\gat^\p - \gat^\pp)\, \pottr^\pp\,(\pottr^\p -
\pottr^\pp)\, dx
\cr
{} & + \,\e^\p_1 \, {\norm{\pottr^\p - \pottr^\pp}}^2 _{\sob
{1,2} {} \home} + 2\,\beta\, {\norm{\pottr^\p - \pottr^\pp}}^2_{\L {2} {} \home}
\cr
{} & + \, {C \over \e^\p_1} \, \Bigl(\, {\norm{\excin^\p - \excin^\pp }}^2_{\bigl(\,\sob {1,2} {} \home\,\bigr)^*} + {\norm{\excex^\p - \excex^\pp }}^2_{\bigl(\, \sob {1,2} {} \home\,\bigr)^*} \,\Bigr)\,.
\cr}$$
Now it must be emphasized that $(3.7)$ holds parametrically in $t$ for almost all fixed $t\in (\,0\,,\,T\,)$. In the subsequent applications of the generalized Cauchy inequality this will become important since the parameters $\e_i$ introduced there must be chosen in a time-dependent way.
\par
We apply first the generalized Cauchy inequality with $\e_2(t) > 0$ and subsequently H\"older's inequality to the first term on the right-hand side of (3.7), thus getting
$$\eqalignno{
{} & \hbox{\kern-20pt} 2\, C \int_\Omega \bigl\vert \,\pottr^\p (t) - \pottr^\pp(t)\,\bigr\vert\cdot
\bigl\vert\, \pottr^\p(t) + \pottr^\pp(t)\,\bigr\vert \cdot \bigl\vert\, \pottr^\p (t)-
\pottr^\pp(t)\,\bigr\vert\, dx \,
\cr
{} & \le\, C\,\e_2 (t) \, \Bigl(\, \int_\Omega \bigl\vert\,\pottr^\p +
\pottr^\pp \,\bigr\vert^4\, dx\,\Bigr)^{1/2} \, \Bigl(\, {\norm{\pottr^\p -
\pottr^\pp}}^4_{\L {4} {} \home} \,\Bigr)^{1/2} + {C \over \e_2 (t)}\,
{\norm{\pottr^\p - \pottr^\pp}}^2_{\L {2} {} \home} & (3.8)
\cr
{} & \le\, C\,\e_2(t) \, \Bigl(\,{\norm{\pottr^\p}}^2_{\L {4} {} \home} + {\norm{\pottr^\pp}}^2_{\L {4} {} \home}\,\Bigr) \, {\norm{\pottr^\p - \pottr^\pp}}^2_{\sob {1,2} {} \home} + {C \over \e_2 (t)}\, {\norm{\pottr^\p - \pottr^\pp}}^2_{\L {2} {} \home} \,. & (3.9)
\cr}$$
Inserting now $\e_2(t) = \e^\p_2 / \bigl(\,1 + {\norm{\pottr^\p(t)}}^2_{\L {4} {} \home} + {\norm{\pottr^\pp(t)}}^2_{\L {4} {} \home}\,\bigr)$ with arbitrary $\e^\p_2 > 0$, we get
$$\eqalignno{
{} & \lldots \le \, C\,\e^\p_2 \, {\norm{\pottr^\p - \pottr^\pp}}^2_{\sob {1,2} {} \home} + {C \over \e^\p_2 }\, \Bigl(\,1 + {\norm{\pottr^\p}}^2_{\L {4} {} \home} + {\norm{\pottr^\pp}}^2_{\L {4} {} \home}\,\Bigr)\,{\norm{\pottr^\p - \pottr^\pp}}^2_{\L {2} {} \home}\,. & (3.10)
\cr}$$
In order to estimate the second term on the right-hand side of $(3.7)$, let us write
$$\eqalignno{
{} & \hbox{\kern-20pt} 2\,\int_\Omega \bigl\vert \,\gat^\p\,\bigr\vert \cdot \bigl\vert \,\pottr^\p - \pottr^\pp\,\bigr\vert^2 \,dx \,= \,2\,\int_\Omega \Bigl(\,\bigl\vert \,\gat^\p\,\bigr\vert^{2/3} \cdot \bigl\vert \,\pottr^\p - \pottr^\pp\,\bigr\vert\,\Bigr)\,\Bigl(\,\bigl\vert \,\gat^\p\,\bigr\vert^{1/3} \cdot \bigl\vert \,\pottr^\p - \pottr^\pp\,\bigr\vert\,\Bigr)\, \,dx & (3.11)
\cr
{} & \le \,C\,\e_3(t) \int_\Omega \bigl\vert \,\gat^\p\,\bigr\vert^{4/3} \cdot \bigl\vert \,\pottr^\p - \pottr^\pp\,\bigr\vert^2\,dx + {C \over \e_3 (t)} \int_\Omega \bigl\vert \,\gat^\p\,\bigr\vert^{2/3} \cdot \bigl\vert \,\pottr^\p - \pottr^\pp\,\bigr\vert^2\,dx & (3.12)
\cr
{} & \le C\, \e_3(t) \,\Bigl(\, \int_\Omega \bigl\vert \,\gat^\p\,\bigr\vert^{8/3}\, dx\Bigr)^{1/2} \,\Bigl(\, \int_\Omega \bigl\vert \,\pottr^\p - \pottr^\pp\,\bigr\vert^4\,dx\,\Bigr)^{1/2} & (3.13)
\cr
{} & & \displaystyle +\, {C \over \e_3(t)} \,\Bigl(\, \tilde\e_3(t) \int_\Omega \bigl\vert \,\gat^\p\,\bigr\vert^{4/3} \cdot \bigl\vert \,\pottr^\p - \pottr^\pp\,\bigr\vert^2\,dx + {C \over \tilde\e_3(t)} \int_\Omega \bigl\vert\, \pottr^\p - \pottr^\pp\,\bigr\vert^2\,dx\,\Bigr)
\cr
{} & \le C\, \e_3(t) \,\Bigl(\, \int_\Omega \bigl\vert \,\gat^\p\,\bigr\vert^{8/3}\, dx\Bigr)^{1/2} \,\Bigl(\, \int_\Omega \bigl\vert \,\pottr^\p - \pottr^\pp\,\bigr\vert^4\,dx\,\Bigr)^{1/2}  & (3.14)
\cr
{} & & \displaystyle +\, {C \,\tilde\e_3(t) \over \e_3(t)} \,\Bigl(\, \int_\Omega \bigl\vert \,\gat^\p\,\bigr\vert^{8/3}\, dx\Bigr)^{1/2} \,\Bigl(\, \int_\Omega \bigl\vert \,\pottr^\p - \pottr^\pp\,\bigr\vert^4\,dx\,\Bigr)^{1/2} + {C \over \e_3(t)\, \tilde\e_3(t)} \int_\Omega \bigl\vert\, \pottr^\p - \pottr^\pp\,\bigr\vert^2\,dx
\cr
{} & \le C\, \e_3(t) \,{\norm{\gat^\p}}^{4/3}_{\L {8/3} {} \home} \cdot {\norm{\pottr^\p - \pottr^\pp}}^2_{\sob {1,2} {} \home} & (3.15)
\cr
{} & & \displaystyle +\, {C \,\tilde\e_3(t) \over \e_3(t)}\, {\norm{\gat^\p}}^{4/3}_{\L {8/3} {} \home} \cdot {\norm{\pottr^\p - \pottr^\pp}}^2_{\sob {1,2} {} \home} + {C \over \e_3(t)\, \tilde\e_3(t)}\,{\norm{\pottr^\p - \pottr^\pp}}^2_{\L {2} {} \home}\,.
\cr}$$
By Proposition 2.1., 2), $\gat^\p$ belongs to $\bochnerc {0} {[\,0\,,\,T\,]} {\L {4} {} \home } \hookrightarrow \bochnerl {\infty} {(\,0\,,\,T\,)} {\L {8/3} {} \home} $. Then with $\e_3(t) = \e^\p_3$, $\tilde\e_3(t) = \e^\pp_3$, we get
$$\eqalignno{
{} & \int_\Omega \bigl\vert \,\gat^\p\,\bigr\vert \cdot \bigl\vert \,\pottr^\p - \pottr^\pp\,\bigr\vert^2 \,dx \,\le \, C \Bigl(\, \e^\p_3 + {\e^\pp_3 \over \e^\p_3} \,\Bigr) \,{\norm{\gat^\p}}^{4/3}_{\bochnerc {0} {[\,0\,,\,T\,]} {\L {4} {} \home} } \cdot {\norm{\pottr^\p - \pottr^\pp}}^2_{\sob {1,2} {} \home} & (3.16)
\cr
{} & & \displaystyle +\,{C \over \e^\p_3\, \e^\pp_3}\,{\norm{\pottr^\p - \pottr^\pp}}^2_{\L {2} {} \home}\,.
\cr}$$
For the third term from the right-hand side of (3.7), we find
$$\eqalignno{
{} & \hbox{\kern-20pt} 2\, \int_\Omega (\gat^\p - \gat^\pp)\, \pottr^\pp\,(\pottr^\p -
\pottr^\pp)\, dx
\cr
{} & \le \, C\,\e_4 (t)\,\int_\Omega \bigl(\, \pottr^\pp\,\bigr)^2
\,\bigl\vert \,\pottr^\p - \pottr^\pp\,\bigr\vert^2 \, dx + {C \over \e_4 (t)}
\,{\norm{\gat^\p - \gat^\pp}}^2_{\L {2} {} \home} & (3.17)
\cr
{} & \le \, C\,\e_4 (t) \,\Bigl(\, \int_\Omega \bigl(\,\pottr^\pp\,\bigr)^4\,
dx\,\Bigr)^{1/2} \,\Bigl(\,{\norm{\pottr^\p - \pottr^\pp}}^4_{\L {4} {}
\home}\,\Bigr)^{1/2} + {C \over \e_4 (t)} \,{\norm{\gat^\p - \gat^\pp}}^2_{\L {2} {}
\home}& (3.18)
\cr
{} & \le \, C\, \e_4(t) \,{\norm{\pottr^\pp}}^2_{\L {4} {}
\home}\,{\norm{\pottr^\p - \pottr^\pp}}^2_{\sob {1,2} {} \home} + {C \over \e_4 (t)}
\,{\norm{\gat^\p - \gat^\pp}}^2_{\L {2} {} \home}\,, & (3.19)
\cr}$$
using the (noncompact but) continuous imbedding $\sob {1,2} {} \home \hookrightarrow \L {6} {} \home$ and applying again the generalized Cauchy inequality with $\e_4(t) > 0$. Specifying now $\e_4(t) = \e^\p_4 /\bigl(\, 1+ {\norm {\pottr^\pp (t)}}^2_{\L {4} {} \home} \,\bigr)$ with arbitrary $\e^\p_4 > 0$, we may continue
$$\eqalignno{
{} & \lldots \le\, C\, \e^\p_4 \,{\norm{\pottr^\p - \pottr^\pp}}^2_{\sob {1,2} {} \home} + {C \over \e^\p_4 } \,\Bigl(\, 1+ {\norm {\pottr^\pp (t)}}^2_{\L {4} {} \home}\,\Bigr)\,{\norm{\gat^\p - \gat^\pp}}^2_{\L {2} {} \home}\,. & (3.20)
\cr}$$
Assembling now (3.7), (3.10), (3.16) and (3.20), we arrive at the following inequality:
$$\eqalignno{
{} &  \hbox{\kern-20pt} {d \over dt}\, {\norm{\pottr^\p - \pottr^\pp}}^2_{\L {2} {}
\home} + 2\, \beta\, {\norm{\pottr^\p - \pottr^\pp }}^2_{\sob {1,2} {} \home} &
(3.21)
\cr
{} & \le \, {C \over \e^\p_1} \,\Bigl(\, {\norm{\excin^\p(t) - \excin^\pp(t) }}^2_{\bigl(\,\sob {1,2} {} \home\,\bigr)^*} + {\norm{\excex^\p(t) - \excex^\pp(t) }}^2_{\bigl(\,\sob {1,2} {} \home\,\bigr)^*} \,\Bigr)
\cr
{} & + \,\e^\p_1 \, {\norm{\pottr^\p - \pottr^\pp}}^2 _{\sob
{1,2} {} \home} + 2\,\beta\, {\norm{\pottr^\p - \pottr^\pp}}^2_{\L {2} {} \home}
\cr
{} & +\, C\,\e^\p_2 \, {\norm{\pottr^\p - \pottr^\pp}}^2_{\sob {1,2} {} \home} + {C \over \e^\p_2}\, \Bigl(\, 1 + {\norm{\pottr^\p}}^2_{\L {4} {} \home} + {\norm{\pottr^\pp}}^2_{\L {4} {} \home} \,\Bigr)\, {\norm{\pottr^\p - \pottr^\pp}}^2_{\L {2} {} \home}
\cr
{} & + \, C \Bigl(\, \e^\p_3 + {\e^\pp_3 \over \e^\p_3} \,\Bigr) \cdot {\norm{\pottr^\p - \pottr^\pp}}^2_{\sob {1,2} {} \home} + {C \over \e^\p_3\, \e^\pp_3}\,{\norm{\pottr^\p - \pottr^\pp}}^2_{\L {2} {} \home}\,
\cr
{} & + \,C\, \e^\p_4 \, {\norm{\pottr^\p - \pottr^\pp}}^2_{\sob {1,2} {} \home} + {C \over \e^\p_4} \,\Bigl(\, 1 + {\norm{\pottr^\pp}}^2_{\L {4} {} \home}\,\Bigr)\,{\norm{\gat^\p - \gat^\pp}}^2_{\L {2} {} \home}\,.
\cr}$$
Now we may fix the numbers $\e^\p_1$, $\e^\p_2$, $\e^\p_3$, $\e^\pp_3$, $\e^\p_4 > 0$ in such a way that the terms with ${\norm{\pottr^\p - \pottr^\pp}}^2_{\sob {1,2} {} \home}$ on both sides of $(3.21)$ will be annihilated. We arrive at
$$\eqalignno{
{} &  \hbox{\kern-20pt} {d \over dt}\, {\norm{\pottr^\p - \pottr^\pp}}^2_{\L {2} {} \home} \le C\,\Bigl(\, 2\,\beta +   {1 \over \e^\p_2}\, \bigl(\, 1 + {\norm{\pottr^\p}}^2_{\L {4} {} \home} + {\norm{\pottr^\pp}}^2_{\L {4} {} \home}\,\bigr)+ \, {1 \over \e^\p_3\,\e^\pp_3}\,\Bigr)\, {\norm{\pottr^\p - \pottr^\pp}}^2_{\L {2} {} \home}
& (3.22)
\cr
{} &   + \, {C \over \e^\p_4} \,\Bigl(\, 1 + {\norm{\pottr^\pp}}^2_{\L {4} {} \home}\,\Bigr)\,{\norm{\gat^\p - \gat^\pp}}^2_{\L {2} {} \home}
\cr
{} &  + \, {C \over \e^\p_1} \,\Bigl(\, {\norm{\excin^\p(t) - \excin^\pp(t) }}^2_{\bigl(\,\sob {1,2} {} \home\,\bigr)^*} + {\norm{\excex^\p(t) - \excex^\pp(t) }}^2_{\bigl(\,\sob {1,2} {} \home\,\bigr)^*} \,\Bigr)\,.
\cr}$$
\par
$\bullet$ {\bf Step 2.} {\it The difference of the gating equations.} Inserting into the difference of the gating equations for $(\pottr^\p, \gat^\p)$ and $(\pottr^\pp, \gat^\pp)$,
$$\eqalignno{
{} & \langle \, {d \over dt} \,\bigl(\, \gat^\p (t) - \gat^\pp(t)\,\bigr) \,,\,
\psi\,\rangle \,=\, -\, \e\,\int_\Omega \bigl(\,\gat^\p (t)-
\gat^\pp(t)\,\bigr)\,\psi\, dx + \,\e\,\kappa\, \int_\Omega \bigl(\,\pottr^\p(t)-
\pottr^\pp(t)\,\bigr)\,\psi\, dx & (3.23)
\cr
{} & & \displaystyle  \quad \forall\, \psi \in \L {2} {} \home\,,
\cr}$$
the feasible test function $\psi = \gat^\p(t) - \gat^\pp(t)$ and applying Cauchy's
inequality to the second term, we get the estimate$\,$\footnote{$^{24)}$}{Note that here $\e > 0$ is the given one from $(2.5)$.}
$$\eqalignno{
{} & {d \over dt}\,{\norm{\gat^\p - \gat^\pp}}^2_{\L {2} {} \home} \,\le\,
2\,\bigl(\,2\,\e + \e\,\kappa\,\bigr)\,{\norm{\gat^\p - \gat^\pp}}^2_{\L {2} {} \home}
+ 2\,\e\,\kappa\, {\norm{\pottr^\p - \pottr^\pp}}^2_{\L {2} {} \home}\,. & (3.24)
\cr}$$
\par
$\bullet$ {\bf Step 3.} {\it The estimates for the differences ${\norm{\pottr^\p - \pottr^\pp}}^2_{\bochnerl {\infty} {(\,0\,,\,T\,)} {\L {2} {} \home} }$, ${\norm{\gat^\p - \gat^\pp}}^2_{\bochnerl {2} {(\,0\,,\,T\,)} {\L {2} {} \home} }$ and ${\norm{\gat^\p - \gat^\pp}}^2_{\bochnerl {\infty} {(\,0\,,\,T\,)} {\L {2} {} \home} }$.} After enlarging and equalizing the factors on the right-hand sides, the
in\-equalities (3.22) and (3.24) yield together
$$\eqalignno{
{} & \hbox{\kern-20pt} {d \over dt}\,\Bigl(\, {\norm{\pottr^\p - \pottr^\pp}}^2_{\L
{2} {} \home} + {\norm{\gat^\p - \gat^\pp}}^2_{\L {2} {} \home} \,\Bigr)\,  & (3.25)
\cr
{} & \le\, C\,\Bigl(\, 2\,\beta +   {1 \over \e^\p_2}\, \bigl(\, 1 + {\norm{\pottr^\p(t)}}^2_{\L {4} {} \home} + {\norm{\pottr^\pp(t)}}^2_{\L {4} {} \home}\,\bigr)+ \, {1 \over \e^\p_3\,\e^\pp_3} \,+ 2\,\e\,\kappa \,\Bigr)\, {\norm{\pottr^\p - \pottr^\pp}}^2_{\L {2} {} \home}
\cr
{} & + \, C\,\Bigl(\, {1 \over \e^\p_4} \,\bigl(\, 1 + {\norm{\pottr^\pp(t)}}^2_{\L {4} {} \home}\,\bigr) + 4\,\e + 2\,\e\,\kappa \,\Bigr)\,{\norm{\gat^\p - \gat^\pp}}^2_{\L {2} {} \home}
\cr
{} & + \, {C \over \e^\p_1} \,\Bigl(\, {\norm{\excin^\p(t) - \excin^\pp(t) }}^2_{\bigl(\,\sob {1,2} {} \home\,\bigr)^*} + {\norm{\excex^\p(t) - \excex^\pp(t) }}^2_{\bigl(\,\sob {1,2} {} \home\,\bigr)^*} \,\Bigr)\quad \follows
\cr
{} & \hbox{\kern-20pt} {d \over dt}\,\Bigl(\, {\norm{\pottr^\p - \pottr^\pp}}^2_{\L
{2} {} \home} + {\norm{\gat^\p - \gat^\pp}}^2_{\L {2} {} \home} \,\Bigr)\,\le\, A(t) \cdot\Bigl(\, {\norm{\pottr^\p - \pottr^\pp}}^2_{\L {2} {} \home} + {\norm {\gat^\p - \gat^\pp}}^2_{\L {2} {} \home} \,\Bigr)  & (3.26)
\cr
{} & & \displaystyle  +  \,{C \over \e^\p_1} \, \Bigl(\,  {\norm{\excin^\p(t) - \excin^\pp(t) }}^2_{\bigl(\,\sob {1,2} {} \home\,\bigr)^*} + {\norm{\excex^\p(t) - \excex^\pp(t) }}^2_{\bigl(\,\sob {1,2} {} \home\,\bigr)^*} \,\Bigr)\,
\cr}$$
where
$$\eqalignno{
{} & \hbox{\kern-20pt}A(t) = C\,\Bigl(\, 1+ {1 \over \e^\p_2}\, \bigl(\, 1 + {\norm{\pottr^\p(t)}}^2_{\L {4} {} \home} + {\norm{\pottr^\pp(t)}}^2_{\L {4} {} \home}\,\bigr)  + \, {1 \over \e^\p_3\,\e^\pp_3} +  {1 \over \e^\p_4} \,\bigl(\, 1 + {\norm{\pottr^\pp(t)}}^2_{\L {4} {} \home}\,\bigr)\,\Bigr)\,.& (3.27)
\cr}$$
Gronwall's inequality implies now that
$$\eqalignno{
{} & \hbox{\kern-20pt} {\norm{\pottr^\p(t) - \pottr^\pp(t)}}^2_{\L {2} {} \home} +
{\norm{\gat^\p(t) - \gat^\pp(t)}}^2_{\L {2} {} \home} \,
\le\, e^{\int^t_0 A(s)\, ds}\,\Bigl(\,
{\norm{\pottr^\p(0) - \pottr^\pp(0)}}^2_{\L {2} {} \home} & (3.28)
\cr
{} &  + \, {\norm{\gat^\p(0) - \gat^\pp(0)}}^2_{\L {2} {} \home} + {C \over \e^\p_1} \,
\int^t_0  \bigl(\, {\norm{\excin^\p(\tau) - \excin^\pp(\tau) }}^2_{\bigl(\,\sob
{1,2} {} \home\,\bigr)^*} + {\norm{\excex^\p(\tau) - \excex^\pp(\tau)
}}^2_{\bigl(\,\sob {1,2} {} \home\,\bigr)^*} \,\bigr) \, d\tau\,\Bigr)
\cr
{} & \le \, e^{\tilde A} \, {C \over \e^\p_1} \, \Bigl(\, {\norm{\excin^\p - \excin^\pp
}}^2_{\bochnerl {2} {(\,0\,,\,T\,)} {\bigl(\,\sob {1,2} {} \home\,\bigr)^*} } +
{\norm{\excex^\p - \excex^\pp }}^2_{\bochnerl {2} {(\,0\,,\,T\,)} {\bigl(\,\sob
{1,2} {} \home\,\bigr)^*} } \,\Bigr)  & (3.29)
\cr}$$
with
$$\eqalignno{
{} & \hbox{\kern-20pt} \widetilde A \,=\, \int^T_0 A(s) \, ds \,=\, C\,T + {C \over e^\p_2}\, \Bigl(\, T + \int^T_0 \bigl(\,{\norm{\pottr^\p(s)}}^2_{\L {4} {} \home} + {\norm{\pottr^\pp(s)}}^2_{\L {4} {} \home}\,\bigr) \,ds\,\Bigr) & (3.30)
\cr
{} & & \displaystyle + \, {C \,T \over \e^\p_3\,\e^\pp_3} +  {C \over \e^\p_4} \,\Bigl(\, T + \int^T_0 {\norm{\pottr^\pp(s)}}^2_{\L {4} {} \home}\,ds\,\Bigr)
\cr
{} & \hbox{\kern-20pt} \le \, C\,\Bigl(\, 1 + {\norm{\pottr^\p}}^4_{\L {4} {} (\Omega_T)} + {\norm{\pottr^\pp}}^4_{\L {4} {} (\Omega_T)} \,\Bigr)\, & (3.31)
\cr
{} & \hbox{\kern-20pt} \le \, C\, \Bigl(\,1 + {\norm{\Phi_0}}^2_{\L {2} {} \home} + {\norm{W_0}}^2_{\L {2} {} \home} +  {\norm{\excin^\p}}^2_{\bochnerl {\infty} {(\,0\,,\,T\,)} {\bigl(\,\sob {1,2} {} \home\,\bigr)^*} }  +  {\norm{\excex^\p}}^2 _{\bochnerl {\infty} {(\,0\,,\,T\,)} {\bigl(\,\sob {1,2} {} \home\,\bigr)^*} }  & (3.32)
\cr
{} & & \displaystyle +  \,{\norm{\excin^\pp}}^2_{\bochnerl {\infty} {(\,0\,,\,T\,)} {\bigl(\,\sob {1,2} {} \home\,\bigr)^*} }  +  {\norm{\excex^\pp}}^2 _{\bochnerl {\infty} {(\,0\,,\,T\,)} {\bigl(\,\sob {1,2} {} \home\,\bigr)^*} }  \,\Bigr)
\cr}$$
$$\eqalignno{
{} & \le \, C\, \Bigl(\, 1 + {\norm{\Phi_0}}^2_{\L {2} {} \home} + {\norm{W_0}}^2_{\L {2} {} \home} + 4\, R^2\,\Bigr) & (3.33)
\cr}$$
by (2.13) and the assumption about the uniform bound $R > 0$ for the norms of the inhomogeneities. Summing up, we obtain from inequality $(3.29)$ the following estimates:
$$\eqalignno{
{} & \hbox{\kern-20pt} {\norm{\pottr^\p - \pottr^\pp}}^2_{\bochnerl {\infty}
{(\,0\,,\,T\,)} {\L {2} {} \home} } \,\le\, e^{\tilde A} \, {C \over \e^\p_1}\, \Bigl(\,
{\norm{\excin^\p - \excin^\pp }}^2_{\bochnerl {2} {(\,0\,,\,T\,)} {\bigl(\,\sob
{1,2} {} \home\,\bigr)^*} } \cr
{} & & \displaystyle +\, {\norm{\excex^\p - \excex^\pp }}^2_{\bochnerl {2}
{(\,0\,,\,T\,)} {\bigl(\,\sob {1,2} {} \home\,\bigr)^*} } \,\Bigr) \,; \quad (3.34)
\cr
{} & \hbox{\kern-20pt} {\norm{\gat^\p - \gat^\pp}}^2_{\bochnerl {\infty}
{(\,0\,,\,T\,)} {\L {2} {} \home} } \,\le\, e^{\tilde A} \, {C \over \e^\p_1} \, \Bigl(\,
{\norm{\excin^\p - \excin^\pp }}^2_{\bochnerl {2} {(\,0\,,\,T\,)} {\bigl(\,\sob
{1,2} {} \home\,\bigr)^*} }
\cr
{} & & \displaystyle + \,{\norm{\excex^\p - \excex^\pp }}^2_{\bochnerl {2}
{(\,0\,,\,T\,)} {\bigl(\,\sob {1,2} {} \home\,\bigr)^*} } \,\Bigr) \,; \quad (3.35)
\cr
{} & \hbox{\kern-20pt} {\norm{\gat^\p - \gat^\pp}}^2_{\bochnerl {2} {(\,0\,,\,T\,)}
{\L {2} {} \home} } \,\le\, T\,e^{\tilde A} \, {C \over \e^\p_1} \, \Bigl(\,
{\norm{\excin^\p - \excin^\pp }}^2_{\bochnerl {2} {(\,0\,,\,T\,)} {\bigl(\,\sob
{1,2} {}  \home\,\bigr)^*} }
\cr
{} & & \displaystyle + \,{\norm{\excex^\p - \excex^\pp }}^2_{\bochnerl {2}
{(\,0\,,\,T\,)} {\bigl(\,\sob {1,2} {} \home\,\bigr)^*} } \,\Bigr) \,. \quad (3.36)
\cr}$$
\par
$\bullet$ {\bf Step 4.} {\it The estimate for the difference ${\norm{\pottr^\p -
\pottr^\pp }}^2_{\bochnerl {2} {(\,0\,,\,T\,)} {\sob {1,2} {} \home} }$.} In (3.21),
the numbers $\e^\p_1$, ... ,\break $\e^\p_4 > 0$ may be alternatively chosen in such a way that
$$\eqalignno{
{} & \hbox{\kern-20pt} {d \over dt}\, {\norm{\pottr^\p - \pottr^\pp}}^2_{\L {2} {}
\home} + \beta\, {\norm{\pottr^\p - \pottr^\pp }}^2_{\sob {1,2} {} \home} \,\le \,C\,
\Bigl(\, 1 +   {1 \over \e^\p_2}\, \bigl(\, 1 + {\norm{\pottr^\p}}^2_{\L {4} {} \home} + {\norm{\pottr^\pp}}^2_{\L {4} {} \home}\,\bigr)
& (3.37)
\cr
{} & + \, {1 \over \e^\p_3\,\e^\pp_3} \,\Bigr)\, {\norm{\pottr^\p - \pottr^\pp}}^2_{\L {2} {} \home} + \, {C \over \e^\p_4} \,\Bigl(\, 1 + {\norm{\pottr^\pp}}^2_{\L {4} {} \home}\,\Bigr)\,{\norm{\gat^\p - \gat^\pp}}^2_{\L {2} {} \home}
\cr
{} & + \, {C \over \e^\p_1} \,\Bigl(\, {\norm{\excin^\p(t) - \excin^\pp(t) }}^2_{\bigl(\,\sob {1,2} {} \home\,\bigr)^*} + {\norm{\excex^\p(t) - \excex^\pp(t) }}^2_{\bigl(\,\sob {1,2} {} \home\,\bigr)^*} \,\Bigr)\,.
\cr}$$
This implies the following modification of (3.26):
$$\eqalignno{
{} & \hbox{\kern-20pt} {d \over dt}\,\Bigl(\, {\norm{\pottr^\p - \pottr^\pp}}^2_{\L
{2} {} \home} + {\norm{\gat^\p - \gat^\pp}}^2_{\L {2} {} \home} \,\Bigr) + \beta\, {\norm{\pottr^\p - \pottr^\pp }}^2_{\sob {1,2} {} \home}\, & (3.38)
\cr
{} & \le\, A(t) \cdot\Bigl(\, {\norm{\pottr^\p - \pottr^\pp}}^2_{\L {2} {} \home} + {\norm {\gat^\p - \gat^\pp}}^2_{\L {2} {} \home} \,\Bigr)
\cr
{} & & \displaystyle  +  \,{C \over \e^\p_1} \, \Bigl(\,  {\norm{\excin^\p(t) - \excin^\pp(t) }}^2_{\bigl(\,\sob {1,2} {} \home\,\bigr)^*} + {\norm{\excex^\p(t) - \excex^\pp(t) }}^2_{\bigl(\,\sob {1,2} {} \home\,\bigr)^*} \,\Bigr)\,
\cr}$$
Together with $(3.34)$ and $(3.35)$, we obtain
$$\eqalignno{
{} & \hbox{\kern-20pt} {d \over dt}\,\Bigl(\, {\norm{\pottr^\p - \pottr^\pp}}^2_{\L
{2} {} \home} + {\norm{\gat^\p - \gat^\pp}}^2_{\L {2} {} \home} \,\Bigr) + \beta\, {\norm{\pottr^\p - \pottr^\pp }}^2_{\sob {1,2} {} \home}\, & (3.39)
\cr
{} & \le\,2\, A(t) \,  e^{\tilde A} \, {C \over \e^\p_1}\, \Bigl(\,
{\norm{\excin^\p - \excin^\pp }}^2_{\bochnerl {2} {(\,0\,,\,T\,)} {\bigl(\,\sob
{1,2} {} \home\,\bigr)^*} }  + {\norm{\excex^\p - \excex^\pp }}^2_{\bochnerl {2}
{(\,0\,,\,T\,)} {\bigl(\,\sob {1,2} {} \home\,\bigr)^*} } \,\Bigr)
\cr
{} & & \displaystyle  +  \,{C \over \e^\p_1} \, \Bigl(\,  {\norm{\excin^\p(t) - \excin^\pp(t) }}^2_{\bigl(\,\sob {1,2} {} \home\,\bigr)^*} + {\norm{\excex^\p(t) - \excex^\pp(t) }}^2_{\bigl(\,\sob {1,2} {} \home\,\bigr)^*} \,\Bigr)\,.
\cr}$$
We integrate $(3.39)$ over $[\,0\,,\,T\,]$ and find, inserting the identical initial values
$$\eqalignno{
{} & {\norm{\pottr^\p(T) - \pottr^\pp(T)}}^2_{\L
{2} {} \home} + {\norm{\gat^\p(T) - \gat^\pp(T)}}^2_{\L {2} {} \home}  + \beta\, {\norm{\pottr^\p - \pottr^\pp }}^2_{\bochnerl {2} {(\,0\,,\,T\,)} {\sob {1,2} {} \home} } & (3.40)
\cr
{} & \le \, C \, \int^T_0 A(t)\, dt \cdot \Bigl(\, {\norm{\excin^\p - \excin^\pp }}^2_{\bochnerl {2} {(\,0\,,\,T\,)} {\bigl(\,\sob {1,2} {} \home\,\bigr)^*} }  + {\norm{\excex^\p - \excex^\pp }}^2_{\bochnerl {2}
{(\,0\,,\,T\,)} {\bigl(\,\sob {1,2} {} \home\,\bigr)^*} } \,\Bigr)
\cr
{} & & \displaystyle  +  \,{C \over \e^\p_1} \, \Bigl(\, {\norm{\excin^\p - \excin^\pp }}^2_{\bochnerl {2} {(\,0\,,\,T\,)} {\bigl(\,\sob {1,2} {} \home\,\bigr)^*} }+ {\norm{\excex^\p - \excex^\pp }}^2_{\bochnerl {2} {(\,0\,,\,T\,)} {\bigl(\,\sob {1,2} {} \home\,\bigr)^*} }\,\Bigr)
\cr
{} & \hbox{\kern-20pt} \le C\,\bigl(\,\widetilde A + {1 \over \e^\p_1}\,\bigr)\,\Bigl(\, {\norm{\excin^\p - \excin^\pp }}^2_{\bochnerl {2} {(\,0\,,\,T\,)} {\bigl(\,\sob {1,2} {} \home\,\bigr)^*} }+ {\norm{\excex^\p - \excex^\pp }}^2_{\bochnerl {2} {(\,0\,,\,T\,)} {\bigl(\,\sob {1,2} {} \home\,\bigr)^*} }\,\Bigr)\,. & (3.41)
\cr}$$
This implies the desired estimate
$$\eqalignno{
{} & {\norm{\pottr^\p - \pottr^\pp }}^2_{\bochnerl {2} {(\,0\,,\,T\,)} {\sob {1,2} {} \home} } \, & (3.42)
\cr
{} & \le \, {C\over \beta}\,\bigl(\,\widetilde A + {1 \over \e^\p_1}\,\bigr)\,\Bigl(\, {\norm{\excin^\p(t) - \excin^\pp(t) }}^2_{\bochnerl {2} {(\,0\,,\,T\,)} {\bigl(\,\sob {1,2} {} \home\,\bigr)^*} }+ {\norm{\excex^\p(t) - \excex^\pp(t) }}^2_{\bochnerl {2} {(\,0\,,\,T\,)} {\bigl(\,\sob {1,2} {} \home\,\bigr)^*} }\,\Bigr)\,.
\cr}$$
\par
$\bullet$ {\bf Step 5.} {\it The estimate for the difference ${\norm{\gat^\p - \gat^\pp }}^2_{\bochnersob {1,2} {(\,0\,,\,T\,)} {\L {2} {} \home} }$.} Into equation (2.21), we insert the test function $\psi = (\partial\gat^\p(t)/\partial t) - (\partial\gat^\pp(t)/\partial t)$ which, by Proposition 2.1., 1), belongs to $\L {2} {} (\Omega_T)$ and is therefore admissible. Then we get with the generalized Cauchy inequality $\vphantom{ \displaystyle\int  }$
$$\eqalignno{
{} & \hbox{\kern-20pt} \langle \, {\partial \gat^\p \over \partial t} - {\partial
\gat^\pp \over \partial t} \,,\, {\partial \gat^\p \over \partial t} - {\partial
\gat^\pp \over \partial t} \,\rangle \,=\, {\norm{{\partial \gat^\p \over \partial
t} - {\partial \gat^\pp \over \partial t} }}{}^2_{\L {2} {} \home} & (3.43)
\cr
{} & \le \, 2\e\,\e^\p_5 \, {\norm{{\partial \gat^\p \over \partial t} -
{\partial \gat^\pp \over \partial t} }}{}^2_{\L {2} {} \home}  + {\e \over 2\,
\e^\p_5}\,{\norm{\gat^\p -\gat^\pp}}^2_{\L {2} {} \home}
\cr
{} & & \displaystyle + \,2\,\e\,\kappa\, \e^\p_6 \, {\norm{{\partial \gat^\p
\over \partial t} - {\partial \gat^\pp \over \partial t} }}{}^2_{\L {2} {} \home}  +
{\e\,\kappa \over 2\,\e^\p_6}\,{\norm{\pottr^\p - \pottr^\pp}}^2_{\L {2} {} \home}
\cr}$$
for arbitrary $\e^\p_5$, $\e^\p_6 > 0$. Fixing the numbers $\e_5$ and $\e_6$ in such a way
that $2\,\e \,\e^\p_5  + 2\,\e\,\kappa\, \e^\p_6  = 1/2 $, we find together with (3.34) and (3.35):
$$\eqalignno{
{} & \hbox{\kern-20pt} {\norm{{\partial \gat^\p \over \partial t} - {\partial
\gat^\pp \over \partial t} }}{} ^2_{\L {2} {} \home}
\,\le\, {\e \over 2\,\e^\p_5}\,{\norm{\gat^\p -\gat^\pp}}^2_{\L {2} {} \home} +
{\e\,\kappa \over 2\,\e^\p_6}\,{\norm{\pottr^\p - \pottr^\pp}}^2_{\L {2} {} \home} \, &
(3.44)
\cr
{} & & \displaystyle \le\, C\, \Bigl(\, {\norm{\excin^\p - \excin^\pp
}}^2_{\bochnerl {2} {(\,0\,,\,T\,)} {\bigl(\,\sob {1,2} {}  \home\,\bigr)^*} } +
{\norm{\excex^\p - \excex^\pp }}^2_{\bochnerl {2} {(\,0\,,\,T\,)} {\bigl(\,\sob
{1,2} {} \home\,\bigr)^*} } \,\Bigr) \quad \follows
\cr
{} & \hbox{\kern-20pt} {\norm{{\partial \gat^\p \over \partial t} - {\partial
\gat^\pp \over \partial t} }}{}^2_{\bochnerl {2} {(\,0\,,\,T\,)} {\L {2} {} \home} }
\,\le\, C\,T\, \Bigl(\, {\norm{\excin^\p - \excin^\pp }}^2_{\bochnerl {2}
{(\,0\,,\,T\,)} {\bigl(\,\sob {1,2} {} \home\,\bigr)^*} } &  (3.45)
\cr
{} & & \displaystyle + \,{\norm{\excex^\p - \excex^\pp }}^2_{\bochnerl {2}
{(\,0\,,\,T\,)} {\bigl(\,\sob {1,2} {} \home\,\bigr)^*} } \,\Bigr) \,.
\cr}$$
Combining (3.45) with (3.36), we get finally:
$$\eqalignno{
{} &  {\norm{\gat^\p - \gat^\pp }}^2_{\bochnersob {1,2}
{(\,0\,,\,T\,)} {\L {2} {} \home} }  & (3.46)
\cr
{} & & \displaystyle \le C\, \Bigl(\, {\norm{\excin^\p - \excin^\pp
}}^2_{\bochnerl {2} {(\,0\,,\,T\,)} {\bigl(\,\sob {1,2} {} \home\,\bigr)^*} }  +
{\norm{\excex^\p - \excex^\pp }}^2_{\bochnerl {2} {(\,0\,,\,T\,)} {\bigl(\,\sob
{1,2} {} \home\,\bigr)^*} } \,\Bigr) \,.
\cr}$$
\par
$\bullet$ {\bf Step 6.} Since we already know from (2.12) that $\pottr^\p$ and $\pottr^\pp$ belong even to $\bochnerc {0} {[\,0\,,\,T\,]} {\L {2} {} \home} $, the norm ${\norm{\pottr^\p - \pottr^\pp}}_{\bochnerl {\infty} {(\,0\,,\,T\,)} {\L {2} {} \home} }$ may be replaced by ${\norm{\pottr^\p - \pottr^\pp}}_{\bochnerc {0} {[\,0\,,\,T\,]} {\L {2} {} \home} }$ on the left-hand side of (3.34) . Note that (3.46) implies a bound for ${\norm{\gat^\p - \gat^\pp}}_{\bochnerc {0} {[\,0\,,\,T\,]} {\L {2} {} \home} }$ as well, and the proof is complete.
\medskip
{\it Part B. The FitzHugh-Nagumo model.} If the Rogers-McCulloch model is replaced by the FitzHugh-Nagumo model then the proof can be repeated with some obvious modifications.
\medskip
{\it Part C. The linearized Aliev-Panfilov model.} Let us specify now the linearized Aliev-Panfilov model instead of the Rogers-McCulloch model. Then Proposition 2.4., 2) ensures that $\gat^\p$, $\gat^\pp$ belong still to $\bochnerc {0} {[\,0\,,\,T\,]} {\L {3} {} \home} $.
\medskip
$\bullet$ {\bf Step 1.} Since Proposition 2.4., 2) ensures that $\gat^\p$, $\gat^\pp$ belong still to $\bochnerc {0} {[\,0\,,\,T\,]} {\L {3} {} \home} $, Step 1 from the proof of Theorem 1.1.~can be taken over without alterations.
\medskip
$\bullet$ {\bf Step 2.} The estimates (3.23) ff.~must be replaced as follows:
$$\eqalignno{
{} & \hbox{\kern-20pt}\langle \, {d \over dt} \,\bigl(\, \gat^\p (t) - \gat^\pp(t)\,\bigr) \,,\,
\psi\,\rangle \,=\, -\, \e\,\int_\Omega \bigl(\,\gat^\p (t)-
\gat^\pp(t)\,\bigr)\,\psi\, dx + \,\e\,\kappa\,(a+1) \int_\Omega \bigl(\,\pottr^\p(t)-
\pottr^\pp(t)\,\bigr)\,\psi\, dx & (3.47)
\cr
{} & & \displaystyle  - \,\e\,\kappa \int_\Omega  \bigl(\,\pottr^\p(t) +
\pottr^\pp(t)\,\bigr) \, \bigl(\,\pottr^\p(t) - \pottr^\pp(t)\,\bigr)\,\psi \quad \forall\, \psi \in \L {2} {} \home\,.
\cr}$$
Inserting the feasible test function $\psi = \gat^\p(t) - \gat^\pp(t)$, we obtain
$$\eqalignno{
{} & \hbox{\kern-20pt} {d \over dt}\,\Bigl(\,{\norm{\gat^\p - \gat^\pp}}^2_{\L {2}
{} \home}\,\Bigr)\,\le \, 2\,\e \, {\norm{\gat^\p - \gat^\pp}}^2_{\L {2}
{} \home}  + 2\,\e\,\kappa\,(a+1) \,\Bigl(\, {\norm{\gat^\p - \gat^\pp}}^2_{\L {2}
{} \home} + {\norm{\pottr^\p - \pottr^\pp}}^2_{\L {2} {} \home} \,\Bigr)
\cr
{} & + \,2 \,\e_8 (t) \, \Bigl(\, \int_\Omega \bigl\vert\,\pottr^\p +
\pottr^\pp \,\bigr\vert^4\, dx\,\Bigr)^{1/2} \, \Bigl(\, {\norm{\pottr^\p -
\pottr^\pp}}^4_{\L {4} {} \home} \,\Bigr)^{1/2} + {2\,C \over \e_8 (t)}\,
{\norm{\gat^\p - \gat^\pp}}^2_{\L {2} {} \home}\,.& (3.48)
\cr}$$
Inserting now $\e_8(t) = \e^\p_8 / \bigl(\,1 + {\norm{\pottr^\p (t)}}^2_{\L {4} {} \home} + {\norm{\pottr^\pp(t)}}^2_{\L {4} {} \home}\,\bigr)$ with arbitrary $\e^\p_8 > 0$, we get
$$\eqalignno{
{} & \lldots \le \, 2\, \e \, {\norm{\gat^\p - \gat^\pp}}^2_{\L {2}
{} \home}  + 2\,\e\,\kappa\,(a+1) \,\Bigl(\, {\norm{\gat^\p - \gat^\pp}}^2_{\L {2}
{} \home} + {\norm{\pottr^\p - \pottr^\pp}}^2_{\L {2} {} \home} \,\Bigr) & (3.49)
\cr
{} & + C\,\e^\p_8 \, {\norm{\pottr^\p - \pottr^\pp}}^2_{\sob {1,2} {} \home} + {C \over \e^\p_8 }\, \Bigl(\,1 + {\norm{\pottr^\p}}^2_{\L {4} {} \home} + {\norm{\pottr^\pp}}^2_{\L {4} {} \home}\,\Bigr)\,{\norm{\gat^\p - \gat^\pp}}^2_{\L {2} {} \home} &
\cr}$$
which, after an appropriate choice of $\e^\p_8 > 0$, allows to continue the estimations as above.
\medskip
$\bullet$ {\bf Step 3.} In the case of the linearized Aliev-Panfilov model, we get from (3.29)
$$\eqalignno{
{} & \lldots \le \, C\, \Bigl(\, 1 + {\norm{\Phi_0}}^2_{\L {2} {} \home} + {\norm{W_0}}^{3}_{\L {3} {} \home} + 4\, R^2\,\Bigr)\, & (3.50)
\cr}$$
instead of (3.33).
\medskip
$\bullet$ {\bf Step 4.} This step can be taken over without changes.
\medskip
$\bullet$ {\bf Step 5.} Instead of (3.43), we find by inserting the feasible test function $\psi = (\partial\gat^\p(t)/\partial t) - (\partial\gat^\pp(t)/\partial t)$ into (3.47)
$$\eqalignno{
{} & \hbox{\kern-20pt} {\norm{{\partial \gat^\p \over \partial t} - {\partial
\gat^\pp \over \partial t} }}{}^2_{\L {2} {} \home} \, \le \, C\,\e\,\e^\p_9 \,
{\norm{ {\partial \gat^\p \over \partial t} - {\partial \gat^\pp \over \partial t}
}}{}^2_{\L {2} {} \home}  + {C\,\e \over \e^\p_9}\,{\norm{\gat^\p -\gat^\pp}}^2_{\L
{2} {} \home} & (3.51)
\cr
{} & \quad + \,C\,\e\,\kappa\,(a+1)\,\Bigl(\, \e^\p_{10} \, {\norm{{\partial \gat^\p
\over \partial t} - {\partial \gat^\pp \over \partial t} }}{}^2_{\L {2} {} \home}  +
{1 \over \e^\p_{10}}\,(a+1)\, {\norm{\pottr^\p - \pottr^\pp}}^2_{\L {2} {}
\home}\,\Bigr)
\cr
{} & & \displaystyle + \,C\,\e^\p_{11} \, {\norm{{\partial \gat^\p
\over \partial t} - {\partial \gat^\pp \over \partial t} }}{}^2_{\L {2} {} \home} + {C \over \e^\p_{11} } \, \int_\Omega  \bigl(\,\pottr^\p(t) + \pottr^\pp(t)\,\bigr)^2 \, \bigl(\,\pottr^\p(t) - \pottr^\pp(t)\,\bigr)^2\, dx
\cr
{} & \le \, C\,\bigl(\,\e^\p_9 + \e^\p_{10} + \e^\p_{11}\,\bigr)\, {\norm{ {\partial \gat^\p \over \partial t} - {\partial \gat^\pp \over \partial t} }}{}^2_{\L {2} {} \home}  + {C\over \e^\p_9}\,{\norm{\gat^\p -\gat^\pp}}^2_{\L {2} {} \home} + {C \over \e^\p_{10}}\, {\norm{\pottr^\p - \pottr^\pp}}^2_{\L {2} {}
\home} & (3.52)
\cr
{} & & \displaystyle  + \, {C \over \e^\p_{11}} \,\Bigl(\, \int_\Omega \bigl(\,\pottr^\p(t) + \pottr^\pp(t)\,\bigr)^4 \,dx \,\Bigr)^{1/2}\,\Bigl(\, \int_\Omega \bigl(\,\pottr^\p(t) - \pottr^\pp(t)\,\bigr)^4 \,dx\,\Bigr)^{1/2}
\cr
{} & \le \, C\,\bigl(\,\e^\p_9 + \e^\p_{10} + \e^\p_{11}\,\bigr)\, {\norm{ {\partial \gat^\p \over \partial t} - {\partial \gat^\pp \over \partial t} }}{}^2_{\L {2} {} \home}  + {C\over \e^\p_9}\,{\norm{\gat^\p -\gat^\pp}}^2_{\L {2} {} \home} + {C \over \e^\p_{10}}\, {\norm{\pottr^\p - \pottr^\pp}}^2_{\L {2} {}
\home} & (3.53)
\cr
{} & & \displaystyle  + \, {C \over \e^\p_{11}} \bigl(\, {\norm{\pottr^\p}}^2_{\L {4} {} \home} + {\norm{\pottr^\pp}}^2_{\L {4} {} \home} \,\bigr)\cdot {\norm{\pottr^\p - \pottr^\pp}}^2_{\L {4} {} \home} \,.
\cr}$$
Choosing $\e^\p_9$, $\e^\p_{10}$ and $\e^\p_{11} > 0$ sufficiently small, we obtain
$$\eqalignno{
{} & \hbox{\kern-20pt} {1 \over 2} \,{\norm{{\partial \gat^\p \over \partial t} - {\partial
\gat^\pp \over \partial t} }}{}^2_{\L {2} {} \home}\,\le\,  {C\over \e^\p_9}\,{\norm{\gat^\p -\gat^\pp}}^2_{\L {2} {} \home} + {C \over \e^\p_{10}}\, {\norm{\pottr^\p - \pottr^\pp}}^2_{\L {2} {} \home} & (3.54)
\cr
{} & & \displaystyle + \, {C \over \e^\p_{11}} \bigl(\, {\norm{\pottr^\p}}^2_{\L {4} {} \home} + {\norm{\pottr^\pp}}^2_{\L {4} {} \home} \,\bigr)\cdot {\norm{\pottr^\p - \pottr^\pp}}^2_{\L {4} {} \home} \quad \follows
\cr
{} & \int^T_0 {\norm{{\partial \gat^\p \over \partial t} - {\partial
\gat^\pp \over \partial t} }}{}^2_{\L {2} {} \home}\, dt \,\le\,  {C\over \e^\p_9}\,\int^T_0 {\norm{\gat^\p -\gat^\pp}}^2_{\L {2} {} \home} \, dt \,+\, {C \over \e^\p_{10}}\, \int^T_0 {\norm{\pottr^\p - \pottr^\pp}}^2_{\L {2} {} \home}\, dt & (3.55)
\cr
{} & & \displaystyle +\, \, {C \over \e^\p_{11}} \Bigl(\, \int^T_0 \bigl(\, {\norm{\pottr^\p}}^4_{\L {4} {} \home} + {\norm{\pottr^\pp}}^4_{\L {4} {} \home} \,\bigr)\, dt\,\Bigr)^{1/2} \cdot \Bigl(\,\int^T_0  {\norm{\pottr^\p - \pottr^\pp}}^4_{\L {4} {} \home}\,dt\,\Bigr)^{1/2}\,,
\cr}$$
and the estimates may be continued as before.
\smallskip
$\bullet$ {\bf Step 6.} Now the proof can be finished as above. $\blackbox$
\bigskip
{\bf b) Proof of Theorem 1.2.}
\medskip
The proof of Theorem 1.2.~relies on the complete structural equivalence of the weak bidomain system (2.75)$\,-\,$(2.79) and the reduced bidomain system (2.82)$\,-\,$(2.84). The monodomain form ${\cal M}$ and the bidomain form ${\cal A}$ satisfy analogous norm estimates, the weak solutions of both systems obey the same type of a-priori estimates, and $[\,${\kpt Kunisch/Wagner 13a}$\,]\,$, p.~959 f., Lemma 2.9., yields for arbitrary $\e_0 > 0$ the estimate
$$\eqalignno{
{} &  \bigl\vert\,\langle\, S^\p(t) - S^\pp(t)\,,\,\psi\,\rangle\,\bigr\vert \, & (3.56)
\cr
{} & & \displaystyle \le\, {C \over 2\e_0} \,  \Bigl(\,{\norm{\excin^\p(t) - \excin^\pp (t)}}^2_{\bigl(\,\sob {1,2} {} \home \,\bigr)^\ast} + {\norm{\excex^\p(t)- \excex^\pp(t)}}^2_{\bigl(\,\sob {1,2} {}  \home\,\bigr)^\ast} \,\Bigr) + {3\e_0 \over 4}\, {\norm{\psi}}^2 _{\sob {1,2} {} \home}
\cr}$$
for the difference of the right-hand sides where the constant $C >0$ does not depend on $\e_0$, $\excin$ and $\excex$ and even not on $\Phi_0$ and $W_0$. Consequently, we may carry over Steps 1$\,-\,$6 from the proof of Theorem 1.1. The proof of the estimate for ${\norm{\potex^\p - \potex^\pp }}^2_{\bochnerl {2} {(\,0\,,\,T\,)} {\sob {1,2} {} \home} } $ from the difference of the elliptic equations was not influenced by the error to be corrected. Thus we may take over the respective step from $[\,${\kpt Kunisch/Wagner 13a}$\,]\,$, pp.~971 ff., (2.117)$\,-\,$(2.125). $\blackbox$
\bigskip
{\bf c) Proof of Theorems 1.3.~and 1.4.}
\medskip
The stability estimates from Theorems 1.1.~and 1.2.~yield immediately uniqueness of weak solutions corresponding to right-hand sides and initial values of the assumed regularity. $\blackbox$
\bigskip
{\bf d) Remarks and corollaries.}
\medskip
\baselineskip=17.2pt
Within the previous statements of Theorems 1.1.~and 1.2., we included an estimate for the norm difference\break ${\norm{ \pottr^\p  - \pottr^\pp }}_{\bochnersob {1, q} {(\,0\,,\,T\,)} {\bigl(\,\sob {1,2} {} \home\,\bigr)^*} } $ with certain $q >1$. However, the subsequent analysis of the optimal control problems could be performed without using an estimate of this type. Nevertheless, for sake of completeness we provide its correction here.
\medskip
\baselineskip=14.4pt
{\bf Corollary 3.2.}$\,$\footnote{$^{25)}$}{Additional correction of $[\,${\ninekpt Kunisch/Wagner 12}$\,]\,$, p.~1533, Theorem 3.8.} {\it Under the assumptions of Theorem 1.1., let two weak solutions $(\pottr^\p, \gat^\p)$, $(\pottr^\pp, \gat^\pp)$ of the monodomain system correspond to initial values $\Phi^\p_0 = \Phi^\pp_0 = \Phi_0 \in \L {2} {} \home$, $W^\p_0 = W^\pp_0 = W_0 \in \L {4} {} \home$ and inhomogeneities $\excin^\p$, $\excex^\p$, $\excin^\pp $ and $\excex^\pp \in \L {2} {} \bigl[\, {(\,0\,,\,T\,)}\,,$ $ {\bigl(\,\sob {1,2} {} \home\,\bigr)^*}\,\bigr]\, $, whose norms are bounded by $R > 0$. Then the following estimate holds:
$$\eqalignno{
{} & \hbox{\kern-20pt} {\norm{ \pottr^\p  - \pottr^\pp }}_{\bochnersob {1, 1}
{(\,0\,,\,T\,)} {\bigl(\,\sob {1,2} {} \home\,\bigr)^*} } & (3.57)
\cr
{} & & \displaystyle \le\,C \,\Bigl(\, {\norm{\excin^\p - \excin^\pp }}^2_{\bochnerl
{2} {(\,0\,,\,T\,)}  {\bigl(\,\sob {1,2} {} \home\,\bigr)^*} } + {\norm{\excex^\p - \excex^\pp }}^2_{\bochnerl {2} {(\,0\,,\,T\,)} {\bigl(\,\sob {1,2} {} \home\,\bigr)^*} } \,\Bigr)\,.
\cr}$$
The constant $C > 0$ does not depend on $\excin^\p$, $\excex^\p$, $\excin^\pp$ and $\excex^\pp$ but possibly on $\Omega$, $R$, $\Phi_0$ and $W_0$.}
\medskip
{\bf Corollary 3.3.}$\,$\footnote{$^{26)}$}{Additional correction of $[\,${\ninekpt Kunisch/Wagner 13a}$\,]\,$, p.~959, Theorem 2.7., and $[\,${\ninekpt Kunisch/Wagner 13b}$\,]\,$, p.~1082, Theorem 2.4.} {\it Under the assumptions of Theorem 1.2., let two weak solutions $(\pottr^\p, \potex^\p, \gat^\p)$, $(\pottr^\pp$, $\potex^\pp, \gat^\pp)$ of the bidomain system correspond to initial values $\Phi^\p_0 = \Phi^\pp_0 = \Phi_0 \in \L {2} {} \home$, $W^\p_0 = W^\pp_0 = W_0 \in \L {4} {} \home$ and inhomogeneities $\excin^\p$, $\excex^\p$, $\excin^\pp $ and $\excex^\pp \in \bochnerl {2} {(\,0\,,\,T\,)} {\bigl(\,\sob {1,2} {} \home\,\bigr)^*} $, which satisfy the compatibility conditions $(1.14)$, and whose norms are bounded by $R > 0$. Then the following estimate holds:
$$\eqalignno{
{} & \hbox{\kern-20pt} {\norm{ \pottr^\p  - \pottr^\pp }}_{\bochnersob {1, 1}
{(\,0\,,\,T\,)} {\bigl(\,\sob {1,2} {} \home\,\bigr)^*} } & (5,58)
\cr
{} & & \displaystyle \le\,C \,\Bigl(\, {\norm{\excin^\p - \excin^\pp }}^2_{\bochnerl
{2} {(\,0\,,\,T\,)}  {\bigl(\,\sob {1,2} {} \home\,\bigr)^*} } + {\norm{\excex^\p - \excex^\pp }}^2_{\bochnerl {2} {(\,0\,,\,T\,)} {\bigl(\,\sob {1,2} {} \home\,\bigr)^*} } \,\Bigr)\,.
\cr}$$
The constant $C > 0$ does not depend on $\excin^\p$, $\excex^\p$, $\excin^\pp$ and $\excex^\pp$ but possibly on $\Omega$, $R$, $\Phi_0$ and $W_0$.}
\medskip
{\bf Proof of Corollary 3.2.} First, let us specify within the monodomain system (2.9)$\,-\,$(2.11) the Rogers-McCulloch model. We start with
$$\eqalignno{
{} & \hbox{\kern-20pt} {\Norm{ {\partial \pottr^\p(t)\over \partial t} - {\partial \pottr^\pp (t)\over \partial t}  }}_{\bigl(\,\sob {1,2} {} \home\,\bigr)^*}  \,\le \, {\Norm{ {1 \over 1 + \la}\, \Bigl(\,\la\, \bigl(\, \excin^\p(t) - \excin^\pp(t)\,\bigr) -
\bigl(\,\excex^\p (t)- \excex^\pp(t)\,\bigr)\,\Bigr) }}_{\bigl(\,\sob {1,2} {} \home\,\bigr)^*}  & (3.59)
\cr
{} & & \displaystyle + \, {\Norm{ {\cal M} \bigl(\, \pottr^\p - \pottr^\pp \,,\,\, \cdot\,\,\bigr) }}_{\bigl(\,\sob {1,2} {} \home\,\bigr)^*}  +  {\Norm{ \ioncurr (\pottr^\p, \gat^\p) - \ioncurr(\pottr^\pp, \gat^\pp) }}_{\bigl(\,\sob {1,2} {}
\home\,\bigr)^*}  \quad \follows
\cr}$$
$$\eqalignno{
{} & \hbox{\kern-20pt} \int^T_0 {\Norm{ {\partial \pottr^\p \over \partial t} - {\partial \pottr^\pp \over \partial t}  }}_{\bigl(\,\sob {1,2} {} \home\,\bigr)^*} \, dt \, & (3.60)
\cr
{} & \le\, C \,\Bigl(\, \int^T_0  {\Norm{ {1 \over 1 + \la}\, \Bigl(\,\la\, \bigl(\, \excin^\p(t) - \excin^\pp(t)\,\bigr) -
\bigl(\,\excex^\p (t)- \excex^\pp(t)\,\bigr)\,\Bigr) }}^2_{\bigl(\,\sob {1,2} {} \home\,\bigr)^*}\,dt\,\Bigr)^{1/2} &
\cr
{} &  \hbox{\kern-20pt} + \, C \, \Bigl(\, \int^T_0 \sup_{\psi\,\in\,\sob {1,2} {} \home}\, \bigl\vert\, {\cal M}\bigl(\,
\pottr^\p - \pottr^\pp\,,\, \psi\,\bigr)\,\bigr\vert^2 \,dt\,\Bigr)^{1/2} \!\!
+  C  \int^T_0 {\Norm{ \ioncurr (\pottr^\p, \gat^\p) - \ioncurr(\pottr^\pp, \gat^\pp) }}_{\bigl(\,\sob {1,2} {} \home\,\bigr)^*} \, dt \,,
\cr}$$
estimating the first and second term by using the continuous imbedding $\L {2} {} [\,0\,,\,T\,] \hookrightarrow \L {1} {} [\,0\,,\,T\,]\,$. Now we estimate the three terms on the right-hand side of (3.60) separately. For the first term, we get$\,$\footnote{$^{27)}$}{Note that, in $[\,${\kpt Kunisch/Wagner 12}$\,]\,$, p.~1547 f., (B.44) and (B.47), the forming of the square root has been overlooked.}
$$\eqalignno{
{} & \hbox{\kern-20pt} {\Norm{\, {1 \over 1 + \la}\,\Bigl(\,\la\, \bigl(\, \excin^\p(t) -
\excin^\pp(t)\,\bigr) - \bigl(\,\excex^\p (t)- \excex^\pp(t)\,\bigr)\,\Bigr)
\,}}{}^2_{\bigl(\,\sob {1,2} {} \home\,\bigr)^*} \, & (3.61)
\cr
{} & & \displaystyle \le\, C\,\Bigl(\, {\norm{\excin^\p(t)-
\excin^\pp(t)}}{}^2_{\bigl(\,\sob {1,2} {} \home\,\bigr)^*} + {\norm{\excex^\p(t)-
\excex^\pp(t)}}^2_{\bigl(\,\sob {1,2} {} \home\,\bigr)^*} \,\Bigr) \quad \follows
\cr
{} & \hbox{\kern-20pt} \Bigl(\, \int^T_0 {\Norm{\, {1 \over 1 + \la}\,\Bigl(\,\la\, \bigl(\, \excin^\p(t) -
\excin^\pp(t)\,\bigr) - \bigl(\,\excex^\p (t)- \excex^\pp(t)\,\bigr)\,\Bigr)
\,}}{}^2_{\bigl(\,\sob {1,2} {} \home\,\bigr)^*}  \,dt\,\Bigr)^{1/2} \, & (3.62)
\cr
{} & & \displaystyle \le\, C\, \Bigl(\, {\norm{\excin^\p - \excin^\pp }}{}_{\bochnerl {2}
{(\,0\,,\,T\,)} {\bigl(\,\sob {1,2} {} \home\,\bigr)^*} } + {\norm{\excex^\p -
\excex^\pp }}_{\bochnerl {2} {(\,0\,,\,T\,)} {\bigl(\,\sob {1,2} {}
\home\,\bigr)^*} }\,\Bigr) \,.
\cr}$$
For the second term, we obtain from the continuity of the monodomain bilinear form and $(3.42)$:
$$\eqalignno{
{} & \hbox{\kern-20pt} \bigl\vert\, {\cal M} \bigl(\, \pottr^\p - \pottr^\pp\,,\,
\psi\,\bigr)\,\bigr\vert^2 \,\le\, \gamma^2 \, {\norm{\pottr^\p -
\pottr^\pp}}^2_{\sob {1,2} {} \home} \cdot {\norm{\psi}}^2_{\sob {1,2} {} \home}
\quad \follows \vphantom{\int} & (3.63)
\cr
{} &  \Bigl(\, \int^T_0 \sup_{\lldots}\, \bigl\vert\, {\cal M}\bigl(\, \pottr^\p -
\pottr^\pp\,,\, \psi\,\bigr)\,\bigr\vert^2 \,dt\,\Bigr)^{1/2} \,\le\, C \,\Bigl(\,
{\norm{\pottr^\p - \pottr^\pp}}^2_{\bochnerl {2} {(\,0\,,\,T\,)} {\sob {1,2} {} \home}
} \,\Bigr)^{1/2} \, & (3.64)
\cr
{} & & \displaystyle \le \, C \, \Bigl(\, {\norm{\excin^\p - \excin^\pp
}}{}_{\bochnerl {2} {(\,0\,,\,T\,)} {\bigl(\,\sob {1,2} {} \home\,\bigr)^*} } +
{\norm{\excex^\p - \excex^\pp }}_{\bochnerl {2} {(\,0\,,\,T\,)} {\bigl(\,\sob
{1,2} {} \home\,\bigr)^*} }\,\Bigr)\,.  \quad (3.65)
\cr}$$
In order to estimate the third term, we note first that, by duality, $\L {6/5} {} \home$ is continuously imbedded into $\bigl(\,\sob {1,2} {} \home\,\bigr)^*$ since $\sob {1,2} {} \home$ is continuously imbedded into $\L {6} {} \home$. Consequently, relying on Lemma 3.1., we may write
$$\eqalignno{
{} & \hbox{\kern-20pt} {\norm{\ioncurr (\pottr^\p, \gat^\p ) - \ioncurr(\pottr^\pp,
\gat^\pp) }}_{ \bigl( \,\sob {1,2} {} \home\,\bigr)^*}\,\le\, C\, {\norm{\ioncurr (\pottr^\p, \gat^\p ) - \ioncurr(\pottr^\pp,
\gat^\pp) }}_{ \L {6/5} {} \home} \, \vphantom{\int}
\cr
{} & \le \, C\,{\norm{\,b\,\Bigl(\, (\pottr^\p)^2 + \pottr^\p \,\pottr^\pp + (\pottr^\pp)^2 - (a+1)\,(\pottr^\p + \pottr^\pp) + a\,\Bigr) \,\bigl(\, \pottr^\p - \pottr^\pp\,\bigr)\, }}_{\L {6/5} {} \home} \vphantom{\int} & (3.66)
\cr
{} & & \displaystyle  + \,C\,{\norm{\,(\pottr^\p - \pottr^\pp)\,\gat^\p}}_{\L {6/5} {}
\home} + C\,{\norm{(\gat^\p - \gat^\pp)\, \pottr^\pp}}_{\L {6/5} {} \home} \, = \, J_1 +
J_2 + J_3\,. \quad (3.67)
\cr}$$
Applying to $J_1$ H\"older's inequality with $p_1 = 5/4$ and $p_2 = 5$, we obtain
$$\eqalignno{
{} & \hbox{\kern-20pt} J_1\,=\, C\,\Bigl(\,\int_\Omega \bigl(\, (\pottr^\p)^2 +
\pottr^\p \,\pottr^\pp + (\pottr^\pp)^2 - (a+1)\,(\pottr^\p + \pottr^\pp) +
a\,\bigr)^{6/5} \,\bigl(\, \pottr^\p - \pottr^\pp\,\bigr)^{6/5}\, dx\,\Bigr)^{5/6} & (3.68)
\cr
{} & \hbox{\kern-20pt} \le\, C\, \Bigl(\,\int_\Omega \bigl(\, (\pottr^\p)^2 + \pottr^\p
\,\pottr^\pp + (\pottr^\pp)^2 -\, (a+1)\, (\pottr^\p + \pottr^\pp) + a\,\bigr)^{3/2} \,dx
\,\Bigr)^{2/3} \,\Bigl(\,\int_\Omega \bigl(\,\pottr^\p -
\pottr^\pp\,\bigr)^{6}\,\Bigr)^{1/6} & (3.69)
\cr
{} & \le  \, C\, \bigl(\, 1 + {\norm{\pottr^\p(t)}}^2_{\L {3} {} \home} +
{\norm{\pottr^\pp(t)}}^2_{\L {3} {} \home} \,\bigr) \cdot{\norm{\pottr^\p -
\pottr^\pp}}_{\sob {1,2} {} \home}\,. \vphantom{\int} & (3.70)
\cr}$$
Further, using again H\"older's inequality with $p_1 = 5/4$ and $p_2 = 5$, we get
$$\eqalignno{
{} & \hbox{\kern-20pt} J_2\,=\, C\,\Bigl(\, \int_\Omega \bigl\vert\,\pottr^\p -
\pottr^\pp\,\bigr\vert^{6/5}\, \bigl\vert\,\gat^\p\,\bigr\vert^{6/5}\, dx\,\Bigr)^{5/6}  & (3.71)
\cr
{} & \le \, C\,{\norm{\gat^\p (t)}}_{\L {3/2} {} \home} \,{\norm{\pottr^\p - \pottr^\pp}}_{\L {6} {} \home} \le \, C\,{\norm{\gat^\p}}_{\bochnerc {0} {[\,0\,,\,T\,]} {\L {2} {} \home} }\,{\norm{\pottr^\p - \pottr^\pp}}_{\L {6} {} \home} \,.& (3.72)
\cr}$$
By (2.13), ${\norm{\gat^\p}}_{\bochnerc {0} {[\,0\,,\,T\,]} {\L {2} {} \home} }$ is uniformly bounded by the norms of the initial data and the bound $R$ of the norms of the inhomogeneities. Consequently, we arrive at
$$\eqalignno{
{} & J_2 \le C\,{\norm{\pottr^\p - \pottr^\pp}}_{\sob {1,2} {} \home} \,. & (3.73)
\cr}$$
Finally, through application of H\"older's inequality with $p_1 = 5/3$ and $p_2 = 5/2$, we will estimate
$$\eqalignno{
{} & \hbox{\kern-20pt} J_3 \,=\, C\,\Bigl(\, \int_\Omega \bigl(\, \gat^\p - \gat ^\pp\,\bigr)^{6/5}\, (\pottr^\pp)^{6/5}\, dx\,\Bigr)^{5/6} \le \,C \, {\norm{\gat^\p - \gat^\pp}}_ {\L {2} {} \home} \cdot {\norm{\pottr^\pp}}_{\L {3} {}
\home}  \,. \vphantom{\int} & (3.74)
\cr}$$
Summing up, the estimates (3.70), (3.73) and (3.74) imply for the third term in (3.61)
$$\eqalignno{
{} & \hbox{\kern-20pt} C \int^T_0 {\Norm{ \ioncurr (\pottr^\p, \gat^\p) - \ioncurr(\pottr^\pp, \gat^\pp) }}_{\bigl(\,\sob {1,2} {} \home\,\bigr)^*} \, dt \le \, C \int^T_0 \Bigl(\, \bigl(\, 1 + {\norm{\pottr^\p(t)}}^2_{\L {3} {} \home}& (3.75)
\cr
{} & & \displaystyle +\, {\norm{\pottr^\pp(t)}}^2_{\L {3} {} \home} \,\bigr) \cdot{\norm{\pottr^\p - \pottr^\pp}}_{\sob {1,2} {} \home} + {\norm{\gat^\p - \gat^\pp}}_ {\L {2} {} \home} \cdot {\norm{\pottr^\pp}}_{\L {3} {} \home}  \,\Bigr)\, dt\,,
\cr}$$
and the combination of (3.60), (3.62), (3.65) and (3.75) yields
$$\eqalignno{
{} & \hbox{\kern-20pt} \int^T_0 {\Norm{ {\partial \pottr^\p \over \partial t} - {\partial \pottr^\pp \over \partial t}  }}_{\bigl(\,\sob {1,2} {} \home\,\bigr)^*} \, dt \, & (3.76)
\cr
{} & \le \, C\, \Bigl(\, {\norm{\excin^\p - \excin^\pp }}{}_{\bochnerl {2}
{(\,0\,,\,T\,)} {\bigl(\,\sob {1,2} {} \home\,\bigr)^*} } + {\norm{\excex^\p - \excex^\pp }}_{\bochnerl {2} {(\,0\,,\,T\,)} {\bigl(\,\sob {1,2} {} \home\,\bigr)^*} }\,\Bigr)
\cr
{} & + \, C \, \Bigl(\,\int^T_0 \bigl(\, 1 + {\norm{\pottr^\p(t)}}^4_{\L {3} {} \home} + {\norm{\pottr^\pp(t)}}^4_{\L {3} {} \home} \,\bigr)\, dt\,\Bigr)^{1/2} \cdot \Bigl(\, \int^T_0 {\norm{\pottr^\p - \pottr^\pp}}^2_{\sob {1,2} {} \home}\,dt\,\Bigr)^{1/2}
\cr
{} & + \, C \,\Bigl(\, \int^T_0 {\norm{\pottr^\pp}}^2_{\L {3} {} \home}\, dt\,\bigr)^{1/2} \cdot \Bigl(\,\int^T_0 {\norm{\gat^\p - \gat^\pp}}^2_ {\L {2} {} \home} \,dt\,\Bigr)^{1/2}  \,.
\cr}$$
Observe that
$$\eqalignno{
{} & \Bigl(\,\int^T_0 {\norm{\pottr^\p(t)}}^4_{\L {3} {} \home}\, dt\,\Bigr)^{1/2} \le\, C\,{\norm{\pottr^\p}}^2_{\L {4} {} (\Omega_T)} \,; & (3.77)
\cr
{} & \Bigl(\,\int^T_0 {\norm{\pottr^\pp(t)}}^4_{\L {3} {} \home}\, dt\,\Bigr)^{1/2} \le\, C\,{\norm{\pottr^\pp}}^2_{\L {4} {} (\Omega_T)} \,; & (3.78)
\cr
{} & \Bigl(\, \int^T_0 {\norm{\pottr^\pp}}^2_{\L {3} {} \home}\, dt\,\bigr)^{1/2} = \,C {\norm{ \pottr^\pp}}_{\bochnerl {2} {(\,0\,,\,T\,)} {\L {3} {} \home} } \le \, C\,{\norm{\pottr^\pp}}_{\L {4} {} (\Omega_T)} & (3.79)
\cr}$$
are bounded by (2.13) and the assumption about the norms of $\excin$ and $\excex$. Consequently, assembling (3.36), (3.42) and (3.76), we arrive at the desired estimate
$$\eqalignno{
{} & \hbox{\kern-20pt} \int^T_0 {\Norm{ {\partial \pottr^\p \over \partial t} - {\partial \pottr^\pp \over \partial t}  }}_{\bigl(\,\sob {1,2} {} \home\,\bigr)^*} \, dt \,\le \, C\, \Bigl(\, {\norm{\excin^\p - \excin^\pp }}{}_{\bochnerl {2} {(\,0\,,\,T\,)} {\bigl(\,\sob {1,2} {} \home\,\bigr)^*} }  & (3.80)
\cr
{} & + {\norm{\excex^\p - \excex^\pp }}_{\bochnerl {2} {(\,0\,,\,T\,)} {\bigl(\,\sob {1,2} {} \home\,\bigr)^*} }\,\Bigr) \,+ C\,\Bigl(\, {\norm{\pottr^\p - \pottr^\pp}}_{\bochnerl {2} {(\,0\,,\,T\,)} {\sob {1,2} {} \home} } + {\norm{\gat^\p - \gat^\pp}}_ {\L {2} {} (\Omega_T)}\,\Bigr)
\cr
{} & \hbox{\kern-20pt} \follows \,\,\, {\norm{ \pottr^\p  - \pottr^\pp }}_{\bochnersob {1, 1} {(\,0\,,\,T\,)} {\bigl(\,\sob {1,2} {} \home\,\bigr)^*} } & (3.81)
\cr
{} & & \displaystyle \le \, C\, \Bigl(\, {\norm{\excin^\p - \excin^\pp }}{}_{\bochnerl {2}
{(\,0\,,\,T\,)} {\bigl(\,\sob {1,2} {} \home\,\bigr)^*} } + {\norm{\excex^\p - \excex^\pp }}_{\bochnerl {2} {(\,0\,,\,T\,)} {\bigl(\,\sob {1,2} {} \home\,\bigr)^*} }\,\Bigr) \,.
\cr}$$
\par
If the Rogers-McCulloch model is replaced by the FitzHugh-Nagumo or the linearized Aliev-Panfilov model, the estimates hold accordingly. Note that, in the latter case, the estimation of $J_2$ in $(3.71)\,-\,(3.73)$ is still possible since, by Proposition 2.4., 2), ${\norm{\gat^\p}}_{\bochnerc {0} {[\,0\,,\,T\,]} {\L {3} {} \home} }$ is uniformly bounded by the norms of the initial data and the bound $R$ of the norms of the inhomogeneities. $\blackbox$
\medskip
{\bf Remark 3.4.} In the case of the Rogers-McCulloch or the FitzHugh-Nagumo model, Theorems 1.1.$\,-\,$1.4.~and Corollaries 3.2.$\,-\,$3.3.~remain valid if $W_0$ belongs to $\L {8/3} {} \home$ instead of $\L {4} {} \home$. Indeed, in order to keep the estimation (3.11)$\,-\,$(3.16) true, it suffices to work with $\gat^\p \in \bochnerl {\infty} {(\,0\,,\,T\,)} {\L {8/3} {} \home} $. $W_0 \in \L {8/3} {} \home$, $(2.13)$ and $(2.17)$ imply
$$\eqalignno{
{} & \hbox{\kern-20pt} \int_\Omega \bigl\vert \,\gat (t)\,\bigr\vert^{8/3} \, dx \,\le\, C\,{\norm{W_0}}^{8/3}_{\L {8/3} {} \home} + C\,\e\,\kappa \,{\norm{\pottr }}^{8/3}_{\L {8/3} {} (\Omega_T)} \le\, C\,{\norm{W_0}}^{8/3}_{\L {8/3} {} \home} + C\,\bigl(\, 1 + {\norm{\pottr }}^4_{\L {4} {} (\Omega_T)}\,\bigr)  & (3.82)
\cr
{} & \hbox{\kern-20pt} \le\, C\,\Bigl(\,1 + {\norm{\Phi_0}}^2_{\L {2} {} \home} + {\norm{W_0}}^{8/3}_{\L {8/3} {} \home} \!+  {\norm{\excin^\p}}^2_{\bochnerl {2} {(\,0\,,\,T\,)} {\bigl(\,\sob {1,2} {} \home\,\bigr)^*} }  +  {\norm{\excex^\p}}^2 _{\bochnerl {2} {(\,0\,,\,T\,)} {\bigl(\,\sob {1,2} {} \home\,\bigr)^*} }  \,\Bigr) \,, & (3.83)
\cr}$$
which confirms that $\gat$ belongs to $\bochnerl {\infty} {(\,0\,,\,T\,)} {\L {8/3} {} \home} $.
\medskip
{\bf Remark 3.5.} In the case of the linearized Aliev-Panfilov model, Theorems 1.1.$\,-\,$1.4.~and Corollaries 3.2.$\,-\,$3.3.~remain valid if $W_0$ belongs to $\L {3} {} \home$ instead of $\L {4} {} \home$.
\bigskip\medskip
{\kapitel 4.~Detailed corrections within the previous papers.}
\bigskip
In the following, we report in detail the corrections to be made in $[\,${\kpt Kunisch/Wagner 12}$\,]\,$, $[\,${\kpt Kunisch/ Wag\-ner 13a}$\,]\,$, $[\,${\kpt Kunisch/Wagner 13b}$\,]$ and $[\,${\kpt Kunisch/Nagaiah/Wagner 11}$\,]\,$.
\bigskip
{\bf a) Corrections within $[\,${\kpt Kunisch/Wagner 12}$\,]\,$.}
\medskip
1) $[\,${\kpt Kunisch/Wagner 12}$\,$, p.~1528 f., Theorem 2.2., Lemma 2.3.~and Theorem 2.4.$\,]\,$. These assertions remain true without changes.
\medskip
2) $[\,${\kpt Kunisch/Wagner 12}$\,$, p.~1529, Theorem 2.5.$\,]\,$. In view of Theorem 1.3.~above, the assertion remains true for the mono\-domain system with the Rogers-McCulloch, FitzHugh-Nagumo or the linearized Aliev-Panfilov model.
\medskip
3) $[\,${\kpt Kunisch/Wagner 12}$\,$, p.~1529 f.$\,]$ The analytical framework for the analysis of the optimal control problem $(\P)$ is determined by the additional regularity requirements within the stability estimate and the existence theorem for the adjoint parabolic equation. In view of Theorem 1.5.~above, no changes have to be made.
\medskip
4) $[\,${\kpt Kunisch/Wagner 12}$\,$, p.~1530 f., Proposition 3.1, Proposition 3.2.~and Theorem 3.3.$\,]\,$. The assertions remain true without changes.
\medskip
5) $[\,${\kpt Kunisch/Wagner 12}$\,$, p.~1532, Theorem 3.5., Corollary 3.6.~and 3.7.$\,]\,$. These assertions remain true without changes. However, in order to ensure the claimed regularity of the multiplier $P_1$, we must argue with Theorem 1.7.~instead of Theorem 1.6.
\medskip
{\bf Lemma 4.1.} {\it Under the assumptions of $[\,${\kpt Kunisch/Wagner 12}$\,$, {\rm p.~1532, Theorem 3.5.}$\,]\,$, the adjoint va\-riable $P_1$ admits still $\bochnerl {4} {(\,0\,,\,T\,)} {\sob {1,2} {} \home} $-regularity.}
\medskip
{\bf Proof.} By $[\,${\kpt Kunisch/Wagner 12}$\,$, p.~1532, Theorem 3.5.$\,]\,$, we obtain $P_1 \in \bochnerl {2} {(\,0\,,\,T\,)} {\sob {2,2} {} \home } \cap \, \bochnersob {1,2} {(\,0\,,\,T\,)} {\L {2} {} \home} $. Instead of the Aubin-Dubinskij lemma (Theorem 1.6.), we apply Theorem 1.7., 1) to the spaces $\X_0 = \sob {2,2} {} \home$, $\X = \sob {1,2} {} \home$ and $\X_1 = \L {2} {} \home$. Since the imbedding $\sob {2,2} {} \home \hookrightarrow \L {2} {} \home$ is compact and the norms satisfy the interpolation inequality$\,$\footnote{$^{28)}$}{$[\,${\ninekpt Adams/Fournier 07}$\,]\,$, p.~135, Theorem 5.2., (3), with $n=3$, $j =1$, $p = 2$ and $m = 2$.}
$$\eqalignno{
{} & {\norm{\psi}}_{\sob {1,2} {} \home} \,\le\, C\, {\norm{\psi}}^{1/2}_{\sob {2,2} {} \home} \cdot {\norm{\psi}}^{1/2} _{\L {2} {} \home} \quad \forall\, \psi \in \sob {2,2} {} \home\,, & (4.1)
\cr}$$
we observe with $p = p^\p =2$ and $\vartheta = 1/2$ that $\Delta(p, p^\p, \vartheta) = 0$. Consequently, $P_1$ belongs to all spaces $\bochnerl {q} {(\,0\,,\,T\,)} {\sob {1,2} {} \home} $ with $1 \le q < \infty$. $\blackbox$
\medskip
6) $[\,${\kpt Kunisch/Wagner 12}$\,$, p.~1533, Theorem 3.8.$\,]\,$. The stability estimate must be replaced by Theorem 1.1.~and Corollary 3.2.~above. Note that the inequality $[\,${\kpt Kunisch/Wagner 12}$\,$, p.~1533, (3.28)$\,]\,$, which could not be maintained, has not been used in the following derivations.
\medskip
7)  $[\,${\kpt Kunisch/Wagner 12}$\,$, p.~1533, Theorem 3.9.$\,]\,$. The theorem remains valid without changes.
\medskip
8) $[\,${\kpt Kunisch/Wagner 12}$\,$, p.~1533, Theorem 3.10.$\,]\,$. Within the formulation of the parabolic existence theorem, there is a transcription error to be corrected. Namely, assumption (c) has to be replaced by$\,$\footnote{$^{29)}$}{Cf.~again $[\,${\ninekpt Lady\v zenskaja/Solonnikov/Ural'ceva 88}$\,]\,$, p.~180, Remark 6.3.}
\smallskip
(c)$^\p$ $\quad a_0 \in \bochnerl {r_2} {(\,0\,,\,T\,)} {\L {q_2} {} \home} $ for some $1 \le r_2 < \infty$, $2 < q_2 \le \infty$ satisfying $\displaystyle {1\over r_2} + {n\over 2 q_2} = 1$.
\medskip
9) $[\,${\kpt Kunisch/Wagner 12}$\,$, pp.~1534$\,-\,$1537, Proof of Theorem 3.9.$\,]\,$. The proof with the necessary corrections is repeated here.
\medskip
Under the assumptions of the theorem, the adjoint equations read as follows:
$$\eqalignno{
{} & \hbox{\kern-20pt} - {\partial P_1 \over\partial t} - \nabla \! \cdot \bigl(\, {\la \over 1 + \la }\,\condin\,\nabla P_1\,\bigr) + \Bigl(\, 3\,b\,
(\,\pottrhat\,)^2 \!- 2\,(a+1)\,b\,\pottrhat + a\,b + \gathat\,\Bigr)\, P_1 \,= \,\e\,\kappa\, P_2 - {\partial r\over\partial \varphi} (\pottrhat, \gathat)\,; & (4.2)
\cr
{} & \hbox{\kern-20pt}  - {\partial P_2 \over \partial t} + \e \, P_2\, =\, - \pottrhat \, P_1  - {\partial r\over \partial w} (\pottrhat, \gathat) \,. & (4.3)
\cr}$$
\par
$\bullet$ {\bf Step 1.} {\it Improved regularity of $\pottrhat$ and $\gathat$.} From Proposition 2.1., 2) and Theorem 1.5., 1), we see that
$$\eqalignno{
{} & \gathat (x,t) \,=\, W_0(x) \, e^{-\e t} + \e\,\kappa\, \int^t_0 \pottrhat(x,\tau)\,e^{\e (\tau-t)}\, d\tau & (4.4)
\cr}$$
belongs to $\bochnerc {0} {[\,0\,,\,T\,]} {\L {4} {} \home} $ while $\pottrhat$ gains at least $\bochnerl {4} {(\,0\,,\,T\,)} {\sob {1,2} {} \home} $-regularity.$\,$\footnote{$^{30)}$}{The derivation $[\,${\kpt Kunisch/Wagner 12}$\,$, p.~1534, (3.39)$\,-\,$(3.41)$\,]\,$, holds wrong.}
\medskip
$\bullet$ {\bf Step 2.} {\it For any $\widetilde P_1 \in \L {4} {} (\Omega_T)$, the terminal problem for the adjoint ODE admits a unique (weak or strong) solution $P_2 \in \bochnerc {1} {(\,0\,,\,T\,)} {\L {2} {} \home} \,\cap\, \bochnerc {0} {[\,0\,,\,T\,]} {\L {2} {} \home} $.} It is obvious that the problem
$$\eqalignno{
{} & - {\partial P_2 \over \partial t} + \e \, P_2\, =\, - \pottrhat \, \widetilde P_1  - {\partial r\over \partial w} (\pottrhat, \gathat) \quad \fforall\, (x,t) \in \Omega_T\,, \,\,\, P_2(x, T) \equiv 0 & (4.5)
\cr}$$
admits the unique solution
$$ P_2(x,t) \,=\, - \int^T_t \Bigl(\, \pottrhat \widetilde P_1 + {\partial r \over \partial w} (\pottrhat, \gathat)\,\Bigr)\, e ^{\e\,(t - \tau)}\, d\tau \,,\eqno (4.6)
$$
which is continuous in time on $[\,0\,,\,T\,]$ and even differentiable in time on $(\,0\,,\,T\,)$. In order to confirm the integrability with respect to $x$, we estimate
$$\eqalignno{
{} & \int_\Omega \bigl(\,\pottrhat(t)\,\widetilde P_1(t)\,\bigr)^2\, dx \,\le\, \Bigl(\, \int_\Omega \vert\,\pottrhat(t)\,\vert^4\, dx\,\Bigr)^{1/2} \,\Bigl(\, \int_\Omega \vert\,\widetilde P_1(t) \,\vert^{4}\, dx\,\Bigr)^{1/2} & (4.7)
\cr}$$
where the right-hand side is finite due to the continuous imbedding $\pottrhat(t) \in \sob {1,2} {} \home \hookrightarrow \L {4} {} \home$. Consequently, $P_2$ belongs to the space $\bochnerc {1} {(\,0\,,\,T\,)} {\L {2} {} \home} \,\cap\, \bochnerc {0} {[\,0\,,\,T\,]} {\L {2} {} \home} $.
\medskip
$\bullet$ {\bf Step 3.} {\it For any $\widetilde P_2 \in \L {2} {} (\Omega_T)$, the terminal-boundary value problem for the parabolic adjoint equation admits a unique weak solution $P_1 \in \bochnerl {2} {(\,0\,,\,T\,)} {\sob {2,2} {} \home} \,\cap\, \bochnersob {1,2} {(\,0\,,\,T\,)} {\L {2} {} \home} $.} In order to confirm this claim, we must check whether the assumptions of $[\,${\kpt Kunisch/Wagner 12}$\,$, p.~1533, Theorem 3.10.$\,]$ are satisfied. Concerning (a), (b), (d), (e) and (f), the arguments from $[\,${\kpt Kunisch/Wagner 12}$\,$, p.~1534 f.$\,]$ can be maintained. In view of the regularity discussion in Step 1, the term
$$\eqalignno{
{} & a_0 (x,t) = 3\,b\, (\,\pottrhat\,)^2 \!- 2\,(a+1)\,b\,\pottrhat + a\,b + \gathat & (4.8)
\cr}$$
satisfies condition (c)$^\p$ with $n = 3$, $r_2 = 2$ and $q_2 = 3$. Consequently, the application of $[\,${\kpt Kunisch/Wagner 12}$\,$, p.~1533, Theorem 3.10.$\,]$ within the proof of $[\,${\kpt Kunisch/Wagner 12}$\,$, p.~1533, Theorem 3.9.$\,]$ is still justified.
\medskip
$\bullet$ {\bf Step 4.} {\it For two functions $P^\p_1$, $P^\pp_1 \in \L {4} {} (\Omega_T)$, the corresponding solutions of the terminal problem for the adjoint ODE satisfy$\,$\footnote{$^{31)}$}{$[\,${\kpt Kunisch/Wagner 12}$\,$, p.~1535, (3.47) ff.$\,]\,$: Again here and in the following, it was claimed in error that ${\norm{\pottrhat}}_{\L {4} {} \home}$ is essentially bounded.}
$$\eqalignno{
{} & {\norm{P^\p_2(t) - P^\pp_2(t)}}^2_{\L {2} {} \home} \,\le \, C \cdot \int^T_t {\norm{P^\p_1(\tau) - P^\pp_1(\tau)}}^2_{\L {4} {} \home} \, d\tau\,. & (4.9)
\cr}$$}%
Applying Jensen's integral inequa\-lity and H\"older's inequality, we may argue that
$$\eqalignno{
{} & \hbox{\kern-20pt}\int_\Omega \vert\, P^\p_2  - P^\pp_2\,\vert^2\, dx \,\le\, C \int_\Omega \Bigl(\, \int^T_t \vert\,\pottrhat\,\vert \cdot \vert \,P^\p_1 - P^\pp_1\,\vert\, d\tau \,\Bigr)^2 \, dx \, & (4.10)
\cr
{} & \le\, C \int_\Omega \Bigl(\, \int^T_0 \vert\,\pottrhat\,\vert^2 \,d\tau\,\Bigr)\cdot \Bigl(\, \int^T_t \vert \,P^\p_1 - P^\pp_1\,\vert^2\, d\tau \,\Bigr)\, dx \,\le\, C\, {\norm{\pottrhat}}^2_{\L {4} {} (\Omega_T)} \cdot {\norm{P^\p_1 - P^\pp_1}}^2_{\bochnerl {4} {(\,t\,,\,T\,)} {\L {4} {} \home} }\,.
\cr}$$
\par
$\bullet$ {\bf Step 5.} {\it For two functions $P^\p_2$, $P^\pp_2 \in \L {2} {} (\Omega_T)$, the corresponding solutions of the terminal-boundary value problem for the parabolic adjoint equation satisfy
$$\eqalignno{
{} & {\norm{P^\p_1(t) - P^\pp_1(t)}}^2_{\L {4} {} \home} \,\le\, C \cdot \int^T_t {\norm{ P^\p_2(\tau) - P^\pp_2(\tau)}}^2_{\L {2} {} \home} \, d\tau\,. & (4.11)
\cr}$$}%
Applying $[\,${\kpt Kunisch/Wagner 12}$\,$, p.~1533, Theorem 3.10.$\,]$ to the difference of the linear parabolic equations determining $P^\p_1$ and $P^\pp_1$, we get the a-priori estimate$\,$\footnote{$^{32)}$}{Cf.~$[\,${\kpt Evans 98}$\,]\,$, p.~287, Theorem 3.}
$$\eqalignno{
{} & \hbox{\kern-20pt} C\,\Bigl(\, \int^T_t {\norm{P^\p_2 - P^\pp_2}}^2_{\L {2} {} \home} \, d\tau \,\Bigr)^{1/2} \ge\, {\norm{P^\p_1 - P^\pp_1}}_{\bochnerc {0} {[\, t \,,\,T\,]} {\sob {1,2} {} \home} } \,\ge \, {\norm{ P^\p_1 - P^\pp_1}}_{\bochnerl {4} {(\,t\,,\,T\,)} {\L {4} {} \home} }\,. & (4.12)
\cr}$$
\par
$\bullet$ {\bf Step 6.} {\it Application of Banach's fixed point theorem.} We consider the operator$\,$\footnote{$^{33)}$}{Here we follow still $[\,${\ninekpt Veneroni 09}$\,]\,$, p.~866.}
$$  {\cal I}\,\colon \,\, \Bigl(\, \bochnerl {2} {(\,0\,,\,T\,)} {\L {4} {} \home} \times \L {2} {} (\Omega_T)\,\Bigr) \to \Bigl(\, \L {2} {} \bigl[\,(\,0\,,\,T\,), \L {4} {} \home\,\bigr] \times \L {2} {} (\Omega_T) \,\Bigr) \,, \eqno (4.13)
$$
which assigns to a given pair $(P_1, P_2)$ the new pair $({\cal I}P_1, {\cal I} P_2)$ arising from the solution ${\cal I}P_2$ of the adjoint ODE after insertion of $P_1$ and the solution of the adjoint parabolic problem after insertion of ${\cal I} P_2$. Let us prove now the contractivity of this operator. We start with two pairs $(P^\p_1, P^\p_2)$, $(P^\pp_1, P^\pp_2) \in \L {2} {} \bigl[\, {(\,0\,,\,T\,)}\,,$ $
{\L {4} {} \home}\,\bigr] \times \L {2} {} (\Omega_T)$. From (4.9) and (4.11), we get
$$\eqalignno{
{} & {\norm{{\cal I}P^\p_1 - {\cal I} P^\pp_1}}^2_{\bochnerl {4} {(\,t\,,\,T\,)} {\L {4} {} \home} } \,\le\, C\, \int^T_t {\norm{{\cal I} P^\p_2(\tau) - {\cal I} P^\pp_2(\tau)}}^2_{\L {2} {} \home}\, d\tau\, & (4.14)
\cr
{} & \le\, C \,  \int^T_t {\norm{\pottrhat}}^2_{\L {4} {} (\Omega_T)} \cdot {\norm{P^\p_1 - P^\pp_1}}^2_{\bochnerl {4} {(\,\tau\,,\,T\,)} {\L {4} {} \home} } \, d\tau \,. & (4.15)
\cr}$$
Defining the functions
$$ f(t) = {\norm{{\cal I}P^\p_1  - {\cal I} P^\pp_1}}^2_{\bochnerl {4} {(\,t\,,\,T\,)} {\L {4} {} \home} } \,\,\,\hbox{and} \,\,\, \widetilde f(t) =  {\norm{P^\p_1 - P^\pp_1}}^2_{\bochnerl {4} {(\,t\,,\,T\,)} {\L {4} {} \home} } \,, \eqno (4.16)
$$
this inequality reads as
$$\eqalignno{
{} & 0 \,\le\, f(t) \,\le\, C\,\int^T_t \widetilde f(\vartheta)\, d\vartheta \,\,\, \follows \,\,\, \int^T_0 e^{\la_1\,t} \, f(t)\, dt\, \le \, C \cdot \int^T_0 e^{\la_1 \, t}\,\Bigl(\,\int^T_t \widetilde f(\vartheta)\, d\vartheta\,\Bigr)\, dt  & (4.17)
\cr
{} & =\, C\,\Bigl[\, {1 \over \la_1}\, e^{\la_1\, t} \cdot \int^T_t \widetilde f(\vartheta)\, d\vartheta\,\Bigr]^T_0 + C \,\int^T_0 {1 \over \la_1}\, e^{\la_1\, t}\, \widetilde f(t)\, dt & (4.18)
\cr
{} & =\, {C \over \la_1}\,\Bigl(\, \int^T_0 e^{\la_1\,t} \widetilde f(t)\, dt - \int^T_0 \widetilde f (\vartheta)\, d\vartheta\,\Bigr) \,\le \, {C \over \la_1}\,\int^T_0 e^{\la_1\,t} \widetilde f(t)\, dt & (4.19)
\cr}$$
since the second member within the brackets is positive. If a sufficiently large $\la_1 > C$ is fixed, we get the relation
$$\eqalignno{
{} & \lim_{N \to \infty}  \int^T_0 e^{\la_1 \,t} \cdot {\norm{{\cal I}^N \,P^\p_1  - {\cal I}^N \, P^\pp_1}}^2_{\bochnerl {4} {(\,t\,,\,T\,)} {\L {4} {} \home} } \, dt \,=\,0\,, & (4.20)
\cr}$$
implying
$$\eqalignno{
{} & {\norm{{\cal I}^N \,P^\p_1  - {\cal I}^N \, P^\pp_1}}^2_{\bochnerl {4} {(\,t\,,\,T\,)} {\L {4} {} \home} } \to 0 & (4.21)
\cr}$$
for almost all $0 \le t \le T$. Consequently, the sequence $\folge{{\cal I}^N \,P^\p_1  - {\cal I}^N \, P^\pp_1}$ converges in $\bochnerl {2} {(\,0\,,\,T\,)} {\L {4} {} \home} $-norm to the zero function, and the operator ${\cal I} $ is contractive with respect to its first component on this space. For the contractivity with respect to the second component, the arguments from $[\,${\kpt Kunisch/Wagner 12}, p.~1536 f., (3.65) ff.$\,]$ may be repeated. The proof is complete. $\blackbox$
\medskip
10) $[\,${\kpt Kunisch/Wagner 12}$\,$, p.~1537, Remark (1)$\,]\,$. The remark holds true since, in the case of the FitzHugh-Nagumo model, the variable coefficient in Step 3 above reads as
$$\eqalignno{
{} &  a_0 (x,t) = 3\,(\,\pottrhat\,)^2 \!- 2\,(a+1)\,\pottrhat + a\,. & (4.22)
\cr}$$
By Theorem 1.5., 1), this function belongs to $\bochnerl {2} {(\,0\,,\,T\,)} {\L {3} {} \home} $ as well.
\medskip
11) $[\,${\kpt Kunisch/Wagner 12}$\,$, p.~1537, Remark (2)$\,]\,$. In the case of the linearized Aliev-Panfilov model, the application of $[\,${\kpt Ku\-nisch/Wagner 12}$\,$, p.~1533, Theorem 3.10.$\,]$ is still justified. Indeed, from Proposition 2.4., 2), we get $\gathat \in \bochnerc {0} {[\,0\,,\,T\,]} {\L {3} {} \home} $, and $a_0$ belongs still to $\bochnerl {2} {(\,0\,,\,T\,)} {\L {3} {} \home} $ as required in assumption (c)$^\p$. Assumption (d) about the right-hand side can be satisfied as follows: It holds that
$$\eqalignno{
{} & \int^T_0 \int_\Omega \bigl\vert \,\pottrhat\, {P_2}\,\bigr\vert^2\, dx\, dt \,\le\, \int^T_0 {\norm{\pottrhat (t)}}^2_{\L {6} {} \home} \cdot {\norm{ P_2(t)}}^2_{\L {3} {} \home}\, dt \,\le\, C\,{\norm{\pottrhat}}^2_{\bochnerl {2} {(\,0\,,\,T\,)} {\L {6} {} \home} } & (4.23)
\cr}$$
since $\pottrhat \in \bochnerl {4} {(\,0\,,\,T\,)} {\L {6} {} \home} $ by Theorem 1.5., 1), and $\widetilde P_1$, $P_1 \in \L {6} {} (\Omega_T)$ by Lemma 4.1., which implies even $P_2 \in \bochnerc {0} {[\,0\,,\,T\,]} {\L {3} {} \home} $.
\medskip
12) $[\,${\kpt Kunisch/Wagner 12}$\,$, pp.~1537$\,-\,$1541, Proof of Theorem 3.5.~with Lemma 3.11.~and 3.12., Proofs of Corollaries 3.6.~and 3.7.$\,]\,$. All arguments remain true without changes.
\medskip
13) $[\,${\kpt Kunisch/Wagner 12}$\,$, p.~1542, Propositions A.2, A.3, A.4 and Theorem A.5$\,]\,$. The citations are correct.
\medskip
14) $[\,${\kpt Kunisch/Wagner 12}$\,$, p.~1542, Theorem A.6$\,]\,$. The citation of the theorem is erroneous and must be replaced by Theorem 1.6.~above.
\medskip
15) $[\,${\kpt Kunisch/Wagner 12}$\,$, p.~1542 ff., Appendix B$\,]$ The proof of the stability estimate must be replaced by the proofs of Theorem 1.1.~and Corollary 3.2.~above. The basic error was the inequality $[\,${\kpt Kunisch/Wagner 12}$\,$, p.~1544, (B.14)$\,]\,$, which holds not true.
\bigskip
{\bf b) Corrections within $[\,${\kpt Kunisch/Wagner 13a}$\,]\,$.}
\medskip
1) $[\,${\kpt Kunisch/Wagner 13a}$\,$, pp.~956$\,-\,$958, Theorem 2.4., 2.5.~and 2.6.$\,]\,$. These assertions remain true without changes.
\medskip
2) $[\,${\kpt Kunisch/Wagner 13a}$\,$, p.~959, Theorem 2.7.$\,]\,$. The stability estimate must be replaced by Theorem 1.2.~and Corollary 3.3.~above.
\medskip
3) $[\,${\kpt Kunisch/Wagner 13a}$\,$, p.~959, Theorem 2.8.$\,]\,$. In view of Theorem 1.4.~above, the assertion remains true for the bidomain system with the Rogers-McCulloch,
FitzHugh-Nagumo or the linearized Aliev-Panfilov model.
\medskip
4) $[\,${\kpt Kunisch/Wagner 13a}$\,$, pp.~961 ff., Proof of Theorem 2.7.$\,]\,$. As mentioned in the proofs of Theorem 1.1.~and 1.2., $[\,${\kpt Kunisch/Wagner 13a}$\,$, p.~959 f., Lemma 2.9.$\,]$ and $[\,${\kpt Kunisch/Wagner 13a}$\,$, p.~961, Lemma 2.10.$\,]$ hold true (the latter is identical with Lemma 3.1.~above). Steps 1$\,-\,$6 of the proof must be changed along the lines of the proof of Theorem 1.2. The basic error entered with the inequalities $[\,${\kpt Kunisch/Wagner 13a}$\,$, p.~964, (2.69) and (2.70)$\,]\,$, which do not hold true. Steps 7 and 8 can be maintained without changes.
\medskip
5) $[\,${\kpt Kunisch/Wagner 13a}$\,$, p.~973 f., Remark 1)$\,]\,$. Remains true without changes.
\medskip
6) $[\,${\kpt Kunisch/Wagner 13a}$\,$, p.~973 f., Remark 2), (2.127)$\,-\,$(2.132)$\,]\,$. This derivation has to be changed along the lines of Part C) of the proof of Theorem 1.1.~above.
\medskip
7) $[\,${\kpt Kunisch/Wagner 13a}$\,$, p.~975 f., Theorem 3.2.$\,]\,$. This assertion remains true.
\medskip
8) $[\,${\kpt Kunisch/Wagner 13a}$\,$, p.~976 f., Theorem 3.3., Corollary 3.4.~and 3.5.$\,]\,$. These assertions remain true while its proof must be corrected.
\medskip
9) $[\,${\kpt Kunisch/Wagner 13a}$\,$, p.~978 f., Proposition 3.6.~and 3.7.$\,]\,$. The propositions remain true.
\medskip
10) $[\,${\kpt Kunisch/Wagner 13a}$\,$, pp.~979 ff., Proof of Theorem 3.3.$\,]\,$. The derivation in Step 1 holds wrong since $[\,${\kpt Nagaiah/Kunisch/Plank 11}$\,]\,$, p.~158, (33), cannot be applied. Instead, the claimed $\bochnerl {4} {(\,0\,,\,T\,)} {\L {6} {} \home} $-regularity of the transmembrane potential $\pottr$ within a given weak solution $(\pottr, \potex, \gat)$ of the bidomain system may be derived from Theorem 1.5., 2), which is applicable by $[\,${\kpt Kunisch/Wagner 13a}, p.~975, Assumption 3.1., 4)$\,]\,$, thus still guaranteeing $\ioncurr (\pottr, \gat) \in \bochnerl {4/3} {(\,0\,,\,T\,)} {\L {2} {} \home} $.
\bigbreak
{\bf c) Corrections within $[\,${\kpt Kunisch/Wagner 13b}$\,]\,$.}
\medskip
1) $[\,${\kpt Kunisch/Wagner 13b}$\,$, p.~1082, Theorem 2.3.$\,]\,$. The existence and uniqueness theorem for the bidomain system must be replaced by Theorem 1.4.~above.
\medskip
2) $[\,${\kpt Kunisch/Wagner 13b}$\,$, p.~1082, Theorem 2.4.$\,]\,$. The stability estimate must be replaced by Theorem 1.2.~and Corollary 3.3.~above.
\medskip
3) $[\,${\kpt Kunisch/Wagner 13b}$\,$, p.~1083 f.$\,]\,$. The analytical framework for the analysis of the optimal control problem $(\P)$ must be chosen in accordance with the regularity requirements within the stability estimate and the adjoint existence theorem. In view of Theorem 1.2.~above, we may work with $W_0 \in \L {4} {} \home$ even in the case of the linearized Aliev-Panfilov model.
\medskip
4) $[\,${\kpt Kunisch/Wagner 13b}$\,$, p.~1084 f., Proposition 3.1.$\,-\,$3.3., Theorem 3.4.$\,]\,$. The assertions remain true without changes.
\medskip
5) $[\,${\kpt Kunisch/Wagner 13b}$\,$, p.~1087, Theorem 4.1.$\,]\,$. While Part 1) of the theorem remains true, Part 2) must be changed considerably. The theorem should be replaced by the following
\medskip
{\bf Theorem 4.2.} {\it Consider the optimal control problem $\,(\P)$ given through $[\,${\rm{\kpt Kunisch/Wagner 13b}, p.~1083 f., (3.11)$\,-\,$(3.17)}$\,]$ under the assumptions from Subsection 3.1.~there, and specify the Rogers-McCulloch model. Assume further that the integrand $r(x,t,\varphi,\eta,w)$ is continuously differentiable with respect to $\varphi$, $\eta$ and $w$.
\smallskip
1) {\bf (A-priori estimates for weak solutions of the adjoint system in the bidomain case)} If $\,(\pottrhat, \potexhat$, $\gathat, \excexhat)$ is a feasible solution of $\,(\P)$ with
$$\eqalignno{
{} & {\partial r \over \partial \varphi} (\pottrhat, \potexhat, \gathat)\,,\,\, {\partial r \over \partial \eta} (\pottrhat, \potexhat, \gathat) \,,\, {\partial r \over \partial w} (\pottrhat, \potexhat, \gathat) \in \L {2} {} (\Omega_T)\, & (4.24)
\cr}$$
then every weak solution $(P_1, P_2, P_3) \in \bochnerl {2} {(\,0\,,\,T\,)} {\sob {1,2} {} \home} \,\times \, \bochnerl {2} {(\,0\,,\,T\,)} {\sob {1,2} {} \home} \,$ $ \times \,\L {2} {} (\Omega_T)$ of the related adjoint system $[\,${\rm{\kpt Kunisch/Wagner 13b}, p.~1086, (4.8)$\,-\,$(4.10)}$\,]$ obeys the estimate
$$\eqalignno{
{} & \hbox{\kern-20pt} {\norm{P_1}}^2_{\bochnerl {\infty} {(\,0\,,\,T\,)} {\L {2} {} \home} }  + {\norm{ P_1 }}^2_{\bochnerl {2} {(\,0\,,\,T\,)} {\sob {1,2} {} \home} } + \, {\norm{P_2}}^2_{\bochnerl {2} {(\,0\,,\,T\,)} {\sob {1,2} {} \home} }  +  {\norm{P_3}}^2_{\bochnerl {\infty} {(\,0\,,\,T\,)} {\L {2} {} \home} }
\cr
{} &\le \, C\,\Bigl({\norm{{\partial r\over\partial \varphi } (\pottrhat, \potexhat, \gathat)}} {}^2_{\L {2} {} (\Omega_T) } + \,{\norm{{\partial r\over\partial \eta } (\pottrhat, \potexhat, \gathat)}} {}^2_{\L {2} {} (\Omega_T) } + {\norm{ {\partial r\over \partial w} (\pottrhat, \potexhat, \gathat) }} {}^2_{\L {2} {} (\Omega_T) }\,\Bigr) & (4.25)
\cr}$$
where the constant $C > 0$ does not depend on $P_1$, $P_2$, $P_3$ but on $(\pottrhat, \potexhat, \gathat, \excexhat)$ and the data of $\,(\P)$.
\smallskip
2) {\bf (Gain of regularity for the variables $P_1$ and $P_3$)} If $\,(\pottrhat, \potexhat, \gathat, \excexhat)$ is a feasible solution of $\,(\P)$ such that, besides of $(4.24)$,  the additional regularity
$$\eqalignno{
{} & {\partial r \over \partial \varphi} (\pottrhat, \potexhat, \gathat)\,,\,\, {\partial r \over \partial \eta} (\pottrhat, \potexhat, \gathat) \in \bochnerl {4} {(\,0\,,\,T\,)} {\L {2} {} \home} \, & (4.26)
\cr}$$
is guaranteed then every weak solution $(P_1, P_2, P_3) \in \bochnerl {2} {(\,0\,,\,T\,)} {\sob {1,2} {} \home} \,\times \, \bochnerl {2} {(\,0\,,\,T\,)} {\sob {1,2} {} \home} \,$\break $ \times \,\L {2} {} (\Omega_T)$ of the related adjoint system $[\,${\rm{\kpt Kunisch/Wagner 13b}, p.~1086, (4.8)$\,-\,$(4.10)}$\,]$ obeys the further estimate
$$\eqalignno{
{} & \hbox{\kern-20pt} {\norm{P_1}}^2_{\bochnerc {0} {[\,0\,,\,T\,]} {\L {2} {} \home} }  + {\norm{\partial P_1 / \partial s}}^{4/3}_{\bochnerl {4/3} {(\,0\,,\,T\,)} {\bigl(\, \sob {1,2} {} \home \,\bigr)^*} } + \, {\norm{P_3}}^2_{\bochnerc {0} {[\,0\,,\,T\,]} {\L {2} {} \home} } + {\norm{\partial P_3 / \partial s }}^2_ {\L {2} {} (\Omega_T)}
\cr
{} & \le\,\Bigl(\, 1 + {\norm{{\partial r\over\partial \varphi } (\pottrhat, \potexhat, \gathat)}} {}^{4}_{\bochnerl {4} {(\,0\,,\,T\,)} {\L {2} {} \home } } + {\norm{{\partial r\over\partial \eta } (\pottrhat, \potexhat, \gathat)}} {}^{4}_{\bochnerl {4} {(\,0\,,\,T\,)} {\L {2} {} \home } } & (4.27)
\cr
{} & & \displaystyle + \,{\norm{ {\partial r\over \partial w} (\pottrhat, \potexhat, \gathat) }} {}^2_{\L {2} {} (\Omega_T) }\,\Bigr)
\cr}$$
where the constant $C > 0$ does not depend on $P_1$, $P_2$, $P_3$ but on $(\pottrhat, \potexhat, \gathat, \excexhat)$ and the data of $\,(\P)$.}%
\medskip
6) $[\,${\kpt Kunisch/Wagner 13b}$\,$, p.~1088, Theorem 4.2.$\,]\,$. The theorem should be replaced by the following
\medskip
{\bf Theorem 4.3.~(Existence and uniqueness of weak solutions for the adjoint system in the bidomain case)} {\it Under the assumptions of Theorem 4.2., 2) above, the adjoint system $[\,${\rm{\kpt Kunisch/Wagner 13b}, p.~1086, (4.8)$\,-\,$(4.10)}$\,]$ admits a uniquely determined weak solution $(P_1, P_2, P_3)$ with
$$\eqalignno{
{} & P_1 \,\in \, \bochnerc {0} {[\,0\,,\,T\,]} {\L {2} {} \home} \, \cap \,\bochnerl {4} {(\,0\,,\,T\,)} {\sob {1,2} {} \home} \,\cap\, \bochnersob {1,4/3} {(\,0\,,\,T\,)} {\bigl(\,\sob {1,2} {} \home \,\bigr)^*} \,; & (4.28)
\cr
{} & P_2 \,\in\, \bochnerl {2} {(\,0\,,\,T\,)} {\sob {1,2} {} \home} \,;\quad \int_\Omega P_2(x,t)\,dx = 0\,\,\,\fforall\,t \in (\,0\,,\,T\,)\,; &  (4.29)
\cr
{} & P_3 \,\in \,\bochnerc {0} {[\,0\,,\,T\,]} {\L {2} {} \home} \,\cap\, \bochnersob {1,2} {(\,0\,,\,T\,)} { \L {2} {} \home } \,. & (4.30)
\cr}$$}%
\par
Consequently, the remark after $[\,${\kpt Kunisch/Wagner 13b}$\,$, p.~1088, (4.25)$\,]$ may be dropped. Note that the assumption (4.26) about the higher regularity of $\partial r (\pottrhat, \potexhat, \gathat) / \partial \eta$ is rather restrictive since, in ge\-ne\-ral, only $\potexhat \in \L {2} {} \bigl[\,{(\,0\,,\,T\,)}\,,$ $ {\sob {1,2} {} \home} \,\bigr]$ can be guaranteed.$\,$\footnote{$^{34)}$}{An example of a functional satisfying (4.26) is given by $\displaystyle {1 \over 2} \int_\Omega \Bigl( \, {1 \over 2\,\e} \int^{t + \e}_{t - \e} \potex(x,\tau)\, d\tau - \Phi_{\hbox{\sevenit desired}} (x)\,\Bigr)^2\, dx$.} In fact, the bidomain system does not allow for smoothing of $\potex$ or $\potin$ in time but only of $\pottr$. However, (4.26) is indispensable in order to ensure the $\bochnerl {4} {(\,0\,,\,T\,)} {\sob {1,2} {} \home} $-regularity of the multiplier $P_1$, which cannot be avoided in the pairing $[\,${\kpt Kunisch/Wagner 13b}$\,$, p.~1085, (4.1)$\,]\,$.
\medskip
7) $[\,${\kpt Kunisch/Wagner 13b}$\,$, pp.~1088 ff., Proof of Theorem 4.1.$\,]\,$. The proof with the necessary corrections is repeated here as
\medskip
{\bf Proof of Theorem 4.2.} In analogy to the primal bidomain equations, the weak adjoint system $[\,${\rm{\kpt Kunisch/ Wagner 13b}, p.~1086, (4.8)$\,-\,$(4.10)}$\,]$ can be equivalently rewritten as a reduced system in terms of the bidomain bilinear form ${\cal A}$, which is defined as in Subsection 2.e) above:
$$\eqalignno{
{} & {d \over ds} \langle\, P_1(s)\,,\, \psi\,\rangle + {\cal A}\bigl(\, P_1(s)\,,\, \psi\,\bigr) + \int_\Omega \Bigl(\, {\partial \ioncurr \over \partial \varphi} (\pottrhat, \gathat)\, P_1 + {\partial G \over \partial \varphi} (\pottrhat, \gathat)\, P_3\,\Bigr)\,\psi\, dx \,=\, \langle \, \widetilde S(s)\,,\,\psi\,\rangle & (4.31)
\cr
{} & & \displaystyle  \quad \forall\, \psi \in \sob {1,2} {} \home\,;
\cr
{} & {d \over ds} \langle \, P_3 (s)\,,\, \psi\,\rangle + \int_\Omega \Bigl(\, {\partial \ioncurr \over \partial w} (\pottrhat, \gathat)\, P_1 + {\partial G \over \partial w} (\pottrhat, \gathat)\, P_3\,\Bigr)\,\psi\, dx \,=\, - \langle \, {\partial r\over \partial w} (\pottrhat, \potexhat, \gathat)\,,\, \psi\,\rangle & (4.32)
\cr
{} & & \displaystyle  \quad \forall\, \psi \in \L {2} {} \home\,;
\cr
{} & P_1(x,0) \,= \,0\quad \fforall\, x \in \Omega\,; \quad P_3 (x,0) \,=\, 0  \quad \fforall\, x \in \Omega\, \vphantom{\int} & (4.33)
\cr}$$
on $[\,0\,,\,T\,]$ in distributional sense, cf.~$[\,${\kpt Kunisch/Wagner 13a}, p.~956 f., Theorem 2.4.$\,]\,$.
\medskip
{\it Part 1)\/} The proof is identical with $[\,${\kpt Kunisch/Wagner 13b}, pp.~1088 ff., Steps 1$\,-\,$3 and Step 6$\,]$ of the proof of Theorem 4.1., 1) there.
\par
In the following, we will particularly employ $[\,${\kpt Kunisch/Wagner 13b}, p.~1088, Lemma 4.3.$\,]\,$, and the norm estimate $[\,${\kpt Kunisch/Wagner 13b}, p.~1090, (4.51)$\,]\,$, which is repeated here:
$$\eqalignno{
{} & {\norm{P_1}}^2_{\bochnerl {\infty} {(\,0\,,\,T\,)} {\L {2} {} \home} } +  {\norm{P_3}}^2_{\bochnerl {\infty} {(\,0\,,\,T\,)} {\L {2} {} \home} }\,\le \, C\,\Bigl(\, {\norm{{\partial r\over\partial \varphi } (\pottrhat, \potexhat, \gathat)}} {}^2_{\L {2} {} (\Omega_T) } & (4.34)
\cr
{} & & \displaystyle + \,{\norm{{\partial r\over\partial \eta } (\pottrhat, \potexhat, \gathat)}} {}^2_{\L {2} {} (\Omega_T) } + {\norm{ {\partial r\over \partial w} (\pottrhat, \potexhat, \gathat) }} {}^2_{\L {2} {} (\Omega_T) }\,\Bigr)\,.
\cr}$$
\par
{\it Part 2)\/} Assume now that the additional regularity assumptions (4.26) hold. Note that, by Proposition 2.6., 2), ${\norm{\gathat(s)}}_{\L {4} {} \home}$ is uniformly bounded. Consequently, the same ist true for ${\norm{\gathat(s)}}_{\L {3} {} \home}$.
\smallskip
$\bullet$ {\bf Step 1.} {\it An estimate for $\bigl\vert \, {1 \over 2} \,{d \over ds}\, {\norm{P_1 (s)}}^2_{\L {2} {} \home} \,\bigr\vert$.} Inserting into equation (4.31) the feasible test function $\psi = P_1(s) \in \sob {1,2} {} \home$ and applying the lower estimate for the bidomain form, we obtain
$$\eqalignno{
{} & \hbox{\kern-20pt} \langle\,{d \over ds} \,P_1(s) \,,\, P_1(s)\,\rangle + {\cal A} \bigl(\, P_1(s) \,,\, P_1(s) \,\bigr) \,=\, - \int_\Omega  \Bigl(\, {\partial \ioncurr \over \partial \varphi} (\pottrhat, \gathat)\, P_1(s) + {\partial G \over \partial \varphi} (\pottrhat, \gathat)\, P_3(s)\,\Bigr)\,\,P_1(s) \, dx & (4.35)
\cr
{} & & \displaystyle +\, \langle \, \widetilde S(s)\,,\,P_1(s) \,\rangle \quad \follows
\cr
{} & \hbox{\kern-20pt}{1 \over 2} \,{d \over ds}\, {\norm{P_1}}^2_{\L {2} {} \home} + \beta \,{\norm{P_1}}^2_{\sob {1,2} {} \home} \,\le \, -  \int_\Omega \Bigl(\, 3 \,b\, (\pottrhat)^2\, P_1(s)^2 - 2\,(a+1)\, b\, \pottrhat\, P_1(s)^2 + a\,b\,P_1(s)^2 & (4.36)
\cr
{} & & \displaystyle + \,\gathat \,P_1(s)^2 - \e\,\kappa\, P_1(s)\,P_3(s)\,\Bigr)\, dx
\, +  \,\bigl\vert \,\langle \, \widetilde S(s)\,,\,P_1(s) \,\rangle  \,\bigr\vert  + \beta \,{\norm{P_1}}^2_{\L {2} {} \home} \quad \follows
\cr
{} & \hbox{\kern-20pt}{1 \over 2} \,{d \over ds}\, {\norm{P_1}}^2_{\L {2} {} \home} + \beta \,{\norm{P_1}}^2_{\sob {1,2} {} \home} \, \le \, C\,\int_\Omega \bigl\vert \,P_1\,\bigr\vert^2\,\bigl(\,\bigl\vert \,\pottrhat\,\bigr\vert + \bigl\vert\,\gathat\,\bigr\vert\,\bigr)\, dx & (4.37)
\cr
{} & & \displaystyle + \,C \int_\Omega \,\bigl\vert \,P_1\, P_3\,\bigr\vert \, dx  +  \,\bigl\vert \,\langle \, \widetilde S(s)\,,\,P_1(s) \,\rangle  \,\bigr\vert  + \beta \,{\norm{P_1}}^2_{\L {2} {} \home}\,.
\cr}$$
To the first term on the right-hand side, we apply the generalized Cauchy inequality and subsequently H\"older's inequality, thus obtaining
$$\eqalignno{
{} &  C\, \int_\Omega \bigl\vert \,P_1\,\bigr\vert^2\,\Bigl(\, \bigl\vert \,\pottrhat\,\bigr\vert + \bigl\vert \,\gathat\,\bigr\vert \,\Bigr) \, dx \, \le\,C\,\e_1(s)\, \int_\Omega \bigl\vert \,P_1\,\bigr\vert^2\,\bigl(\,\bigl\vert \,\pottrhat\,\bigr\vert^2 + \bigl\vert\,\gathat\,\bigr\vert^2\,\bigr)\, dx + {C \over \e_1(s)} \,{\norm {P_1}}^2_{\L {2} {} \home}  & (4.38)
\cr
{} &  \le\, C\,\e_1(s)\, {\norm{P_1}}^2_{\L {4} {} \home} \, {\norm{\pottrhat}}^2_{\L {4} {} \home} + C\,\e_1(s)\, {\norm{P_1}}^2_{\L {6} {} \home} {\norm{\gathat}}^2_{\L {3} {} \home} + {C \over \e_1(s)} \, {\norm {P_1}}^2_{\L {2} {} \home}  & (4.39)
\cr}$$
for arbitrary $\e_1(s) > 0$. Specifying within (4.39) $\e_1(s) = \e^\p_1 / \bigl(\,1 + {\norm{\pottrhat}}^2_{\L {4} {} \home} + {\norm{\gathat}}^2_{\L {3} {} \home} \,\bigr)$ and noticing the (almost) uniform bound (4.34) for ${\norm {P_1 (s) }}^2_{\L {2} {} \home}$, we arrive at
$$\eqalignno{
{} & C\, \int_\Omega \bigl\vert \,P_1 \,\bigr\vert^2\,\Bigl(\, \bigl\vert \,\pottrhat\,\bigr\vert + \bigl\vert \,\gathat\,\bigr\vert \,\Bigr) \, dx \, \le \,C\,\e^\p_1  \,{\norm{P_1}}^2_{\sob {1,2} {} \home} + {C \over \e^\p_1 } \,\Bigl(\,1 + {\norm{\pottrhat}}^2_{\L {4} {} \home} + {\norm{\gathat}}^2_{\L {3} {} \home}\,\Bigr)  \,. & (4.40)
\cr}$$
Further, we observe that
$$\eqalignno{
{} & C \int_\Omega \,\bigl\vert \,P_1\, P_3\,\bigr\vert \, dx \,\le\, C\, \Bigl(\,{\norm{P_1}}^2_{\L {2} {} \home} + {\norm{P_3}}^2_{\L {2} {} \home}\,\Bigr)  & (4.41)
\cr}$$
is (almost) uniformly bounded again by (4.34), $[\,${\kpt Kunisch/Wagner 13b}, p.~1088, Lemma 4.3.$\,]$ yields the estimate
$$\eqalignno{
{} & \bigl\vert\,\langle\, \widetilde S(s) \,,\,P_1\,\rangle\,\bigr\vert \,\le\, C\,\e^\p_0\, {\norm{P_1}}^2 _{\sob {1,2} {} \home} + {C \over \e^\p_0}\,\Bigl(\, {\norm{{\partial r\over\partial \varphi } (\pottrhat, \potexhat, \gathat)}} {}^2_{\L {2} {} \home} +  {\norm{{\partial r\over\partial \eta } (\pottrhat, \potexhat, \gathat)}} {}^2_{\L {2} {} \home} \,\Bigr)  & (4.42)
\cr}$$
for arbitrary $\e^\p_0 > 0$, and the last summand on the right-hand side of (4.37) is (almost) uniformly bounded again by (4.34). Consequently, (4.37), (4.40), (4.41) and (4.42) imply together
$$\eqalignno{
{} & \hbox{\kern-20pt}{1 \over 2} \,{d \over ds}\, {\norm{P_1}}^2_{\L {2} {} \home} + \beta \,{\norm{P_1}}^2_{\sob {1,2} {} \home} \, \le\,C\,\bigl(\,\e^\p_0 + \e^\p_1\,\bigr)\, {\norm{P_1}}^2 _{\sob {1,2} {} \home} & (4.43)
\cr
{} & & \displaystyle + \,{C \over \e^\p_0}\,\Bigl(\, {\norm{{\partial r\over\partial \varphi } (\pottrhat, \potexhat, \gathat)}} {}^2_{\L {2} {} \home} +  {\norm{{\partial r\over\partial \eta } (\pottrhat, \potexhat, \gathat)}} {}^2_{\L {2} {} \home} \,\Bigr)  + {C \over \e^\p_1 } \,\Bigl(\,1 + {\norm{\pottrhat}}^2_{\L {4} {} \home} + {\norm{\gathat}}^2_{\L {3} {} \home}\,\Bigr) \,.
\cr}$$
Choosing now $\e^\p_0$, $\e^\p_1 > 0$ in such a way that the terms with ${\norm{P_1}}^2_{\sob {1,2} {} \home}$ on both sides of (4.43) annihilate, we obtain the inequality
$$\eqalignno{
{} & \hbox{\kern-20pt}{1 \over 2} \,{d \over ds}\, {\norm{P_1}}^2_{\L {2} {} \home} \, \le \,{C \over \e^\p_0}\,\Bigl(\, {\norm{{\partial r\over\partial \varphi } (\pottrhat, \potexhat, \gathat)}} {}^2_{\L {2} {} \home} +  {\norm{{\partial r\over\partial \eta } (\pottrhat, \potexhat, \gathat)}} {}^2_{\L {2} {} \home} \,\Bigr)  & (4.44)
\cr
{} & & \displaystyle + \,{C \over \e^\p_1 } \,\Bigl(\,1 + {\norm{\pottrhat}}^2_{\L {4} {} \home} + {\norm{\gathat}}^2_{\L {3} {} \home}\,\Bigr) \,.
\cr}$$
Inserting now into (4.31) the reverse test function $\psi = - P_1(s)$, we get instead
$$\eqalignno{
{} & \hbox{\kern-20pt} - \langle\,{d \over ds} \,P_1(s) \,,\, P_1(s),\rangle - {\cal A} \bigl(\, P_1(s)  \,,\, P_1(s) \,\bigr) \,=\, \int_\Omega \Bigl(\, {\partial \ioncurr \over \partial \varphi} (\pottrhat, \gathat)\, P_1(s) & (4.45)
\cr
{} & & \displaystyle + \,{\partial G \over \partial \varphi} (\pottrhat, \gathat)\, P_3(s)\,\Bigr)\,\,P_1(s) \, dx -\, \langle \, \widetilde S(s)\,,\,P_1(s) \,\rangle \quad \follows
\cr
{} & \hbox{\kern-20pt}- {1 \over 2} \,{d \over ds}\, {\norm{P_1}}^2_{\L {2} {} \home} -\beta \,{\norm{P_1}}^2_{\sob {1,2} {} \home} \,\ge \, \int_\Omega \Bigl(\, 3 \,b\, (\pottrhat)^2\, P_1(s)^2 - 2\,(a+1)\, b\, \pottrhat\, P_1(s)^2 + a\,b\,P_1(s)^2 & (4.46)
\cr
{} & & \displaystyle + \,\gathat \,P_1(s)^2 - \e\,\kappa\, P_1(s)\,P_3(s)\,\Bigr)\, dx
\, -  \,\bigl\vert \,\langle \, \widetilde S(s)\,,\,P_1(s) \,\rangle  \,\bigr\vert  -\beta \,{\norm{P_1}}^2_{\L {2} {} \home} \quad \follows
\cr
{} & \hbox{\kern-20pt}- {1 \over 2} \,{d \over ds}\, {\norm{P_1}}^2_{\L {2} {} \home} - \beta \,{\norm{P_1}}^2_{\sob {1,2} {} \home} \, \ge \, - C\,\int_\Omega \bigl\vert \,P_1\,\bigr\vert^2\,\bigl(\,\bigl\vert \,\pottrhat\,\bigr\vert + \bigl\vert\,\gathat\,\bigr\vert\,\bigr)\, dx & (4.47)
\cr
{} & & \displaystyle - \,C \int_\Omega \,\bigl\vert \,P_1\, P_3\,\bigr\vert \, dx -
\,\bigl\vert \,\langle \, \widetilde S(s)\,,\,P_1(s) \,\rangle  \,\bigr\vert  - \beta \,{\norm{P_1}}^2_{\L {2} {} \home}\,.
\cr}$$
Using again (4.15) and (4.40)$\,-\,$(4.42), (4.47) implies the reverse inequality
$$\eqalignno{
{} & \hbox{\kern-20pt}- {1 \over 2} \,{d \over ds}\, {\norm{P_1}}^2_{\L {2} {} \home} - \beta \,{\norm{P_1}}^2_{\sob {1,2} {} \home} \, \le\, -C\,\bigl(\,\e^\p_0 + \e^\p_1\,\bigr)\, {\norm{P_1}}^2 _{\sob {1,2} {} \home} & (4.48)
\cr
{} & & \displaystyle - \,{C \over \e^\p_0}\,\Bigl(\, {\norm{{\partial r\over\partial \varphi } (\pottrhat, \potexhat, \gathat)}} {}^2_{\L {2} {} \home} +  {\norm{{\partial r\over\partial \eta } (\pottrhat, \potexhat, \gathat)}} {}^2_{\L {2} {} \home} \,\Bigr) - {C \over \e^\p_1 } \,\Bigl(\,1 + {\norm{\pottrhat}}^2_{\L {4} {} \home} + {\norm{\gathat}}^2_{\L {3} {} \home}\,\Bigr) \,,
\cr}$$
and by an appropriate choice of $\e^\p_0$, $\e^\p_1 > 0$, the summands with $- {\norm{P_1}}^2_{\sob {1,2} {} \home}$ annihilate. Thus (4.44) is reversed as
$$\eqalignno{
{} & \hbox{\kern-20pt}- {1 \over 2} \,{d \over ds}\, {\norm{P_1}}^2_{\L {2} {} \home} \, \ge \,- {C \over \e^\p_0}\,\Bigl(\, {\norm{{\partial r\over\partial \varphi } (\pottrhat, \potexhat, \gathat)}} {}^2_{\L {2} {} \home} + {\norm{{\partial r\over\partial \eta } (\pottrhat, \potexhat, \gathat)}} {}^2_{\L {2} {} \home} \,\Bigr) & (4.49)
\cr
{} & & \displaystyle - \,{C \over \e^\p_1 } \,\Bigl(\,1 + {\norm{\pottrhat}}^2_{\L {4} {} \home} + {\norm{\gathat}}^2_{\L {3} {} \home}\,\Bigr) \,,
\cr}$$
and we arrive at the desired estimate
$$\eqalignno{
{} & \Bigl\vert\,{1 \over 2} \,{d \over ds}\, {\norm{P_1(s))}}^2_{\L {2} {} \home} \,\Bigr\vert\,\le \, {C \over \e^\p_0}\,\Bigl(\, {\norm{{\partial r\over\partial \varphi } (\pottrhat, \potexhat, \gathat)}} {}^2_{\L {2} {} \home} + {\norm{{\partial r\over\partial \eta } (\pottrhat, \potexhat, \gathat)}} {}^2_{\L {2} {} \home} \,\Bigr) & (4.50)
\cr
{} & & \displaystyle + \,{C \over \e^\p_1 } \,\Bigl(\,1 + {\norm{\pottrhat}}^2_{\L {4} {} \home} + {\norm{\gathat}}^2_{\L {3} {} \home}\,\Bigr) \,. \cr}$$
\par
$\bullet$ {\bf Step 2.} {\it An estimate for ${\norm{P_1}}^4_{\bochnerl {4} {(\,0\,,\,T\,)} {\sob {1,2} {} \home} }$.} We return to (4.43) and choose now $\e^\p_0$, $\e^\p_0 > 0$ in such a way that $C\,\e^\p_0 + C\,\e^\p_1 = \beta/2$. Then the inequality (4.43) and (4.50) imply
$$\eqalignno{
{} & {\beta\over 2} \,{\norm{P_1}}^2_{\sob {1,2} {} \home} \, \le\,- {1 \over 2} \,{d \over ds}\, {\norm{P_1}}^2_{\L {2} {} \home} + \,{C \over \e^\p_0}\,\Bigl(\, {\norm{{\partial r\over\partial \varphi } (\pottrhat, \potexhat, \gathat)}} {}^2_{\L {2} {} \home} +  {\norm{{\partial r\over\partial \eta } (\pottrhat, \potexhat, \gathat)}} {}^2_{\L {2} {} \home} \,\Bigr)  & (4.51)
\cr
{} & & \displaystyle +\, {C \over \e^\p_1 } \,\Bigl(\,1 + {\norm{\pottrhat}}^2_{\L {4} {} \home} + {\norm{\gathat}}^2_{\L {3} {} \home}\,\Bigr)
\cr
{} & \le\, \Bigl\vert\, {1 \over 2} \,{d \over ds}\, {\norm{P_1}}^2_{\L {2} {} \home}\,\Bigr\vert + \,{C \over \e^\p_0}\,\Bigl(\, {\norm{{\partial r\over\partial \varphi } (\pottrhat, \potexhat, \gathat)}} {}^2_{\L {2} {} \home} +  {\norm{{\partial r\over\partial \eta } (\pottrhat, \potexhat, \gathat)}} {}^2_{\L {2} {} \home} \,\Bigr)  & (4.52)
\cr
{} & & \displaystyle +\, {C \over \e^\p_1 } \,\Bigl(\,1 + {\norm{\pottrhat}}^2_{\L {4} {} \home} + {\norm{\gathat}}^2_{\L {3} {} \home}\,\Bigr)
\cr
{} & \le\, 2 \,{C \over \e^\p_0}\, \Bigl(\, {\norm{{\partial r\over\partial \varphi } (\pottrhat, \potexhat, \gathat)}} {}^2_{\L {2} {} \home} +  {\norm{{\partial r\over\partial \eta } (\pottrhat, \potexhat, \gathat)}} {}^2_{\L {2} {} \home} \,\Bigr) + 2\, {C \over \e^\p_1 } \,\Bigl(\, 1 + {\norm{\pottrhat}}^2_{\L {4} {} \home} + {\norm{\gathat}}^2_{\L {3} {} \home}\,\Bigr)\,.
\cr}$$
Consequently, we find \hfill (4.53)
$$\eqalignno{
{} & {\norm{P_1}}^4_{\bochnerl {4} {(\,0\,,\,T\,)} {\sob {1,2} {} \home} } \,=\, \int^T_0 {\norm{P_1 (s)}}^4_{\sob {1,2} {} \home} \, ds & (4.54)
\cr
{} & \le \,C \int^T_0 \! \Bigl(\,1 + {\norm{\pottrhat}}^4_{\L {4} {} \home} + {\norm{\gathat}}^4_{\L {3} {} \home} + {\norm{{\partial r\over\partial \varphi } (\pottrhat, \potexhat, \gathat)}} {}^4_{\L {2} {} \home} +  {\norm{{\partial r\over\partial \eta } (\pottrhat, \potexhat, \gathat)}} {}^4_{\L {2} {} \home} \,\Bigr)\,ds & (4.55)
\cr
{} & \le\,C \,\Bigl(\,1 + {\norm{\pottrhat}}^4_{\L {4} {} (\Omega_T) } + {\norm{\gathat}}^4_{\bochnerc {0} {[\,0\,,\,T\,]} {\L {3} {} \home} } & (4.56)
\cr
{} & & \displaystyle + \,{\norm{{\partial r\over\partial \varphi } (\pottrhat, \potexhat, \gathat)}} {}^4 _{\bochnerl {4} {(\,0\,,\,T\,)} {\L {2} {} \home} }  + {\norm{{\partial r\over\partial \eta } (\pottrhat, \potexhat, \gathat)}} {}^4_{\bochnerl {4} {(\,0\,,\,T\,)} {\L {2} {} \home} }  \,\Bigr)\, ,
\cr}$$
and the right-hand side is bounded by (2.13), Proposition 2.6., 2) and the additional regularity assumptions about $\partial r (\pottrhat, \potexhat, \gathat) / \partial \varphi$ and $\partial r  (\pottrhat, \potexhat, \gathat) / \partial \eta$.
\medskip
$\bullet$ {\bf Step 3.} {\it The estimate for ${\norm{\partial P_1 / \partial s}}^{4/3}_{\bochnerl {4/3} {(\,0\,,\,T\,)} {\bigl(\, \sob {1,2} {} \home \,\bigr)^*} }$.} Exploiting the definition of the dual norm, we start with
$$\eqalignno{
{} & \hbox{\kern-20pt} {\norm{\partial P_1 / \partial s}}_{\bigl(\, \sob {1,2} {} \home \,\bigr)^*} \, =\, \sup_{ {\norm{\psi}}_{\sob {1,2} {} \home}\, = \, 1} \bigl\vert \,\langle \, \partial P_1 (s)/ \partial s\,,\, \psi\,\rangle\,\bigr\vert & (4.57)
\cr
{} & \le \, \sup_{\lldots} \,\bigl\vert \, {\cal A} (P_1, \psi)\,\bigr\vert +  \sup_{\lldots}  \int_\Omega \,\Bigl\vert\, {\partial \ioncurr \over \partial \varphi} (\pottrhat, \gathat)\, P_1 + {\partial G \over \partial \varphi} (\pottrhat, \gathat)\, P_3\,\Bigr\vert\,\bigl\vert\,\psi\,\bigr\vert\, dx +  \sup_{\lldots} \,\bigl\vert \,\langle \, \widetilde S(s)\,,\,\psi\,\rangle \,\bigr\vert & (4.58)
\cr
{} & \le \, \sup_{\lldots} \,\bigl\vert \, {\cal A} (P_1, \psi)\,\bigr\vert +  \sup_{\lldots}  \,C \int_\Omega \Bigl(\,1 + \bigl\vert \,\pottrhat \,\bigr\vert + \bigl\vert \,\gathat\,\bigr\vert + \bigl\vert\,\pottrhat\,\bigr\vert^2 \,\Bigr) \,\bigl\vert \, P_1\,\bigr\vert\,\bigl\vert \,\psi\,\bigr\vert \, dx & (4.59)
\cr
{} & & \displaystyle +\, \sup_{\lldots}\, C \int_\Omega \bigl\vert \, P_3\,\bigr\vert \,\bigl\vert \,\psi\,\bigr\vert \, dx +  \sup_{\lldots} \,\bigl\vert \,\langle \, \widetilde S(s)\,,\,\psi\,\rangle \,\bigr\vert \,.
\cr}$$
The four terms on the right-hand side of (4.59) will be estimated separately. For the first term, we get with the upper estimate for the bidomain form$\,$\footnote{$^{35)}$}{Cf.~$[\,${\ninekpt Kunisch/Wagner 13a}, p.~957, Theorem 2.4., 2), (2.27)$\,]\,$.}
$$\eqalignno{
{} & \sup_{\lldots} \,\bigl\vert \, {\cal A} (P_1, \psi)\,\bigr\vert \, \le\, \sup_{\lldots} \,\gamma \, {\norm{P_1(s)}}_{\sob {1,2} {} \home}\,{\norm{\psi}}_{\sob {1,2} {} \home}\, \le \, \gamma \, {\norm{P_1}}_{\sob {1,2} {} \home} \, . & (4.60)
\cr}$$
For the second term, we obtain
$$\eqalignno{
{} & \sup_{\lldots} \,C \int_\Omega \Bigl(\,1 + \bigl\vert \,\pottrhat \,\bigr\vert + \bigl\vert \,\gathat\,\bigr\vert + \bigl\vert\,\pottrhat\,\bigr\vert^2 \,\Bigr) \,\bigl\vert \, P_1\,\bigr\vert\,\bigl\vert \,\psi\,\bigr\vert \, dx \, & (4.61)
\cr
{} & \le \,  C + \sup_{\lldots} \, C  \int_\Omega \bigl\vert\, \pottrhat\, P_1\,\psi\,\bigr\vert \, dx + \sup_{\lldots} \, C  \int_\Omega \bigl\vert\, \gathat \, P_1\,\psi\,\big\vert\, dx  + \sup_{\lldots} \, C  \int_\Omega \bigl\vert\, \pottrhat\,\bigr\vert^2 \cdot \bigl\vert \, P_1\,\psi\,\bigr\vert \, dx
\cr}$$
$$\eqalignno{
{} & \le\, C + \sup_{\lldots} \, C \,\Bigl(\, \int_\Omega \bigl\vert\, \pottrhat\, \bigr\vert^4\, dx\, \Bigr)^{1/4} \, \Bigl(\, \int_\Omega \bigl\vert\, P_1\,\bigr\vert^4 \, dx\,\Bigr)^{1/4} \, \Bigl(\, \int_\Omega \bigl\vert \, \psi\,\bigr\vert^2\, dx \,\Bigr)^{1/2}  + \sup_{\lldots} \, C \,\Bigl(\, \int_\Omega \bigl\vert\, \gathat\,\bigr\vert^3\, dx\,\bigr)^{1/3} \cdot & (4.62)
\cr
{} & & \displaystyle \Bigl(\,\int_\Omega \bigl\vert\, P_1\,\bigr\vert^3\, dx\,\Bigr)^{1/3} \,\Bigl(\, \int_\Omega \bigl\vert \,\psi \, \bigr\vert^3 \, dx \,\Bigr)^{1/3}
+ \sup_{\lldots} \,C\, \Bigl(\, \int_\Omega \bigl\vert \,\pottrhat\,\bigr\vert^4\,dx\,\Bigr)^{1/2} \,\Bigl(\, \int_\Omega \bigl\vert \, P_1\,\bigr\vert^4\,dx\,\Bigr)^{1/4} \,\Bigl(\, \int_\Omega \bigl\vert \,\psi\,\bigr\vert^4\,dx\,\Bigr)^{1/4}
\cr
{} & \le\, C + \sup_{\lldots} \, C \cdot {\norm{\pottrhat}}_{\L {4} {} \home} \cdot {\norm{P_1}}_{\L {4} {} \home} \cdot {\norm{\psi}}_{\sob {1,2} {} \home} \phantom{\int} & (4.63)
\cr
{} & & \displaystyle + \,\sup_{\lldots} \, C \cdot {\norm{\gathat}}_{\L {3} {} \home} \cdot {\norm{P_1}}_{\L {3} {} \home} \cdot {\norm{\psi}}_{\sob {1,2} {} \home}
+ \sup_{\lldots} \,C\cdot {\norm{\pottrhat}}^2_{\L {4} {} \home} \cdot {\norm{P_1}}_{\L {4} {} \home} \cdot {\norm{\psi}}_{\sob {1,2} {} \home}
\cr
{} & \le \, C\, \Bigl(\, 1 \, + \,{\norm{\pottrhat}}_{\L {4} {} \home} + {\norm{\gathat}}_{\L {3} {} \home} + {\norm{\pottrhat}}^2_{\L {4} {} \home} \,\Bigr) \cdot {\norm{P_1}}_{\L {4} {} \home}  \,.& (4.64)
\cr}$$
For the third summand, we get
$$\eqalignno{
{} &  \sup_{\lldots}\, C \int_\Omega \bigl\vert \, P_3\,\bigr\vert \,\bigl\vert \,\psi\,\bigr\vert \, dx \,\le \,  \sup_{\lldots}\, C \,\Bigl(\, \int_\Omega \bigl\vert \,P_3\,\bigr\vert^2\,dx\,\Bigr)^{1/2} \,\Bigl(\, \int_\Omega \bigl\vert \,\psi\,\bigr\vert^2\,dx\,\Bigr)^{1/2} \, & (4.65)
\cr
{} & & \displaystyle \le \sup_{\lldots}\,  C\cdot {\norm{P_3}}_{\L {2} {} \home} \cdot {\norm{\psi}}_{\L {2} {} \home} \,\le\, C \,{\norm{P_3}}_{\L {2} {} \home} \,, \quad (4.66)
\cr}$$
and for the last summand, $[\,${\kpt Kunisch/Wagner 13b}, p.~1088, Lemma 4.3.$\,]$ (with $\e^\p_0 = 1$) yields the estimate
$$\eqalignno{
{} & \hbox{\kern-20pt} \sup_{\lldots}\, \bigl\vert \,\langle \, \widetilde S(s)\,,\,\psi\,\rangle \,\bigr\vert \,\le\, \sup_{\lldots}\,  C\,{\norm{\psi}}^2 _{\sob {1,2} {} \home} +  C \, \Bigl(\, {\norm{{\partial r\over\partial \varphi } (\pottrhat, \potexhat, \gathat)}} {}^2_{\L {2} {} \home} +  {\norm{{\partial r\over\partial \eta } (\pottrhat, \potexhat, \gathat)}} {}^2_{\L {2} {} \home} \,\Bigr) & (4.67)
\cr
{} & \le \, C\, \Bigl(\,1 + {\norm{{\partial r\over\partial \varphi } (\pottrhat, \potexhat, \gathat)}} {}^2_{\L {2} {} \home} +  {\norm{{\partial r\over\partial \eta } (\pottrhat, \potexhat, \gathat)}} {}^2_{\L {2} {} \home} \,\Bigr) \,. & (4.68)
\cr}$$
Inserting (4.60), (4.64), (4.66) and (4.68) into (4.57), we arrive at
$$\eqalignno{
{} & \hbox{\kern-20pt} {\norm{\partial P_1 / \partial s}}_{\bigl(\, \sob {1,2} {} \home \,\bigr)^*} \, \le\, C\, \Bigl(\ 1 + \bigl(\,1 + {\norm{\pottrhat}}_{\L {4} {} \home}  + {\norm{\gathat}}_{\L {3} {} \home} + {\norm{\pottrhat}}^2_{\L {4} {} \home} \,\bigr) \cdot {\norm{P_1}}_{\sob {1,2} {} \home}\, & (4.69)
\cr
{} & + \,{\norm{P_3}}_{\L {2} {} \home} + {\norm{{\partial r\over\partial \varphi } (\pottrhat, \potexhat, \gathat)}} {}^2_{\L {2} {} \home} +  {\norm{{\partial r\over\partial \eta } (\pottrhat, \potexhat, \gathat)}} {}^2_{\L {2} {} \home} \,\Bigr) \,.
\cr}$$
Consequently, we find the estimate
$$\eqalignno{
{} & \hbox{\kern-20pt}{\norm{\partial P_1 / \partial s}}^{4/3}_{\bochnerl {4/3} {(\,0\,,\,T\,)} {\bigl(\, \sob {1,2} {} \home \,\bigr)^*} } \,=\, \int^T_0 {\norm{\partial P_1 (s) / \partial s}}^{4/3}_{\bigl(\, \sob {1,2} {} \home \,\bigr)^*} \, ds & (4.70)
\cr
{} & \le\,C +  C\, \int^T_0 \Bigl(\, 1 + {\norm{\pottrhat}}^{4/3}_{\L {4} {} \home} + {\norm{\gathat}}^{4/3}_{\L {3} {} \home} + {\norm{\pottrhat}}^{8/3}_{\L {4} {} \home}\,\Bigr) \cdot {\norm{P_1}}^{4/3}_{\sob {1,2} {} \home}\, ds  & (4.71)
\cr
{} & & \displaystyle + \, C \int^T_0 {\norm{P_3}}^{4/3}_{\L {2} {} \home} \, ds  +  C\,\int^T_0 \Bigl(\, {\norm{{\partial r\over\partial \varphi } (\pottrhat, \potexhat, \gathat)}} {}^{8/3}_{\L {2} {} \home} +  {\norm{{\partial r\over\partial \eta } (\pottrhat, \potexhat, \gathat)}} {}^{8/3}_{\L {2} {} \home} \,\Bigr) \, ds
\cr
{} & \hbox{\kern-20pt} \le \, C\,\Bigl(\, \int^T_0 \Bigl(\, 1 + {\norm{\pottrhat}}^{4/3}_{\L {4} {} \home} + {\norm{\gathat}}^{4/3}_{\L {3} {} \home} + {\norm{\pottrhat}}^{8/3}_{\L {4} {} \home}\,\Bigr)^{3/2} \,ds\,\Bigr)^{2/3} \, \Bigl(\, \int^T_0 {\norm{P_1}}^{4}_{\sob {1,2} {} \home}\, ds \,\Bigr)^{1/3} + \lldots & (4.72)
\cr
{} & \hbox{\kern-20pt} \le \, C\,\Bigl(\, 1 + {\norm{\pottrhat}}^{4/3}_{\bochnerl {2} {(\,0\,,\,T\,)} {\L {4} {} \home} } + {\norm{\gathat}}^{4/3}_{\bochnerl {2} {(\,0\,\,T\,)} {\L {3} {} \home} } + {\norm{\pottrhat}}^{8/3}_{\L {4} {} (\Omega_T)}\,\Bigr) \cdot {\norm{P_1}}^{4/3}_{\bochnerl {4} {(\,0\,,\,T\,)} {\sob {1,2} {} \home} } & (4.73)
\cr
{} & + C\,\Bigl( \,1 +  {\norm{P_3}}^{4/3}_{\bochnerl {4/3} {(\,0\,,\,T\,)} {\L {2} {} \home} } + {\norm{{\partial r\over\partial \varphi } (\pottrhat, \potexhat, \gathat)}} {}^{8/3}_{\bochnerl {8/3} {(\,0\,,\,T\,)} {\L {2} {} \home} }
\cr
{} & & \displaystyle + \,{\norm{{\partial r\over\partial \eta } (\pottrhat, \potexhat, \gathat)}} {}^{8/3}_{\bochnerl {8/3} {(\,0\,,\,T\,)} {\L {2} {} \home} }\,\Bigr)\,.
\cr}$$
The first summand is bounded by (2.13), (4.56) and Proposition 2.6., 2), the second one is bounded by (4.34) and the additional regularity assumptions about $\partial r (\pottrhat, \potexhat, \gathat) / \partial \varphi$ and $\partial r (\pottrhat, \potexhat, \gathat) / \partial \eta$. Summing up, we obtain the claimed estimate
$$\eqalignno{
{} & {\norm{\partial P_1 / \partial s}}^{4/3}_{\bochnerl {4/3} {(\,0\,,\,T\,)} {\bigl(\, \sob {1,2} {} \home \,\bigr)^*} } \,
\cr
{} & \hbox{\kern-20pt} \le\, C\,\Bigl( \, 1 + \,{\norm{{\partial r\over\partial \varphi } (\pottrhat, \potexhat, \gathat)}} {}^{4/3} _{\bochnerl {4} {(\,0\,,\,T\,)} {\L {2} {} \home} }  + \, {\norm{{\partial r\over\partial \eta } (\pottrhat, \potexhat, \gathat)}} {}^{4/3}_{\bochnerl {4} {(\,0\,,\,T\,)} {\L {2} {} \home} } & (4.74)
\cr
{} & + \, {\norm{{\partial r\over\partial \varphi } (\pottrhat, \potexhat, \gathat)}} {}^{4/3}_{\L {2} {} (\Omega_T) } + \,{\norm{{\partial r\over\partial \eta } (\pottrhat, \potexhat, \gathat)}} {}^{4/3}_{\L {2} {} (\Omega_T) } + {\norm{ {\partial r\over \partial w} (\pottrhat, \potexhat, \gathat) }} {}^{4/3}_{\L {2} {} (\Omega_T) }
\cr
{} &  + \,{\norm{{\partial r\over\partial \varphi } (\pottrhat, \potexhat, \gathat)}} {}^{8/3}_{\bochnerl {8/3} {(\,0\,,\,T\,)} {\L {2} {} \home} } + {\norm{{\partial r\over\partial \eta } (\pottrhat, \potexhat, \gathat)}} {}^{8/3}_{\bochnerl {8/3} {(\,0\,,\,T\,)} {\L {2} {} \home} }\,\Bigr)
\cr
{} & \hbox{\kern-20pt} \le \, C\,\Bigl( \, 1 + \,{\norm{{\partial r\over\partial \varphi } (\pottrhat, \potexhat, \gathat)}} {}^{8/3} _{\bochnerl {4} {(\,0\,,\,T\,)} {\L {2} {} \home} }  + \, {\norm{{\partial r\over\partial \eta } (\pottrhat, \potexhat, \gathat)}} {}^{8/3}_{\bochnerl {4} {(\,0\,,\,T\,)} {\L {2} {} \home} } & (4.75)
\cr
{} & & \displaystyle + \,{\norm{ {\partial r\over \partial w} (\pottrhat, \potexhat, \gathat) }} {}^{4/3}_{\L {2} {} (\Omega_T) } \,\Bigr)\,.
\cr}$$
\par
$\bullet$ {\bf Step 4.} {\it The estimate for $ {\norm{\partial P_3 / \partial s}}^2_ {\bochnerl {2} {(\,0\,,\,T\,)} {\bigl(\,\L {2} {} \home\,\bigr)^* } } = {\norm{\partial P_3 / \partial s}}^2_ {\L {2} {} (\Omega_T)}$.} We start again with the calculation of the dual norm
$$\eqalignno{
{} & \hbox{\kern-20pt} {\norm{\partial P_3 / \partial s}}_{ \bigl(\, \L {2} {} \home \,\bigr)^*}  \,=\, \sup_{ {\norm{\psi}}_{\L {2} {} \home}\, = \, 1} \bigl\vert \,\langle \, \partial P_3 (s)/ \partial s\,,\, \psi\,\rangle\,\bigr\vert  & (4.76)
\cr
{} & =\, \sup_{\lldots}\, \int_\Omega \,\Bigl\vert\,{\partial \ioncurr \over \partial w} (\pottrhat, \gathat)\, P_1 + {\partial G \over \partial w} (\pottrhat, \gathat)\, P_3 + {\partial r\over \partial w} (\pottrhat, \potexhat, \gathat)\,\Bigr\vert \cdot \bigl\vert\,\psi\,\bigr\vert\, dx & (4.77)
\cr
{} & \le\, \sup_{\lldots}\, \int_\Omega \bigl\vert \, \pottrhat\, P_1\, \psi\,\bigr\vert \, dx + \sup_{\lldots} \,\int_\Omega \e\,\bigl\vert \,P_3\,\psi\,\bigr\vert \, dx + \sup_{\lldots} \, \int_\Omega \,\Bigl\vert\, {\partial r\over \partial w} (\pottrhat, \potexhat, \gathat)\,\Bigr\vert \cdot \bigl\vert\,\psi\,\bigr\vert\, dx \,. & (4.78)
\cr}$$
Relying on (4.61)$\,-\,$(4.66), we obtain
$$\eqalignno{
{} & \hbox{\kern-20pt} {\norm{\partial P_3 (s) / \partial s}}_{\bigl(\, \L {2} {} \home \,\bigr)^*}  \,
\cr
{} & \le\, \sup_{\lldots} \, C \,\Bigl(\, {\norm{\pottrhat}}_{\L {4} {} \home} \cdot {\norm{P_1}}_{\L {4} {} \home} + {\norm{P_3}}_{\L {2} {} \home}+  \norm{{\partial r\over \partial w} (\pottrhat, \potexhat, \gathat)} {}_{\L {2} {} \home} \,\Bigr) \cdot {\norm{\psi}}_{\L {2} {} \home} & (4.79)
\cr
{} & \le \, C \, \Bigl(\, {\norm{\pottrhat}}_{\L {4} {} \home} \cdot {\norm{P_1}}_{\L {4} {} \home} + {\norm{P_3}}_{\L {2} {} \home} + \norm{{\partial r\over \partial w} (\pottrhat, \potexhat, \gathat)} {}_{\L {2} {} \home} \,\Bigr)\,. & (4.80)
\cr}$$
Together with (4.34) and (4.56), we arrive at the claimed estimate
$$\eqalignno{
{} & \hbox{\kern-20pt}{\norm{\partial P_3 / \partial s}}^2_{\L {2} {} (\Omega_T) }\,=\, \int^T_0  {\norm{\partial P_3 (s) / \partial s}}_{\bigl(\, \L {2} {} \home \,\bigr)^*} \, ds & (4.81)
\cr
{} & \le \, C \int^T_0 \Bigl(\, {\norm{\pottrhat}}^2_{\L {4} {} \home} \cdot {\norm{P_1}}^2_{\L {4} {} \home} + {\norm{P_3}}^2_{\L {2} {} \home} + \norm{{\partial r\over \partial w} (\pottrhat, \potexhat, \gathat)} {}^2_{\L {2} {} \home} \,\Bigr) \, ds & (4.82)
\cr
{} & \le\, C\,\Bigl(\, {\norm{\pottrhat}}^2_{\L {4} {} (\Omega_T)} \cdot {\norm{P_1}}^2_{\bochnerl {4} {(\,0\,,\,T\,)} {\sob {1,2} {} \home} } + {\norm{P_3}}^2_{\L {2} {} (\Omega_T)} + \norm{{\partial r\over \partial w} (\pottrhat, \potexhat, \gathat)} {}^2_{\L {2} {} (\Omega_T)} \,\Bigr)  & (4.83)
\cr
{} & \le\, C\,\Bigl(\, {\norm{{\partial r\over\partial \varphi } (\pottrhat, \potexhat, \gathat)}} {}^4 _{\bochnerl {4} {(\,0\,,\,T\,)} {\L {2} {} \home} }  + {\norm{{\partial r\over\partial \eta } (\pottrhat, \potexhat, \gathat)}} {}^4_{\bochnerl {4} {(\,0\,,\,T\,)} {\L {2} {} \home} }  & (4.84)
\cr
{} & & \displaystyle  + \,{\norm{{\partial r\over\partial \varphi } (\pottrhat, \potexhat, \gathat)}} {}^2_{\L {2} {} (\Omega_T) } +  {\norm{{\partial r\over\partial \eta } (\pottrhat, \potexhat, \gathat)}} {}^2_{\L {2} {} (\Omega_T) } + 2\, {\norm{ {\partial r\over \partial w} (\pottrhat, \potexhat, \gathat) }} {}^2_{\L {2} {} (\Omega_T) }\,\Bigr)
\cr
{} & \le \, C\,\Bigl(\, {\norm{{\partial r\over\partial \varphi } (\pottrhat, \potexhat, \gathat)}} {}^4 _{\bochnerl {4} {(\,0\,,\,T\,)} {\L {2} {} \home} }  + {\norm{{\partial r\over\partial \eta } (\pottrhat, \potexhat, \gathat)}} {}^4_{\bochnerl {4} {(\,0\,,\,T\,)} {\L {2} {} \home} } & (4.85)
\cr
{} & & \displaystyle + \,{\norm{ {\partial r\over \partial w} (\pottrhat, \potexhat, \gathat) }} {}^2_{\L {2} {} (\Omega_T) }\,\Bigr) \,.
\cr}$$
\par
$\bullet$ {\bf Step 5.} {\it Conclusion of the proof.} The fact that $P_1$ belongs even to $\bochnerc {0} {[\,0\,,\,T\,]} {\L {2} {} \home} $ can be confirmed analogously to $[\,${\kpt Bourgault/Coudi\`ere/Pierre 09}$\,]\,$, p.~478, Subsection 5.3. As a consequence of the imbedding theorem $[\,${\kpt Evans 98}$ \,]\,$, p.~286, Theorem 2, $P_3 \in \bochnerc {0} {[\,0\,,\,T\,]} {\L {2} {} \home} $ holds true as well. Consequently, the norms on the left-hand side of (4.34) can be replaced by $\bochnerc {0} {[\,0\,,\,T\,]} {\L {2} {} \home} $-norms, and the proof is complete. $\blackbox$
\medskip
8) $[\,${\kpt Kunisch/Wagner 13b}$\,$, pp.~1096 ff., Proof of Theorem 4.2.$\,]\,$. The proof with the necessary corrections is sketched here as
\medskip
{\bf Proof of Theorem 4.3.} Note first that $[\,${\kpt Kunisch/Wagner 13b}$\,$, p.~1097, Lemma 4.4.$\,]$ must be replaced by the following
\medskip
{\bf Lemma 4.4.} {\it Let the assumptions of Theorem 4.2., 2) above hold for the data of $\,(\P)$ and a feasible solution $(\pottrhat, \potexhat,$ $ \gathat, \excexhat)$. Then for all $N \in {\hbox{\blackboard N}}_0$, the functions $P^N_1$, $P^N_3$ satisfy the estimate
$$\eqalignno{
{} & \hbox{\kern-20pt} {\norm{P^N_1}}^2_{\bochnerc {0} {[\,0\,,\,T\,]} {\L {2} {} \home} }  + {\norm{ P^N_1 }}^4_{\bochnerl {4} {(\,0\,,\,T\,)} {\sob {1,2} {} \home} }   + {\norm{\partial P^N_1 / \partial s}}^{4/3}_{\bochnerl {473} {(\,0\,,\,T\,)} {\bigl(\, \sob {1,2} {} \home \,\bigr)^*} } & (4.86)
\cr
{} & & \displaystyle  + \, {\norm{P^N_3}}^2_{\bochnerc {0} {[\,0\,,\,T\,]} {\L {2} {} \home} } + {\norm{\partial P^N_3 / \partial s }}^2_{\bochnerl {2} {(\,0\,,\,T\,)} {\L {2} {} \home} } \, \le\,C
\cr}$$
for a constant $C > 0$ independent of $N$.}
\medskip
Then $[\,${\kpt Kunisch/Wagner 13b}$\,$, pp.~1096 ff., Proof of Theorem 4.2.$\,]$ runs as before with the appropriate modifications, caused by the usage of Theorem 4.2.~above instead of $[\,${\kpt Kunisch/Wagner 13b}$\,$, pp.~1087, Theorem 4.1.$\,]\,$. $\blackbox$
\medskip
9) $[\,${\kpt Kunisch/Wagner 13b}$\,$, pp.~1098, Remark 4.5., 1)$\,]\,$. The arguments remain true without changes. Consequently, Theorems 4.2.~and 4.3.~hold true for the FitzHugh-Nagumo model as well.
\medskip
10) $[\,${\kpt Kunisch/Wagner 13b}$\,$, pp.~1098 f, Remark 4.5., 2)$\,]\,$. In the case of the linearized Aliev-Panfilov model, $W_0 \in \L {4} {} \home$ implies $\gathat \in \C {0} {} \bigl[ \,{[\,0\,,\,T\,]}\,,$ $ {\L {3} {} \home}\,\bigr] $ by Proposition 2.7.~above. Then the derivations $[\,${\kpt Kunisch/Wagner 13b}$\,$, p.~1098 f., (4.132)$\,-$ (4.136)$\,]$ can be maintained but $[\,${\kpt Kunisch/Wagner 13b}$\,$, p.~1099, (4.137)$\,-\,$(4.140)$\,]$ must be changed in the following way.
\par
In Step 3 of the proof of Theorem 4.2., the following summand must be added to the right-hand side of (4.59):
$$\eqalignno{
{} & \hbox{\kern-20pt} \sup_{\lldots } \,\int_\Omega \bigl\vert \,\pottrhat \,P_3\,\psi\,\bigr\vert \, dx  \,\le\, \sup_{\lldots } \, C\cdot {\norm{\pottrhat}}_{\L {4} {} \home} \cdot {\norm{P_3}}_{\L {2} {} \home} \cdot {\norm{\psi}}_{\L {4} {} \home} \,\le \,C\cdot {\norm{\pottrhat}}_{\L {4} {} \home} \cdot {\norm{P_3}}_{\L {2} {} \home} & (4.90)
\cr
{} & \follows \, \,\, C\,\int^T_0 {\norm{\pottrhat}}^{4/3} _{\L {4} {} \home} \cdot {\norm{P_3}}^{4/3}_{\L {2} {} \home} \,ds \,\le \,C\, {\norm{\pottrhat}}^{4/3} _{\L {4} {} (\Omega_T) } \cdot {\norm{P_3}}^{4/3}_{\L {2} {} (\Omega_T)} \,. & (4.91)
\cr}$$
Consequently, Step 3 can be maintained in the case of the linearized Aliev-Panfilov model as well. In Steps 4 and 5, no changes have to be made. We conclude that Theorems 4.2.~and 4.3.~remain true even if the linearized Aliev-Panfilov model is specified.
\medskip
11) $[\,${\kpt Kunisch/Wagner 13b}$\,$, pp.~1100, Theorem 5.2.$\,]\,$. The assumptions of this theorem can be substantially weakened. In fact, Theorems 4.2.~and 4.3.~allow to settle the optimality conditions within a framework of {\it weak solutions of the bidomain system\/} instead of strong ones. In particular, the conditions apply now even to optimal controls $\excexhat \in \bochnerl {\infty} {(\,0\,,\,T\,)} {\L {2} {} \home } \setminus \bochnersob {1,2} {(\,0\,,\,T\,)} {\L {2} {} \home} $, which cause no additional time smoothing of the state. Thus the theorem should be replaced by the following
\medskip
{\bf Theorem 4.5.~(First-order necessary optimality conditions for the control problem $(\P)$ in the bidomain case)} {\it We consider the optimal control problem $\,(\P)$ given through $[\,${\rm{\kpt Kunisch/Wagner 13b}, p.~1083 f., (3.11)$\,-\,$(3.17)}$\,]$ under the assumptions of Subsection 3.1.~there with either the Rogers-McCulloch, the FitzHugh-Nagumo or the linearized Aliev-Panfilov model. Assume further that the integrand $r(x,t,\varphi, \eta, w)$ is continuously differentiable with respect to $\varphi$, $\eta$ and $w$. Let $\,(\pottrhat, \potexhat, \gathat , \excexhat)$ be a weak local minimizer of $\,(\P)$ such that
$$\eqalignno{
{} & {\partial r \over \partial \varphi} (\pottrhat, \potexhat, \gathat)\,,\,\, {\partial r \over \partial \eta} (\pottrhat, \potexhat, \gathat) \in \bochnerl {4} {(\,0\,,\,T\,)} {\L {2} {} \home} \,;\,\, {\partial r \over \partial w} (\pottrhat, \potexhat, \gathat)  \in \L {2} {} (\Omega_T) \,.& (4.92)
\cr}$$
Then there exist multipliers $P_1 \,\in \, \bochnerl {4} {(\,0\,,\,T\,)} {\sob {1,2} {} \home} $, $P_2 \,\in\, \bochnerl {2} {(\,0\,,\,T\,)} {\sob {1,2} {} \home} \,\cap \, \bigl\{\, Z \,\,\big\vert\, \int_\Omega\, Z (x,t)\, dx$ $=\, 0 \,\,\, \fforall\, t \in (\,0\,,\,T\,)\,\bigr\}$ and $P_3 \,\in\, \L {2} {} (\Omega_T)$, satisfying together with $(\pottrhat, \potexhat, \gathat , \excexhat)$ the optimality condition
$$\eqalignno{
{} & \int^T_0 \int_{\Omega_{\hbox{\sevenit con}}} \Bigl(\,\mu\, \excexhat - Q\, P_2 \,\Bigr)\cdot \bigl(\, \excex - \excexhat\,\bigr) \, dx \, dt \,\ge \, 0 \quad \forall \, \excex \in {\cal C}  & (4.93)
\cr}$$
where the control domain ${\cal C}  \subset \bochnerl {\infty} {(\,0\,,\,T\,)} {\L {2} {} \home} $ and the operator $Q$ are described through $[\,${\rm{\kpt Kunisch/ Wagner 13b}, p.~1078, (1.9)$\,-\,$(1.10)}$\,]\,$, as well as the adjoint equations
$$\eqalignno{
{} & \hbox{\kern-20pt} \int^T_0 \int_\Omega \Bigl(\, - {\partial P_1 \over \partial t} + {\partial \ioncurr\over \partial \varphi} (\pottrhat, \gathat)\,P_1 + {\partial G \over \partial \varphi} (\pottrhat, \gathat)\, P_3\,\Bigr)\,\psi\,dx\, dt \,+\, \int^T_0 \int_\Omega \nabla \psi^{\T} \,\condin\,\bigl(\, \nabla P_1 + \nabla P_2\,\bigr) \, dx\, dt
\cr
{} & =\, - \int^T_0 \!\!\int_\Omega \Bigl(\, {\partial r\over\partial \varphi } (\pottrhat, \potexhat, \gathat)  \,\Bigr)\, \psi\, dx\, dt \quad \forall\,\psi \in \bochnerl {2} {(\,0\,,\,T\,)} {\sob {1,2} {} \home} \,, \,\, P_1(x, T) \equiv 0\,; & (4.94)
\cr}$$
$$\eqalignno{
{} & \hbox{\kern-20pt} \int^T_0 \int_\Omega \nabla \psi^{\T} \,\condin\,\nabla P_1 \,dx \, dt \,+ \,\int^T_0 \int_\Omega \nabla \psi^{\T} \,(\condin+ \condex) \,\nabla P_2 \,dx \, dt\, =\, - \int^T_0 \int_\Omega {\partial r\over\partial \eta } (\pottrhat, \potexhat, \gathat) \,\psi\, dx\, dt & (4.95)
\cr
{} & & \displaystyle \quad \forall\,\psi \in \bochnerl {2} {(\,0\,,\,T\,)} {\sob {1,2} {} \home} \,\,\,\hbox{with}\,\,\, \int_\Omega \psi(x,t)\, dx = 0 \quad \fforall\, t \in (\,0\,,\,T\,) \,, \,\, \int_ \Omega P_2(x, t)\, dx = 0\quad \fforall \, t \in (\,0\,,\, T\,)\,;
\cr
{} & \hbox{\kern-20pt} \int^T_0 \!\int_\Omega \Bigl(\, - {\partial P_3 \over \partial t} + {\partial \ioncurr \over \partial w} (\pottrhat, \gathat)\, P_1 + {\partial G \over \partial w} (\pottrhat, \gathat)\, P_3\,\Bigr)\, \psi\, dx \,dt \, =\, - \int^T_0 \!\int_\Omega \Bigl(\, {\partial r\over \partial w} (\pottrhat, \potexhat, \gathat) \,\Bigr)\, \psi\, dx\, dt& (4.96)
\cr
{} &  &\displaystyle  \quad \forall\, \psi \in \bochnerl {2} {(\,0\,,\,T\,)} {\L {2} {} \home} \,,\,\, P_3(x, T ) \equiv 0\,,
\cr}$$
which are solved in weak sense. The multipliers admit the additional regularity described in $(4.28\,-\,(4.30)$.}%
\medskip
{\bf Proof.} $[\,${\kpt Kunisch/Wagner 13b}, pp.~1101 ff., Proof of Theorem 5.2.$\,]$ may be repeated with obvious modifications caused by the usage of Theorems 4.2.~and 4.3.~above instead of $[\,${\kpt Kunisch/Wagner 13b}, p.~1087 f., Theorem 4.1.~and 4.2.$\,]\,$. Note that, even in the case of the linearized Aliev-Panfilov model, the assumptions of Theorem 4.2.~and 4.3.~remain unchanged. $\blackbox$
\medskip
12) $[\,${\kpt Kunisch/Wagner 13b}$\,$, pp.~1100 f., Corollary 5.3.$\,-\,$5.5.$\,]\,$. These assertions hold true even under the assumptions of Theorem 4.5.~above.
\bigskip
{\bf d) Corrections within $[\,${\kpt Kunisch/Nagaiah/Wagner 11}$\,]\,$.}
\medskip
1) $[\,${\kpt Kunisch/Nagaiah/Wagner 11}, p.~260, Theorem 2.2.$\,]\,$. In view of Theorems 1.1.~and 1.3.~above, the assertion remains true without changes.
\bigskip\medskip
{\kapitel References.}
\bigskip
\parindent=0pt \def\ref #1 {{$[\,${\kpt #1}$\,]$\ }}
{\ninepoint \font\blackboard=dsrom10 scaled 900 \font\bf=cmbx9
\def\cyr{\ninecyr\cyracc}\def\cyrit{\ninecyrit\cyracc}
\baselineskip=12pt
\hangindent=17pt\hangafter=1
$\phantom{0}$1.\kern5pt\ref {Adams/Fournier 07} Adams, R.~A.; Fournier, J.~J.~F.: {\it Sobolev Spaces.} Academic Press / Elsevier; Amsterdam etc.~2007, 2nd ed.
\smallskip\hangindent=17pt\hangafter=1
$\phantom{0}$2.\kern5pt\ref {Aliev/Panfilov 96} Aliev, R.~R.; Panfilov, A.~V.: {\it A simple two-variable model of cardiac excitation.} Chaos, Solitons \&\ Fractals {\bf 7} (1996), 293 -- 301
\smallskip\hangindent=17pt\hangafter=1
$\phantom{0}$3.\kern5pt\ref {Barrett/S\"uli 12} Barrett, J.~W.; S\"uli, E.: {\it Reflections on Dubinski\u\i{}'s nonlinear compact embedding theorem.} Publ.~Inst.~Math.~(Beograd) (N.S.) {\bf 91} (2012), 95 -- 110
\smallskip\hangindent=17pt\hangafter=1
$\phantom{0}$4.\kern5pt\ref {Bourgault/Coudi\`ere/Pierre 09} Bourgault, Y.; Coudi\`ere, Y.; Pierre, C.: {\it Existence and uniqueness of the solution for the bidomain model used in cardiac electrophysiology.} Nonlinear Analysis: Real World Appl.~{\bf 10} (2009), 458 -- 482
\smallskip\hangindent=17pt\hangafter=1
$\phantom{0}$5.\kern5pt\ref {Dubinskij 65} Dubinskij, Ju.~A.~{\cyr (Dubinski\u\i, Yu.~A.):} {\cyrit Slabaya skhodimost\cprime\ v neline\u\i nykh \`ellipticheskikh i parabolicheskikh uravneniyakh.} {\cyr Mat.~Sbornik} {\bf 67} (1965) 4, 609 -- 642
\smallskip\hangindent=17pt\hangafter=1
$\phantom{0}$6.\kern5pt\ref {Evans 98} Evans, L.~C.: {\it Partial Differential Equations.}
American Mathematical Society; Providence 1998
\smallskip\hangindent=17pt\hangafter=1
$\phantom{0}$7.\kern5pt\ref {FitzHugh 61} FitzHugh, R.: {\it Impulses and physiological states in theoretical models of nerve membrane.} Biophysical J.~{\bf 1} (1961), 445 -- 466
\smallskip\hangindent=17pt\hangafter=1
$\phantom{0}$8.\kern5pt\ref {Fursikov 00} Fursikov, A.~V.: {\it Optimal Control of Distributed Systems. Theory and Applications.} American Mathematical Society; Providence 2000 (Translations of Mathematical Monographs 187)
\smallskip\hangindent=17pt\hangafter=1
$\phantom{0}$9.\kern5pt\ref {Kunisch/Nagaiah/Wagner 11} Kunisch, K.; Nagaiah, C.; Wagner, M.: {\it A parallel Newton-Krylov method for optimal control of the monodomain model in cardiac electrophysiology.} Computing and Visualization in Science {\bf 14} (2011), 257 -- 269
\smallskip\hangindent=17pt\hangafter=1
10.\kern5pt\ref {Kunisch/Wagner 12} Kunisch, K.; Wagner, M.: {\it Optimal control of the bidomain system (I): The monodomain approximation with the Rogers-McCulloch model.} Nonlinear Analysis: Real World Applications  {\bf 13} (2012), 1525 -- 1550
\smallskip\hangindent=17pt\hangafter=1
11.\kern5pt\ref {Kunisch/Wagner 13a} Kunisch, K.; Wagner, M.: {\it Optimal control of the bidomain system (II): Uniqueness and regularity theorems for weak solutions.} Ann.~Mat.~Pura Appl.~{\ninebf 192} (2013), 951 -- 986
\smallskip\hangindent=17pt\hangafter=1
12.\kern5pt\ref {Kunisch/Wagner 13b} Kunisch, K.; Wagner, M.: {\it Optimal control of the bidomain system (III): Existence of minimizers and first-order optimality conditions.} ESAIM: Mathematical Modelling and Numerical Analysis {\bf 47} (2013), 1077 -- 1106
\smallskip\hangindent=17pt\hangafter=1
13.\kern5pt\ref {Lady\v zenskaja/Solonnikov/Ural'ceva 88} Lady\v zenskaja, O.~A.; Solonnikov, V.~A.; Ural'ceva, N.~N.: {\it Linear and Quasi-linear Equations of Parabolic Type.} American Mathematical Society; Providence 1988
\smallskip\hangindent=17pt\hangafter=1
14.\kern5pt\ref {Nagumo/Arimoto/Yoshizawa 62} Nagumo, J.; Arimoto, S.; Yoshizawa, S.: {\it An active pulse transmission line simulating nerve axon.} Proceedings of the Institute of Radio Engineers {\bf 50} (1962), 2061 -- 2070
\smallskip\hangindent=17pt\hangafter=1
15.\kern5pt\ref {Nagaiah/Kunisch/Plank 11} Nagaiah, C.; Kunisch, K.; Plank, G.: {\it Numerical solution for optimal control of the reaction-diffusion equations in cardiac electrophysiology.} Comput.~Optim.~Appl.~{\ninebf 49} (2011), 149 -- 178
\smallskip\hangindent=17pt\hangafter=1
16.\kern5pt\ref {Rogers/McCulloch 94} Rogers, J.~M.; McCulloch, A.~D.: {\it A collocation-Galerkin finite element model of cardiac action potential propagation.} IEEE Transactions of Biomedical Engineering {\bf 41} (1994), 743 -- 757
\smallskip\hangindent=17pt\hangafter=1
17.\kern5pt\ref {Simon 87} Simon, J.: {\it Compact sets in the space $L^p (0, T ; B)$.} Ann.~Mat.~Pura Appl.~{\bf 146} (1987), 65 -- 96
\smallskip\hangindent=17pt\hangafter=1
18.\kern5pt\ref {Stein 70} Stein, E.~M.: {\it Singular Integrals and Differentiability Properties of Functions.} Princeton University Press; Princeton 1970
\smallskip\hangindent=17pt\hangafter=1
19.\kern5pt\ref {Sundnes/Lines/Cai/Nielsen/Mardal/Tveito 06} Sundnes, J.; Lines, G.~T.; Cai, X.; Nielsen, B.~F.; Mardal, K.-A.; Tveito, A.: {\it Computing the Electrical Activity in the Heart.} Springer; Berlin etc.~2006
\smallskip\hangindent=17pt\hangafter=1
20.\kern5pt\ref {Tung 78} Tung, L.: {\it A Bi-Domain Model for Describing Ischemic Myocardial D-C Potentials.} PhD thesis. Massachusetts Institute of Technology 1978
\smallskip\hangindent=17pt\hangafter=1
21.\kern5pt\ref {Veneroni 09} Veneroni, M.: {\it Reaction-diffusion systems for the macroscopic bidomain model of the cardiac electric field.} Nonlinear Analysis: Real World Applications {\bf 10} (2009), 849 -- 868
\smallskip\hangindent=17pt\hangafter=1
\vfill
\baselineskip=13.8pt
{\ninebf Last modification:} 18.$\,$08.$\,$2015
\smallskip\noindent
{\ninebf Authors' addresses$\,$/$\,$e-mail.} {\it Karl Kunisch:} University of Graz, Institute for Mathematics and Scientific Compu\-ting, Heinrichstra{\ss}e 36, A-8010 Graz, Austria. e-mail: {\ninett karl.kunisch$\,$@$\,$uni-graz.at}
\smallskip
{\it Marcus Wagner:} Max Planck Institute for Mathematics in the Sciences, Inselstra{\ss}e 22, D-04103 Leipzig, Germany. Homepage$\,$/$\,$e-mail: {\ninett www.thecitytocome.de} / {\ninett marcus.wagner$\,$@$\,$mis.mpg.de}
\par}%
\supereject\end